\documentclass[a4paper,11pt,french,leqno]{smfart}

\usepackage{hyperref}
\usepackage{amsmath}
\usepackage{amssymb}
\usepackage{amsthm}
\usepackage{amsopn}
\usepackage{amscd}
\usepackage{fullpage}
\usepackage[utf8]{inputenc}
\usepackage{graphics, xcolor}
\usepackage{textcomp} 
\usepackage{lscape}
\usepackage{url}
\usepackage{verbatim}
\usepackage{mathrsfs}
  \usepackage[francais]{babel}

\usepackage{bm}

\def\commentaire#1{\ifx\cachercommentaires\undefined \textcolor{red}{#1}\else \fi } 
\newcommand{\scrC}{\mathscr{C}}

\newcommand{\frb}{\mathfrak{b}}

\newcommand{\frg}{\mathfrak{g}}

\newcommand{\frk}{\mathfrak{k}}
\newcommand{\frl}{\mathfrak{l}}

\newcommand{\frn}{\mathfrak{n}}
\newcommand{\fro}{\mathfrak{o}}
\newcommand{\frp}{\mathfrak{p}}
\newcommand{\frqqq}{\mathfrak{q}}

\newcommand{\frs}{\mathfrak{s}}
\newcommand{\frt}{\mathfrak{t}}
\newcommand{\fru}{\mathfrak{u}}

\newcommand{\frz}{\mathfrak{z}}

\newcommand{\bbA}{\mathbb{A}}

\newcommand{\bbC}{\mathbb{C}}

\newcommand{\bbN}{\mathbb{N}}

\newcommand{\bbQ}{\mathbb{Q}}
\newcommand{\bbR}{\mathbb{R}}

\newcommand{\bbZ}{\mathbb{Z}}

\newcommand{\caM}{\mathcal{M}}

\newcommand{\caO}{\mathcal{O}}

\newcommand{\caQ}{\mathcal{Q}}
\newcommand{\caR}{\mathcal{R}}

\newcommand{\caV}{\mathcal{V}}

\def\U{\mathbf{U}}

\newcommand{\GL}{\mathbf{GL}}

\newcommand{\Hom}{\mathrm{Hom}}

\newcommand{\Ind}{\mathrm{Ind}}

\renewcommand{\Ind}{\mathrm{Ind}}

\newcommand{\SL}{\mathbf{SL}}

\newcommand{\SO}{\mathbf{SO}}
\newcommand{\Sp}{\mathbf{Sp}}
\newcommand{\Or}{\mathbf{O}}

\newcommand{\bil}[2]{\langle  #1,#2 \rangle }

\newcommand{\sgn}{\mathbf{sgn}}
\newcommand{\Triv}{\mathbf{Triv}}

\theoremstyle{plain}
\newtheorem{thm}{Théorème}[section]
\newtheorem{lemme}[thm]{Lemme}
\newtheorem{cor}[thm]{Corollaire}
\newtheorem{prop}[thm]{Proposition}

\newtheorem*{thmI}{Théorème}

\theoremstyle{definition}
 \newtheorem{defi}[thm]{Définition}

\newtheorem{rmq}[thm]{Remarque}

\newtheorem{notation}[thm]{Notation}
\newtheorem{convention}[thm]{Convention}

\def \dem {\noindent \underline{\sl Démonstration}. }

\begin{document}

\numberwithin{equation}{section}

\title{Sur les paquets d'Arthur de $\Sp(2n,\bbR)$ contenant des modules unitaires de plus haut poids, scalaires}
  \author{Colette Moeglin}
 \address{CNRS, Institut Mathématique de Jussieu } 
 \email{colette.moeglin@imj-prg.fr}

  \author{David Renard  }
 \address{Centre de Mathématiques
 Laurent Schwartz,  Ecole Polytechnique} 
\email{david.renard@polytechnique.edu}

\date{\today}

\begin{abstract} Soit $\pi$ un module de plus haut poids unitaire du groupe $G=\Sp(2n,\bbR)$.
On s'intéresse aux paquets d'Arthur contenant $\pi$. Lorsque le plus haut poids est scalaire, on détermine 
 les paramètres de  ces paquets, on établit la  propriété de multiplicité un  de $\pi$ dans le paquet, et l'on calcule le caractère 
 $\rho_\pi$ (du groupe des composantes connexes du centralisateur du paramètre dans le groupe dual)  associé à $\pi$ et qui joue un grand rôle dans la théorie d'Arthur.
 On fait de même pour certains modules  de plus haut poids unitaires unipotents $\sigma_{n,k}$, ou bien lorsque le caractère infinitésimal est régulier.
\medskip
 
\noindent {\bf  Abstract.}\, ---\, Let $\pi$ be an irreducible unitary highest weight module for $G=\Sp(2,\bbR)$. 
We would like to determine  the Arthur packets containing $\pi$. When the highest weight is scalar, we determine
the Arthur parameter of these packets, we establish the multiplicity one property of $\pi$ in the packet and we compute
the character  $\rho_\pi$ (of the group of connected components of the centralizer of $\psi$ in the dual group) associated to $\pi$ 
which plays an important role in Arthur's theory. We also deal with the case of some unipotent unitary  highest weight modules  $\sigma_{n,k}$, or 
when the infinitesimal character is regular.
\end{abstract}

\maketitle

\section{Introduction} 

Soient $G=\Sp(2n,\bbR)$ le groupe symplectique réel de rang $n$, et $K$ un sous-groupe compact maximal de $G$;
il est isomorphe à $\U(n)$. Les modules de plus haut poids unitaires de $G$ ont été déterminés dans 
\cite {KV} (voir aussi \cite{EHW} qui montre qu'il n'y en a pas d'autres). Ils sont paramétrés par 
des $n$-uplets  décroissants d'entiers $\mu=(m_1,\ldots,m_n)$, soumis à certaines conditions (voir théorème \ref{hw}).
Notons $\pi(\mu)$ le module de plus haut poids unitaire  correspondant; son plus haut poids est donné par  $(-m_n, \ldots,-m_1)$. 
Lorsque $m_n>n$ (resp. $m_n=n$), $\pi(\mu)$ est une série discrète holomorphe  \footnote{Avec nos conventions, d'autres diraient antiholomorphe.}
(resp. une limite de séries  discrètes holomorphes). En général, ce sont des représentations unitaires très intéressantes et largement  étudiées
dans la littérature. Elles apparaissent dans de nombreux problèmes, en particulier dans la théorie des formes de Siegel, nous y reviendrons à la fin de cette introduction.

Pour décrire le spectre automorphe discret d'un groupe réductif  $\mathbf H$ sur un corps de nombres $k$, 
J. Arthur a introduit certains paramètres $\bm \psi$ et leurs localisations $\bm \psi_v$ en toute place $v$ de $k$. 
Fixons une place   $v$ et posons $\psi=\bm \psi_v$. A ce paramètre doit   correspondre un certain paquet fini   $\Pi(\psi)$
 de représentations unitaires du groupe $\mathbf H(k_v)$.
Lorsque $\mathbf H(k_v)$ est un groupe classique quasi-déployé, par exemple un groupe symplectique qui est le cas qui nous intéresse ici, 
 il détermine  dans \cite{Art13} ces paquets en les caractérisant par des identités de transfert endoscopique. 
 Notons $A(\psi)$ le groupe des composantes connexes du centralisateur de $\psi$ dans le groupe dual $\widehat H$.
 Dans le cas des groupes classiques, c'est un groupe abélien, et toujours dans ce cas,  Arthur associe à une représentation $\pi$ dans le paquet
  $\Pi(\psi)$ une représentation de dimension finie  $\rho_\pi$ de $A(\psi)$. 
  Ces représentations  $\rho_\pi$ interviennent dans les identités endoscopiques et aussi  dans les formules de multiplicité globales.
  On conjecture que ces représentations $\rho_\pi$ sont de dimension un, c'est la propriété de multiplicité un qui déjà 
  établie dans de nombreux cas.

Dans \cite{MR3} et \cite{MR5}, nous nous sommes intéressé au cas des places archimédiennes réelles, où nous avons donné des constructions
des représentations d'un paquet.  Lorsque le paramètre $\psi$ est unipotent, le paquet correspondant est étudié dans \cite{pourhowe} par des méthodes globales
avec des séries thêta, ce qui localement se traduit par des correspondances de Howe. Les paquets généraux sont obtenus
à partir des paquets unipotents en utilisant  l'induction cohomologique et l'induction parabolique.

Nos buts dans cet article  sont, premièrement,   étant donné un module de plus haut poids unitaire
$\pi(\mu)$ pour $G=\Sp(2n,\bbR)$ comme ci-dessus, de déterminer les paquets d'Arthur $\Pi(\psi)$ le contenant.
Deuxièmement, d'établir la propriété de multiplicité $1$, et troisièmement de calculer  le caractère $\rho_{\pi(\mu)}$ de $A(\psi)$.

Une condition nécessaire évidente pour que $\pi(\mu)$ appartienne à $\Pi(\psi)$ est que  le caractère infinitésimal de $\psi$ soit celui de $\pi(\mu)$, qui est entier.
Si $m_n>n$, comme nous l'avons dit, $\pi(\mu)$ est une série discrète holomorphe et son caractère infinitésimal est de plus régulier.
Les résultats de \cite{AMR} montrent alors que les paquets $\Pi(\psi)$ contenant $\pi(\mu)$ sont exactement les paquets d'Adams-Johnson
\cite{AdJo} construits par induction cohomologique,   ayant le bon caractère infinitésimal et dont la partie unipotente du paramètre est de dimension $1$.
 Dans la suite, nous supposons donc $m_n\leq n$.
Nous donnons  alors une réponse complète aux questions posées ci-dessus dans trois  cas. 
Premièrement,  lorsque le caractère infinitésimal est régulier. C'est un cas facile, les paquets  d'Arthur sont alors des paquets d'Adams-Johnson et les résultats
découlent de \cite{AMR}, et \cite{MR3}. Deuxièmement, 
lorsque le plus haut poids est scalaire, c'est-à-dire
$\mu=(m,\ldots,m)$, et dans ce cas, on note aussi la représentation  correspondante $\pi_n(m)$.   Troisièmement, lorsque 
$\mu$ est de la forme
 $\mu=(\underbrace{k+1,\ldots,k+1}_{2k}, \underbrace{k,\ldots,k}_{n-2k})$,  avec     $2\leq 2k\leq n$, 
et l'on note alors $\sigma_{n,k}$ la représentation correspondante.

Pour pouvoir énoncer nos résultats, nous avons besoin d'introduire quelques notations concernant les paramètres d'Arthur $\psi$.
Rappelons qu'un tel paramètre est un morphisme 
$ \psi: W_\bbR \times \SL_2(\bbC) \rightarrow {}^LG$
vérifiant certaines propriétés, et où $W_\bbR$ désigne le groupe de Weil de $\bbR$.
Ici, comme $G$ est déployé, on peut prendre ${}^LG=\widehat G=\SO(2n+1,\bbC)$, et en composant avec la représentation 
standard de $\SO(2n+1,\bbC)$ dans $\GL(2n+1,\bbC)$, on voit $\psi$ comme une représentation de dimension $2n+1$
de $W_\bbR \times \SL_2(\bbC)$. Cette représentation est complètement réductible. 
Pour tout $a\in \bbN^\times$, notons
$R[a]$ la représentation algébrique de $\SL_2(\bbC)$ de dimension $a$, et pour tout entier $t$ strictement positif, notons
$\delta_t$ la représentation irréductible de $W_\bbR$ de dimension $2$ qui est le paramètre de Langlands de la série
discrète de $\GL_2(\bbR)$ de caractère infinitésimal $(t/2,-t/2)$.
En tenant compte de la condition mentionnée plus haut sur le caractère infinitésimal, on peut écrire la décomposition en irréductibles de $\psi$ sous la forme :
\begin{equation}\label{formpsi}  \psi= \psi_u\oplus \psi_d=\left( \bigoplus_{i=1}^r   \eta_i \boxtimes R[a'_i] \right)
 \oplus \left( \bigoplus_{j=1}^s   \delta_{t_j} \boxtimes R[a_j]\right).
 \end{equation}
Dans la première   somme qui constitue la partie unipotente 
$\psi_u$ du paramètre,  $\eta_i$ désigne un caractère quadratique de $W_\bbR$ (le caractère trivial
que nous notons $1_{W_\bbR}$  ou  le caractère signe que nous notons $\sgn_{W_\bbR}$), les $a'_i$ sont impairs et $r=1$ ou $3$.
Dans la deuxième somme, qui constitue la partie discrète $\psi_d$ du paramètre,
  les $t_j$ sont des entiers strictement positifs  et  $t_j+a_j$ est impair. On définit  $a(\psi_u)=\max_{i=1,\ldots,r} (a'_i)$.
Notre résultat est le suivant (Théorème \ref{mainthm} du texte).

\begin{thmI} Soient $\psi$ un paramètre d'Arthur  pour $G=\Sp(2n,\bbR)$ et 
 $m$ un entier, avec  $0\leq m\leq n$. On suppose que  le caractère infinitésimal
associé à $\psi$ est  celui de $\pi_n(m)$.  Alors $\psi$ se décompose comme en (\ref{formpsi}) et 
 le paquet $\Pi(\psi)$ contient $\pi_n(m)$ si et seulement si l'on est dans un des cas suivants: 

$(i)$ $\dim \psi_u=1$, $2m >n+1$, et quels que soient $i<j$ entre $1$ et $s$, 
\begin{equation*} \left[ \frac{t_i-a_i+1}{2},  \frac{t_i+a_i-1}{2}\right]\cap \left[ \frac{t_j-a_j+1}{2},  \frac{t_j+a_j-1}{2}\right]=\emptyset , \end{equation*}

$(ii)$  $\psi$ s'écrit 
 \[ \psi=(\sgn^{\frac{2n+1-a(\psi_u)}{2}}  \boxtimes R[a(\psi_u)] )\oplus \psi' \] 
où $\psi'$ est un $A$-paramètre pour le groupe  orthogonal pair compact $\Or(0, 2n+1-a(\psi_u))$ tel que le paquet 
$\Pi(\psi')$ contienne une representation de dimension  finie $E_{\psi'}$
et $\pi_n(m)$ est image de $E_{\psi'}$ par la correspondance de Howe.

Dans ces deux cas,  la multiplicité de $\pi_n(m)$ dans le paquet $\Pi(\psi)$ est $1$.
\end{thmI}

Précisons ce que l'on veut dire par paramètre et paquet pour un groupe orthogonal pair, ce qui est assez standard. Si $\psi'$ est un paramètre
pour un  groupe spécial orthogonal pair, on considère son conjugué par l'automorphisme extérieur. On obtient un paramètre 
$\psi''$ (qui peut être équivalent ou non à $\psi$). On voit ces deux paramètres comme équivalents pour le groupe orthogonal, et le paquet associé
pour le groupe orthogonal  est constitué des représentations irréductibles dont les composantes de la restriction au groupe spécial orthogonal sont 
dans l'union des deux paquets $\Pi(\psi')$ et $\Pi(\psi'')$.

On complète ces résultats dans le théorème \ref{complement} :  les valeurs possibles pour $a(\psi_u)$ dans $(ii)$ sont 
$2(n-m)+1$, auquel cas  $E_{\psi'}$ est la représentation triviale, ou bien  $2(n-m)+3$, seulement si $2m\geq (n+2)$ et on détermine $E_{\psi'}$ qui n'est pas un caractère en général.

D'autre part, on donne  conformément aux résultats de \cite{MR3}
une construction de $\pi_n(m)$
comme composante d'une  induite cohomologique à partir d'un $c$-Levi de $G$ isomorphe à un
produit de groupes unitaires ici nécessairement compacts et d'un groupe symplectique de rang plus petit. On induit un produit tensoriel extérieur de caractères
de ces groupes unitaires compacts, uniquement déterminé par le caractère infinitésimal de $\psi$ et par $a(\psi_u)$,  et  d'une  représentation  unipotente  du petit groupe symplectique que l'on détermine. Cette représentation unipotente est aussi de plus haut poids.

Le cas du caractère infinitésimal régulier, beaucoup plus simple,  est traité rapidement  dans la section \ref{carinfreg}
et le résultat (théorème \ref{thmcasreg}) est le suivant. On fixe un caractère infinitésimal régulier, et ceci détermine un entier
$a_{max}$ tel que les modules  de plus haut poids unitaires  soient indexés par  $a\in \{0,\ldots,a_{max}\}$
(le cas $a=0$ est celui des séries discrètes holomorphes). Notons $\pi_a$ le module indexé par $a$. 
Les paquets d'Arthur $\Pi(\psi)$ ayant ce même caractère infinitésimal sont des paquets d'Adams-Johnson, et la partie unipotente $\psi_u$
du paramètre est irréductible, de dimension $2 a(\psi_u)+1$. Le paquet  $\Pi(\psi)$  contient alors $\pi_a$ si et seulement si $a=a(\psi_u)$.

Pour les représentations $\sigma_{n,k}$, l'énoncé est le suivant (Théorème \ref{Thmsnk}). Remarquons
que $\sigma_{2k,k}=\pi_{2k}(k+1)$ et que ce cas est déjà traité dans le théorème ci-dessus.
\begin{thmI} On suppose que $n-1\geq 2k$.  Soit $\psi$ un $A$-paramètre pour $G=\Sp(2n,\bbR)$ 
dont le caractère infinitésimal est celui de $\sigma_{n,k}$.  Alors $\psi$ se décompose comme en (\ref{formpsi}) et 
 le paquet $\Pi(\psi)$ contient 
 $\sigma_{n,k}$ si et seulement si $a(\psi_u)=2(n-k)+1$ et $\psi$ contient $\sgn^k\boxtimes R[2(n-k)+1]$.
\end{thmI}

J. Adams a décrit tous les modules de plus haut poids comme induites cohomologiques dans \cite{Ad87} .
Avec l'analyse du cas où les paramètres sont unipotents, ceci permet de conclure assez rapidement que  les conditions
énoncées dans les théorèmes sont suffisantes.  L'essentiel de notre travail consiste donc  à montrer 
 qu'elles sont nécessaires.  Ensuite, pour ces représentations, nous établissons  la propriété de multiplicité un. Enfin, 
  si $\pi=\pi_n(m)$ ou $\pi=\sigma_{n,k}$ et si  $\psi$ est un $A$-paramètre tel que $\pi \in \Pi(\psi)$, 
  nous calculons le caractère $\rho_{\pi}$ de $A(\psi)$  dans les propositions \ref{Apsipinm} et \ref{Apsisnm}.

Nous allons utiliser de manière essentielle certains arguments globaux, qui nécessitent l'introduction de notations
spécifiques et le rappel de résultats de la littérature, en particulier 
 \cite{pourhowe} et surtout la formule de Siegel-Weil de Kudla et Rallis ainsi que ses généralisations,
ce qui est fait dans la deuxième moitié de l'article.  Certaines démonstrations sont de ce fait différées.
Le calcul des caractères $\rho_\pi$ se fait par ces méthodes globales : on se ramène à certains 
cas non archimédiens connus grâce à une formule de produit (Prop. \ref{formprod}) qui n'est rien d'autre que la formule de multiplicité d'Arthur.

\

Donnons un aperçu du contenu de cet article et de la démonstration de ces résultats.
La section \ref{Sec2} introduit des notations et fait quelques rappels concernant  l'induction cohomologique, 
que l'on particularise au cas du groupe symplectique. Une condition nécessaire pour l'occurence d'un $K$-type
dans une induite cohomologique est donnée sous la forme de l'inégalité (\ref{ineqfond}). Ceci nous servira
à montrer que les représentations $\pi(\mu)$ ne peuvent apparaitre dans certaines induites cohomologiques
et permet de contrôler leur appartenance à certains paquets (ou leur multiplicité dans ceux-ci).
Dans la section \ref{MUPHP}, on introduit les modules de plus haut poids pour le groupe symplectique $\Sp(2n,\bbR)$, dont on rappelle la classification
(\cite{EHW}) et la description en terme de correspondance de Howe (\cite{KV}):  ce sont les images par celles-ci de représentations de dimension finie de 
groupes orthogonaux compacts $\Or(0,2\ell)$, avec $\ell\leq n$.
Dans la section \ref{PAQART}, on rappelle ce que sont les paquets d'Arthur pour le groupe symplectique (\cite{Art13})
et la description de ceux-ci donnée dans \cite{MR3}. On cherche donc à déterminer tous les couples $(\pi,\psi)$, où 
$\pi$ est un module unitaire de plus haut poids et $\psi$ est un paramètre d'Arthur tels que $\pi$ appartienne au paquet d'Arthur 
$\Pi(\psi)$. 
Une  condition nécessaire évidente sur le caractère infinitésimal contraint fortement les paramètres possibles  (voir  Equation (\ref{decApar}) et le lemme \ref{longpsiu}).
Ceux qui vont nous intéresser sont construits à partir des paquets unipotents de groupes symplectiques
 plus petits par induction cohomologique, le cas des paquets unipotents étant traité dans \cite{pourhowe}.
Dans la section \ref{ParUnip}, on considère les couples  $(\pi,\psi)$ comme ci-dessus, avec de plus $\psi$ unipotent.
Ils sont obtenus par une recette bien connue des experts. On enlève de la décomposition en irréductibles du paramètre son plus gros bloc 
(le paramètre comporte au plus trois blocs), et l'on obtient un paramètre d'Arthur qui doit être celui de la représentation triviale pour un groupe spécial orthogonal compact 
$\SO(0,2\ell)$. La représentation $\pi$ doit alors être l'image par la correspondance de Howe de la représentation triviale ou bien du déterminant 
de ce groupe orthogonal compact (théorème \ref{thmunip}). Ce résultat est un cas particulier du théorème \ref{lemp} qui sera démontré par des méthodes globales.
Dans le premier cas, $\pi=\pi_n(m)$, avec $2m\leq n+1$ et dans le second
$\pi=\sigma_{n,k}$, $2k\leq n$. Dans la section \ref{carinfreg}, on traite le cas du caractère infinitésimal régulier.
 Dans la section \ref{scal}, on considère le cas des modules unitaires  $\pi_n(m)$ de plus haut poids scalaire $-m$, et l'on énonce le théorème principal
\ref{mainthm} qui donne la liste des paramètres d'Arthur  $\psi$ tel que $\pi_n(m)\in \Pi(\psi)$. Ceci est complété par le théorème 
\ref{complement} qui décrit $\pi_n(m)$ comme induite cohomologique et comme image par la correspondence de Howe à partir de groupes orthogonaux compacts.
Les représentations unipotentes $\sigma_{n,k}$ sont étudiées  plus avant dans la section \ref{sigmana}. On 
utilise les résultats de \cite{Ad87} pour donner des réalisations de ces représentations comme constituants d'induites cohomologiques, ce qui montre l'appartenance
de celles-ci à certains paquets (corollaire \ref{corsna}, proposition \ref{propsna}). On énonce ensuite le théorème \ref{Thmsnk}
qui affirme que la liste de ces paquets est complète.
Dans la section \ref{langlands}, on donne  les paramètres de Langlands des $\pi_n(m)$ 
et des $\sigma_{n,k}$. Les démonstrations sont globales, et sont données dans la section \ref{demparlang}.
On en tire le  corollaire \ref{borne} en utilisant  un résultat de \cite{Art13}
sur les exposants des représentations à l'intérieur d'un paquet, tiré de \cite{Art13}, et 
que nous utiliserons de manière essentielle.  
La démonstration des théorème \ref{mainthm} et \ref{complement}  proprement dite occupe les sections \ref{reduc} à \ref{condnec}.
La section \ref{reduc} énonce le résultat de réduction qui  permet une   récurrence.
 L'énoncé et la démonstration sont assez techniques, et  ne s'étendent pas 
au cas des modules unitaires de plus haut poids non scalaire; c'est là l'obstacle principal qui nous empêche d'obtenir des résultats complets pour
 tous les modules de plus haut poids.
On démontre ensuite  dans la section \ref{sensdirect} que les conditions du théorème sont suffisantes, et dans la section \ref{condnec}
qu'elles sont nécessaires. On utilise de manière cruciale analogue local de   Siegel-Weil   qui étend celui obtenu dans \cite{pourhowe}
pour les représentations unipotentes   (corollaire \ref{SW})  pour conclure, et qui est établi  par des méthodes globales dans la section \ref{siegelweil}.
 Dans la section \ref{Multun}, on vérifie que $\pi_n(m)$ apparaît dans un paquet $\Pi(\psi)$ avec multiplicité au plus un.
Dans la  section \ref{demThmsnk}, nous vérifions grâce au corollaire \ref{borne} que la liste des paquets contenant une représentation unipotente  $\sigma_{n,k}$
du théorème   \ref{Thmsnk} est complète, ainsi que la propriété de multiplicité un.

Dans la deuxième partie de l'article, nous utilisons des arguments globaux, et en particulier
des résultats établis aux places non archimédiennes. Dans la section \ref{discethasse},
nous rappelons les définitions du discriminant et de l'invariant de Hasse d'un espace quadratique $(V,Q)$ sur un corps local
et nous en donnons une normalisation adaptée à notre usage.
La section \ref{caracO} commence par quelques rappels sur les correspondances de Howe locale et globale, en particulier dans le cas d'une paire 
duale orthogonale paire/symplectique. On paramètre ensuite les caractères des groupes orthogonaux pairs
et l'on s'intéresse à leur image par la correspondance de Howe. Le résultat principal de cette section est le théorème \ref{lemp}
qui affirme l'appartenance de cette image à certains paquets d'Arthur unipotent. Dans le cas archimedien, on obtient en particulier 
l'énoncé du théorème \ref{thmunip}.
Dans la section \ref{calcrhopi}, on calcule le caractère $\rho_\pi$ de $A(\psi)$ attaché à une représentation 
$\pi$ d'un paquet $\Pi(\psi)$ avec $\pi$ et $\psi$ comme dans le théorème \ref{lemp} dans certains cas,
 en particulier dans le cas archimédien, l'outil principal étant la formule de produit de la proposition \ref{formprod}
 et les calculs dans le cas non archimédien de \cite{elementaire} rappelés dans le  lemme \ref{Epique}.
Enfin, dans la section \ref{rhopiR}, nous traduisons ces résultats dans le cas particulier des modules 
de plus haut poids $\pi_n(m)$ et $\sigma_{n,k}$ et des paquets unipotents, puis nous donnons la formule 
pour un paquet quelconque  en utilisant \cite{MR3}.

Terminons cette introduction par quelques mots concernant des applications possibles de nos résultats.
Soit $S_{\underline k}(\Sp(2n,\bbZ))$ l'espace des formes modulaires paraboliques  (pour le groupe  $\Sp(2n,\bbZ)$)  de Siegel holomorphes
à valeurs dans la représentation algébrique irréductible  de $\GL_n(\bbC)$ de plus haut poids $\underline k=(k_1\geq k_2\geq \cdots \geq k_n)$.
 La dimension de ces espaces comme fonction de $\underline k$ est un résultat classique de la théorie des formes modulaires
 lorsque $n=1$. Lorsque $\underline k$ est scalaire et $n=2$ ou $n=3$, ce sont des résultats dus  respectivement à   Igusa \cite{Igu} et Tsuyumine \cite{Tsuy}, 
  et lorsque $n=2$ et que $\underline k=(k_1,k_2)$ avec $k_2\geq 5$, c'est  la formule de Tsushima \cite{Tsu}. 
  Tsuyumine donne aussi une formule pour certaines familles de poids   non scalaires  particulières et  $k_n>n$.
  Récemment, la dimension de  ces espaces  à été déterminée de manière algorithmique pour de nombreuses
autres valeurs de $\underline k$ en rang $3$ (\cite{CR}), puis  en rang $n\leq 7$ par O. Taïbi \cite{taibi2}.
Ces travaux sont basés sur la classification du spectre automorphe discret des groupes classiques par Arthur \cite{Art13}
qui relient les espaces  $S_{\underline k}(\Sp(2n,\bbZ))$ et certaines  représentations automorphes cuspidales 
algébriques auto-duales des groupes généraux linéaires. Les poids $\underline k$ considérés
vérifient $k_n>n$, et en conséquence une forme modulaire propre pour les opérateurs de Hecke dans $S_{\underline k}(\Sp(2n,\bbZ))$
engendre une représentation  $\pi(\underline k)$ de $\Sp(2n,\bbR)$ qui est une série discrète holomorphe, et donc en particulier un module unitaire de plus bas poids,
donné par $\underline k$, vu ici comme plus haut poids d'une représentation du compact maximal $K\simeq \U(n)$ de $G=\Sp(2n,\bbR)$.
Les calculs de \cite{CR} et \cite{taibi2} utilisent les formules de multiplicité d'Arthur, et il faut en particulier connaître
la liste des paramètres d'Arthur $\psi$  tels que le paquet associé $\Pi(\psi)$ contienne 
$\pi(\underline k)$ et pour chacun d'eux, la représentation $\rho_{\pi(\underline k)}$ du groupe
des composantes du centralisateur de $\psi$ dans $\widehat G$ attachée à $\pi(\underline k)$ et à $\psi$ par Arthur.
Comme nous l'avons remarqué plus haut, 
les résultats de \cite{AMR} montrent alors que les paquets $\Pi(\psi)$ contenant $\pi(\underline k)$ sont exactement les paquets d'Adams-Johnson
\cite{AdJo} construits par induction cohomologique  ayant le bon caractère infinitésimal et un paramètre de partie unipotente de dimension 1.
 La représentation $\rho_{\pi(\underline k)}$
(dans ce cas, un caractère) est déterminée
dans \cite{CR}, Chapter 9.
En général, on peut se demander ce que les résultats de cet article peuvent nous apprendre sur les espaces  $S_{\underline k}(\Sp(2n,\bbZ))$, 
si l'on ne  suppose plus $k_n>n$
(le passage des modules de plus bas poids à ceux de plus haut poids considérés dans l'article n'est qu'une affaire de convention, 
ou bien l'on passe aux contragrédients).
 La représentation  $\pi(\underline k)$ de $\Sp(2n,\bbR)$ engendrée par une   forme modulaire propre pour les opérateurs de Hecke
 dans $S_{\underline k}(\Sp(2n,\bbZ))$  est alors le  module unitaire irréductible de plus bas poids $\underline k$.
On peut alors envisager  d'utiliser la formule de multiplicité d'Arthur comme dans le cas 
régulier $k_n>n$  pour relier ces espaces de formes de Siegel aux représentations automorphes cuspidales 
algébriques auto-duales des groupes généraux linéaires. Dans un esprit légèrement différent, Chenevier et Lannes montrent dans \cite{CL}
par des techniques très détournées de fonctions $L$  que pour un poids scalaire $\underline k=(k,\ldots,k)$,  $k\leq 12$, l'espace 
$S_{\underline k}(\Sp(2n,\bbZ))$ est nul si $k\leq n\leq 2k$ (sauf peut-être dans pour  $k=12$ et $g=24$).
Les résultats de cet article permettent de simplifier notablement les démonstrations de \cite{CL} et d'aborder les cas  suivants.

\ 

{\bf Remerciements}.
Nous remercions G. Chenevier qui a attiré notre attention sur le  problème considéré dans cet article, et qui en  relu avec attention
les versions successives avec attention en signalant de nombreuses coquilles, en particulier une dans les formules de la proposition 
 \ref{Apsipinm}.
Le deuxième auteur a bénéficié d'une aide de  l'agence nationale de la recherche 
ANR-13-BS01-0012 FERPLAY.

\section{Généralités et notations}\label{Sec2}

\subsection{Paires paraboliques}\label{parthetastable}
Soit  $G$  le groupe des points réels d'un groupe algébrique connexe réductif défini sur $\bbR$. On fixe une involution de 
Cartan $\theta$ de $G$, et l'on note $K$ le sous-groupe des points fixes de $\theta$ : c'est un sous-groupe compact maximal de $G$.
On suppose que  $G$ et $K$ sont de même rang;    autrement dit, $G$ possède un sous-groupe de Cartan $T$
inclus dans $K$ et donc compact.
On note $\frt_0$, $\frk_0$  et $\frg_0$ les  algèbres de Lie respectives de $T$,$K$  et $G$ et $\frt$, $\frk$  et $\frg$ leur complexifiées.
Les  sous-algèbres paraboliques $\theta$-stables de $\frg$ sont obtenues de la manière suivante. On fixe un 
élément $\nu\in \sqrt{-1} \frt_0^*$, et l'on pose :
\[\frl=\frg^\nu=\frt\oplus \left( \bigoplus_{\alpha\in \Delta(\frg,\frt), \bil{\nu}{\alpha}=0}\frg^\alpha\right), 
\quad \fru =\bigoplus_{\alpha\in \Delta(\frg,\frt), \\ \bil{\nu}{\alpha}>0}\frg^\alpha, \quad \frqqq=\frl\oplus \fru,  \quad L=\mathrm{Norm}_G(\frqqq).
\]
Dans cet article, nous appellerons paire parabolique 
une paire $(\frqqq,L)$ obtenue comme ci-dessus, avec $\frqqq$ sous-algèbre parabolique $\theta$-stable de $\frg$. Le  
 sous-groupe $L$ de $G$  sera appelé $c$-Levi de $G$ (terminologie de Shelstad \cite{Sh15}).

\medskip 

On note $\caR_{\frqqq,L,G}^i$ le foncteur d'induction cohomologique de Vogan-Zuckerman ({\sl cf.} \cite{Vgreen}, \S 6.3.1)
en degré $i$, de la catégorie des $(\frl,K\cap L)$-modules vers la catégorie des $(\frg,K)$-modules.
Dans ce contexte, le degré qui nous intéresse particulièrement, et même exclusivement,  est $S=\dim(\fru\cap \frk)$, et dans l'article, nous écrirons
$\caR_{\frqqq,L,G}^S$ sans préciser de nouveau ce qu'est $S$.

Si $\Lambda$ est un caractère unitaire  de $L$, on note  $\lambda$ sa différentielle, que l'on voit comme un élément
 de $i\frt_0^*$.
On pose alors 
$A_\frqqq(\Lambda)=\caR^S_{\frqqq, L,G}(\Lambda)$.
Si le groupe $L$ est  connexe, $\lambda$ détermine $\Lambda$ et l'on note alors cette représentation $A_\frqqq(\lambda)$.

Sous certaines conditions sur le caractère infinitésimal de la représentation $\sigma$ de $L$ que l'on induit, on a des résultats d'annulation, d'irréductibilité et d'unitarité
des modules $\caR_{\frqqq,L,G}^i(\sigma)$. Nous renvoyons à \cite{KnVo} pour les définitions du (weakly) good range, du  (weakly) fair  range
et des représentations faiblement unipotentes, pour lesquels on a les résultats
suivants : dans le weakly good range les modules $\caR_{\frqqq,L,G}^i(\sigma)$ sont nuls si $i\neq S=\dim(\fru\cap \frk)$.
Si $\sigma$ est irréductible et dans le good range  (resp. weakly good range), $\caR_{\frqqq,L,G}^S(\sigma)$ est irréductible (resp. irréductible ou nul).
Si $\sigma$ est unitaire  et dans le weakly good range, $\caR_{\frqqq,L,G}^S(\sigma)$ est unitaire.
Si $\sigma$ est faiblement unipotente et dans le weakly fair range, alors les modules $\caR_{\frqqq,L,G}^i(\sigma)$ sont nuls si $i\neq S$ et 
$\caR_{\frqqq,L,G}^S(\sigma)$ est unitaire si $\sigma$ est de plus unitaire. 
En revanche, on n'a pas de résultat d'irréductibilité en général dans le weakly fair range, ni même le fair range.

\subsection{Le groupe symplectique}
On suppose $\bbR^{2n}$ (identifié à $\caM_{2n,1}(\bbR)$, les matrices colonnes) muni de sa forme symplectique usuelle, c'est-à-dire, si $X,Y\in \bbR^{2n}$ 
\[ (X\vert Y)={}^tX J Y \text{ où }  J=\begin{pmatrix}0_n&I_n\\-I_n&0_n\end{pmatrix}.\]
Soit $G=\Sp(2n,\bbR)$ le groupe des isomorphismes de $(\bbR^{2n},  (.\vert .))$, que l'on muni
de l'involution de Cartan $\theta: g\mapsto{}^tg^{-1}$. Le sous-groupe de $G$ des  points fixes sous $\theta$
est un sous-groupe compact maximal de $G$, que l'on note $K$, et qui est isomorphe
au groupe unitaire $\U(n)$. On note $\frg_0$ et $\frk_0$ les sous-algèbres de Lie respectives de 
$G$ et $K$, réalisées comme sous-algèbre de Lie de $\caM_{2n}(\bbR)$.
Pour tout $(a_1,\ldots ,a_n)\in \bbR^n$, on pose :
\[t(a_1,\ldots,a_n)=\left(  \begin{array}{llll|llll}   & & & & a_1&&&\\
& & && & a_2 &&\\
&&& &&&\ddots&\\
&&& &&&& a_n\\
\hline
-a_1 & && & & &&\\
       &-a_2 &&& &&&\\
    & &\ddots &&   &&&\\
    &&& -a_n&&&&
  \end{array}\right). \]
Alors $\frt_0:=\{ t(a_1,\ldots,a_n), \, (a_1,\ldots ,a_n)\in \bbR^n \}$ est une sous-algèbre de Cartan de $\frk_0$ 
et aussi de $\frg_0$.

Notons $\frg,\frk,\frt$ les complexifications des algèbres de Lie $\frg_0,\frk_0,\frt_0$,
respectivement. Soient $\Delta(\frg,\frt)$, $\Delta(\frk,\frt)$ les systèmes de racines de 
$\frg$ et $\frk$ respectivement, relativement à la sous-algèbre de Cartan $\frt$.
On a 
\[\Delta(\frg,\frt)=\{ \pm(e_i \pm e_j), \, 1\leq i<j \leq n\} \cup \{ \pm 2e_i,  1\leq i \leq n\} , \]
\[\Delta(\frk,\frt)=\{ \pm(e_i - e_j), \, 1\leq i<j \leq n\}, \]
où $e_i\in  \sqrt{-1}\,  \frt_0^*\subset \frt^*$ est la forme linéaire $t(a_1,\ldots,a_n)\mapsto   \sqrt{-1} \, a_i$.
On fixe les systèmes de racines positives 
\[\Delta^+(\frg,\frt)=\{ (e_i \pm e_j), \, 1\leq i<j \leq n\} \cup \{  2e_i,  1\leq i \leq n\} , \]
\[\Delta^+(\frk,\frt)=\{ (e_i - e_j), \, 1\leq i<j \leq n\}. \]
On identifie $\frt^*$ et $\bbC^n$ grâce à la base $(e_i)_{1\leq i\leq }$ de $\frt^*$, et de même pour $\frt$ grâce à la base duale. 
On note encore $\theta$ la différentielle de l'involution de Cartan de $G$, et l'on pose 
\[\frg_0=\frk_0\stackrel{\theta}{\oplus}\frp_0; \quad \frg=\frk\stackrel{\theta}{\oplus}\frp. \]

On a alors 
\[\frp=\frp^+\oplus \frp^-, \quad \frp^+=\bigoplus_{ 1\leq i<j \leq n} \frg_{e_i+e_j}\oplus \bigoplus_{1\leq i \leq n }\frg_{2e_i} ,\quad
\frp^-=\bigoplus_{ 1\leq i<j \leq n} \frg_{-e_i-e_j}\oplus \bigoplus_{1\leq i \leq n }\frg_{-2e_i}.\]

\begin{defi}\label{qholo}
On dit qu'une sous-algèbre parabolique $\theta$-stable $\frqqq=\frl\oplus \fru$ de $\frg$  est {\sl holomorphe} si 
$\fru\cap \frp \subset \frp^-$.
\end{defi}

\begin{rmq}\label{antihol}
On aurait pu (dû?) appeler de telles sous algèbre paraboliques anti-holomorphes, mais la dichotomie 
holomorphe/anti-holomorphe n'est qu'une affaire de convention. Nos séries discrètes holomorphes
sont celles qui sont plus communément appelées anti-holomorphes dans la littérature.
\end{rmq}

On note $W$ le groupe de Weyl du système de racine $\Delta(\frg,\frt)$. Il agit sur $\frt^*$, identifié comme expliqué ci-dessus à $\bbC^n$,
par permutations et changements de signes des coordonnées. Via l'isomorphisme d'Harish-Chandra, le  
caractère infinitésimal d'une représentation de $G$ s'identifie à une $W$-orbite dans $\frt^*$.
Nous ne considérerons que des caractères infinitésimaux entiers. 

\begin{convention}\label{conv} Il sera commode de voir le caractère infinitésimal d'une représentation de $G=\Sp(2n,\bbR)$
de la façon suivante. On identifie  $\frt^*$ à  une sous-algèbre de Cartan  ${}^L\frt$ de l'algèbre de Lie duale ${}^L\frg=\frs \fro(2n+1,\bbC)$.
Via le plongement standard ${}^L\frg=\frs \fro(2n+1,\bbC) \hookrightarrow  \frg\frl(2n+1,\bbC)$, nous verrons un caractère infinitésimal comme un $(2n+1)$-uplet
décroissant, symétrique, où $0$ a une multiplicité impaire. 
\end{convention}

\subsection{Paires paraboliques maximales et induction cohomologique}\label{ppm}
Dans ce paragraphe,  $G=\Sp(2n,\bbR)$ et l'on continue avec les notations de du paragraphe précédent.
Soit $(\frqqq,L)$ une paire parabolique pour $G$, et l'on suppose que la sous-algèbre parabolique $\theta$-stable $\frqqq=\frl\oplus\fru$ est 
maximale.  Une telle sous-algèbre est obtenue en prenant un élément de $\frt$ de la forme
\[  t_{p,q}=(\underbrace{1,\ldots,1}_{p}, \underbrace{0,\ldots,0}_{n-p-q},  \underbrace{-1,\ldots,-1}_{q}) \]
avec $p+q\leq n$.
On pose alors 
\begin{equation}\label{qpq}  \frl=\frl_{p,q}=\frg^{t_{p,q}}=\frt\oplus\left(  \bigoplus_{\alpha\in \Delta(\frg,\frt)\vert \alpha(t)=0 } \frg_\alpha  \right), 
\quad \fru=\fru_{p,q}= \bigoplus_{\alpha\in \Delta(\frg,\frt)\vert \alpha(t)>0 } \frg_\alpha.\end{equation}
Les racines de $\frt$ dans $\frl$ sont :
\begin{align*}
&\pm(e_i-e_j), &\,  1\leq i<j\leq p  \text{ ou }   n-q+1 \leq i<j\leq n,\\
&  \pm(e_i\pm e_j), \, \pm 2e_i, &\, p+1\leq i<j \leq n-q, \\
&\pm(e_i+e_j),& \,   1\leq i \leq p\leq  n-q+1\leq j\leq n. 
 \end{align*}
 On choisit comme système de racines positives :
 \[\Delta^+(\frl,\frt)= \left\{  \begin{array}{ll}  -( e_i-e_j),  &  1\leq i<j\leq p, \\
   e_i-e_j , &  n-q+1 \leq i<j\leq n \\
  e_i- e_j, \, -(e_i+e_j) , - 2e_i,  &\, p+1\leq i<j \leq n-q\\
  -(e_i+e_j), & \,   1\leq i \leq p\leq  n-q+1\leq j\leq n
    \end{array} \right\}.\]
 On pose 
 $\delta(\frl)=\frac{1}{2}\sum_{\alpha\in \Delta^+(\frl,\frt)}\alpha$
 \[{\tiny =\frac{1}{2} (
 \underbrace{-q-p+1,-q-p+3 , \ldots, -q+p-1}_{p} ,\underbrace{ -2,-4 ,\ldots -2(n-p-q)}_{n-p-q},\underbrace{  -p+q-1,\ldots,  -p-q+3,\ldots ,-p-q+1}_{q}  )}. \]
Les racines de $\frt$ dans $\fru\cap \frp$ sont : 
\begin{align*}
& e_i+e_j, \,   2\epsilon_i, \, &\, 1\leq i<j  \leq p   \\
 & -e_i-e_j, \, -2\epsilon_i,\, &\, n-q+1\leq i<j\leq n, \\
& e_i+e_j, \,  &\, 1\leq i  \leq p<j\leq n-q   \\
 & -e_i-e_j, \,  &\,  p+1\leq i\leq  n-q <  j\leq n  . 
\end{align*} 
Les racines de $\frt$ dans $\fru\cap \frk$ sont : $e_i-e_j $,  $   1\leq i \leq p< j\leq n$    ou    $ p+1\leq i\leq n-q<j\leq n$.
On pose :
\begin{align*}& \delta(\fru\cap \frp)= \frac{1}{2}\sum_{\alpha\in \Delta(\fru\cap \frp)}\alpha
={\scriptsize  \frac{1}{2} (\underbrace{n-q+1,\ldots,n-q+1}_{p}, \underbrace{p-q,\ldots, p-q}_{n-p-q}, \underbrace{-n+p-1, \ldots -n+p-1}_{q})},\\
&\delta(\fru\cap \frk)= \frac{1}{2}\sum_{\alpha\in \Delta(\fru\cap \frk)}\alpha 
= {\scriptsize  \frac{1}{2} (\underbrace{n-p,\ldots,n-p}_{p}, \underbrace{-p+q,\ldots, -p+q}_{n-p-q},
\underbrace{-(n-q), \ldots -(n-q)}_{q})},   \\
&\delta(\fru)=\delta(\fru\cap \frp)+ \delta(\fru\cap \frk)={\scriptsize \frac{1}{2}(2n-p-q+1) (\underbrace{1,\ldots,1}_{p}, \underbrace{0,\ldots, 0}_{n-p-q},
\underbrace{-1, \ldots -1}_{q})},\\
& \delta_{p,q}=\delta(\frl)+\delta(\fru)\\
 &={\scriptsize (\underbrace{n-p-q+1, n-p-q+2,\ldots,n-q}_{p},  
\underbrace{-1, -2,\ldots ,-(n-p-q)}_{n-p-q}, \underbrace{-(n-q+1), -(n-q+2) \ldots, -n}_{q}  ) }.
\end{align*}
De plus  $L$ est isomorphe à $\Sp(2(n-p-q),\bbR)\times\U(p,q)$. 
Nous allons considérer des induites  cohomologiques  à partir de cette paire $(\frqqq,L)$ 
de la forme 
\begin{equation}\label{zalum} \caR^S_{\frqqq,L,G}( \rho \boxtimes \Lambda).\end{equation}
 où  $\Lambda$ est un caractère unitaire du facteur $\U(p,q)$ de $L$, et $\rho$ une représentation du facteur 
 $\Sp(2(n-p-q),\bbR)$.

Le centre $\frz$ de $\frl$ est contenu dans $\frt$, il est de dimension $1$ et engendré par $t_{p,q}$.
Soit
\begin{equation}\label{lambda1}\lambda=(\underbrace{y,\ldots,y}_{p}, \underbrace{0,\ldots, 0}_{n-p-q},
\underbrace{-y, \ldots ,-y}_{q}) \end{equation}
la différentielle de $\Lambda$.
Le caractère infinitésimal de $\rho \boxtimes \Lambda$ est donné par  un élément $\nu$ de $\frt^*$ de la forme 
\[\nu=(\underbrace{y,\ldots,y}_{p}, \underbrace{*,\ldots, *}_{n-p-q},
\underbrace{-y, \ldots ,-y}_{q})+\delta(\frl). \]

 La condition pour que l'induction soit dans le weakly   fair range est  que pour tout $\alpha\in \Delta(\fru,\frt)$,  
\[   \bil{\nu+\delta(\fru)}{\alpha_{\vert \frz}} \geq 0.\]
Or pour un tel $\alpha$,  $\bil{\delta(\frl)}{\alpha_{\vert \frz}} =0$, et donc la condition est,  pour tout $\alpha\in \Delta(\fru,\frt)$,
\[  y+  \frac{1}{2}(2n-p-q+1) \bil{ (\underbrace{1,\ldots,1}_{p}, \underbrace{*,\ldots, *}_{n-p-q},
\underbrace{-1, \ldots,  -1}_{q}) }{\alpha_{\vert \frz}} \geq 0.\]
Ceci est équivalent à $ y+  \frac{1}{2}(2n-p-q+1) \geq 0$. Posons $t=2y+n+1-p-q$, de sorte que 
\begin{equation}\label{y} \lambda={\scriptsize (\underbrace{\frac{t+p+q-1}{2}-n,\ldots,\frac{t+p+q-1}{2}-n}_{p}, \underbrace{0,\ldots, 0}_{n-p-q},
\underbrace{n-\frac{t+p+q-1}{2}, \ldots , n-\frac{t+p+q-1}{2}}_{q})} 
\end{equation}
 et la condition de weakly fair range est  $t\geq 0$.

 \medskip
 {\bf Une inégalité fondamentale}.
 Nous allons maintenant rappeler un résultat de \cite{KnVo} qui donne des informations sur les $K$-types de 
$ \caR^S_{\frqqq,L,G}( \rho \boxtimes \Lambda)$
 Dans \cite{KnVo},  Equation (4.71) p. 272,  il est introduit 
un élément de $\frt_0$,  que l'on note ici $h_{\frqqq}$ et qui est dual de 
$\delta(\fru)$, c'est-à-dire 
\begin{equation}
h_{\frqqq}= 2\times (2n-p-q+1)^{-1}\times (\underbrace{1,\ldots,1}_{p}, \underbrace{0,\ldots, 0}_{n-p-q},
\underbrace{-1, \ldots ,-1}_{q}).
\end{equation}

L'élément $h_{\frqqq}$ agit par un scalaire $c$  dans la représentation  $ \rho \boxtimes \Lambda$
On fixe une représentation irréductible de $K$ de plus haut poids $\mu=(-m_n,\ldots, -m_1)$ et l'on suppose que cette
représentation est un $K$-type de $\caR^S_{\frqqq,L,G}(\rho\boxtimes \Lambda)$.
L'inégalité de la troisième ligne de la page 369 de \cite{KnVo} est 
\begin{equation}\label{p369}
\mu(h_{\frqqq}) \geq c+2\delta(\fru\cap \frp)(h_{\frqqq}) . \end{equation}
 Le terme de gauche est 
\[  \mu(h_{\frqqq}) =2\times (2n-p-q+1)^{-1}\times  \left(- \sum_{j=i}^p m_{n-i+1} +  \sum_{j=1}^q  m_j \right). \]
On a aussi 
\[ 2\delta(\fru\cap \frp)(h_{\frqqq}) =  2\times (2n-p-q+1)^{-1}  \times (p(n-q+1)+q(n-p+1)) . \]

Le scalaire $c$ est alors donné par 
\[c=\bil{\nu}{h_{\frqqq'} }= 2 (2n-p-q+1)^{-1}  \left( (p+q)y + \bil{\delta(\frl)}{h_{\frqqq'}}  \right)=
 2 (2n-p-q+1)^{-1}   (p+q)y. \]

L'inégalité (\ref{p369}) s'écrit donc en simplifiant par $2\times (2n-p-q+1)^{-1}$ : 
\begin{align*}  - \sum_{j=i}^p m_{n-i+1} +  \sum_{j=1}^q  m_j &\geq  (p+q)\left(\frac{t+p+q-1}{2}-n\right)+ (p(n-q+1)+q(n-p+1))
\\
&\geq  (p+q)\left(\frac{t+p+q+1}{2}\right)-2pq.
 \end{align*}
Si tous les $m_j$ sont égaux à $m$, on obtient 
\begin{equation}\label{ineqfond} m(q-p)\geq  (p+q)\left(\frac{t+p+q+1}{2}\right)-2pq.  \end{equation}

\section{Modules unitaires de plus haut poids pour $\Sp(2n,\bbR)$}\label{MUPHP}

\subsection{Modules de plus haut poids}\label{ss:modphp}

On appelle module de plus haut poids de $G=\Sp(2n,\bbR)$ un $(\frg,K)$-module irréductible $\pi$
admettant un vecteur non nul annulé par une sous-algèbre de Borel de $\frg$.
Si $\pi$ n'est pas de dimension finie, nécessairement une telle sous-algèbre de Borel 
$\frb$ est  à conjugaison par $G$ près 
\[\frb^+=\frb_\frk\oplus \frp^+  \text{ ou bien } \frb^-=\frb_\frk\oplus \frp^-, \]
où $\frb_\frk=\frt\oplus \bigoplus_{\alpha\in \Delta^+(\frk,\frt)} \frg_\alpha$ est une sous-algèbre de Borel de $\frk$.
Dans le premier cas, on dit que $\pi$ est holomorphe, et dans le second, on dit que $\pi$ est antiholomorphe ({\sl cf.} Rmq \ref{antihol}).

Comme un  $(\frg,K)$-module irréductible de dimension finie possède pour toute sous-algèbre de Borel $\frb$ de $\frg$ 
un vecteur non nul annulé par $\frb$, un tel module est à la fois holomorphe et antiholomorphe, mais si $\pi$ est de dimension 
infinie, les deux possibilités sont exclusives.

Dans cet article, nous ne nous intéressons qu'aux modules $(\frg,K)$-holomorphes, les antiholomorphes  
étant leur contragrédient.
Les modules unitaires holomorphes de $G$ sont décrits dans \cite {KV} (voir aussi  \cite{EHW} et \cite{Jak}). 
Donnons cette description.

Si  $(\delta,V)$ une représentation irréductible de dimension finie de $K$, sa différentielle, encore notée
$\delta$, est une représentation   irréductible  $\frk$, qui détermine complètement $(\delta,V)$ puisque 
$K$ est connexe.
Soit $\mu=\mu_\delta=(m_1,\ldots,m_n)\in \frt^*$ le plus haut  poids de $(\delta,V)$, relativement à $\Delta^+(\frk,\frt)$, on a donc
 $m_1\geq m_2\geq .... \geq m_n$ et les $m_i$ dans $\bbZ$.

\begin{thm} \label{hw}Soit $\mu=(m_1,\ldots,m_n)\in \frt^*$ avec  $m_1\geq m_2\geq .... \geq m_n$ et les $m_i$ entiers.
Notons  $u$ le nombre de $m_i$ égaux à $m_n$ et $v$ le nombre de $m_i$ égaux à $m_n +1$.
Supposons que $m_n\geq n-(u+v/2)$.
Notons $(\delta,V)=(\delta_{\mu^*}, V_{\mu^*})$ la représentation irréductible de  dimension finie de $K$ de plus haut poids $\mu^*=
(-m_n,\ldots ,-m_1)$. 
Alors il existe un unique module unitaire holomorphe $(\pi(\mu),W_\mu)$ contenant le $K$-type $(\delta, V)$
avec multiplicité un (et l'on peut alors supposer que $V\subset W_\mu$),
 de tel sorte qu'un vecteur non nul de plus haut poids $v$ dans $V$ soit annulé
par $\frb_\frk\oplus \frp^+ $. Tout module unitaire holomorphe de $G$ est  caractérisé ainsi.
\end{thm}

\begin{rmq}   \label{carinfhw}
Avec la convention \ref{conv}, le caractère infinitésimal est constitué des entiers $m_1-1,\ldots,m_n-n$, de leur opposés, et de $0$, 
réarrangés pour les mettre dans l'ordre décroissant.
\end{rmq}

\begin{notation} \label{pimn} Nous allons particulièrement étudier le cas des modules holomorphes de plus haut poids scalaire, 
c'est-à-dire les $\pi(\mu)$ avec $\mu=(\underbrace{m,\ldots,m}_{n})$. Comme dans l'introduction, nous notons ces représentations 
$\pi_n(m)$ ou même $\pi(m)$ si le rang $n$ du groupe symplectique que nous considérons est donné clairement  par le contexte.
 \end{notation}

\subsection{ Description par la correspondance de Howe} 
Les modules unitaires holomorphes $\pi(\mu)$
du théorème ci-dessus sont dans l'image de la correspondance de Howe (\cite{Ho})
pour des paires duales de la forme $(\Or(0,2\ell),\Sp(2n,\bbR))$ avec $\ell\leq n$.

\begin{rmq} \label{rmqpqqp}Dans cet article,  nous allons utiliser de manière essentielle la  correspondance de Howe  ({\sl cf.} \cite{Ho})
entre groupe symplectiques et groupes orthogonaux pairs. Celle ci est définie  en partant d'un espace symplectique  $W$ de
dimension $2n$, de groupe d'automorphismes $\Sp(W)$
  et d'un espace $V$ de dimension $N=2\ell$  muni d'une forme
quadratique non dégénérée  $Q$  et de groupe d'automorphismes $\Or(V,Q)$.
Il faut aussi fixer un caractère additif du corps local sur lequel on travaille, ici $\bbR$, et l'on 
choisit le caractère $\psi_{\bbR,1}: x\mapsto \exp(2i\pi x)$.
Comme tous les espaces symplectiques de dimension $2n$ sont isomorphes, on peut prendre  $\Sp(W)=\Sp(2n,\bbR)$.
La classe d'isomorphisme de l'espace quadratique $(V,Q)$ est elle déterminée par la signature $(p,q)$ de
la forme quadratique $Q$ et l'on peut prendre  $\Or(V,Q)=\Or(p,q)$.
La correspondance de Howe dépend de $(V,Q)$ 
et pas simplement de son groupe d'isomorphisme,  et ainsi même si $\Or(q,p)=\Or(V,-Q)=\Or(V,Q)=\Or(p,q)$,
il importe de bien les distinguer.
\end{rmq}

Introduisons des notations pour les représentations irréductibles des groupes compacts $\Or(0,2\ell)$.
Les représentations de $\SO(0,2\ell)$, groupe compact et connexe, sont déterminées par leur plus haut poids, 
que l'on peut voir comme un $n$-uplet $(\nu_1,\ldots,\nu_\ell)$, où les $\nu_i$ sont des entiers, avec 
$\nu_1\geq \nu_2\geq \ldots \geq \vert \nu_{\ell}\vert$. 
Si est une représentation irréductible de $\Or(0,2\ell)$, alors soit sa restriction à $\SO(0,2\ell)$ est réductible, 
somme de deux représentations irréductibles de plus haut poids respectifs $(\nu_1,\ldots,\nu_\ell)$ et 
$(\nu_1,\ldots,-\nu_\ell)$, avec $\nu_\ell>0$, soit sa restriction à $\SO(0,2\ell)$ est irréductible, de plus haut poids 
$(\nu_1,\ldots,\nu_{\ell-1},\nu_{\ell}=0)$. Dans le premier cas, on note $[\nu_1,\ldots,\nu_\ell]_+$ cette représentation, et 
dans le second cas, il y a deux extensions possibles de 
la représentation de $\SO(0,2\ell)$ à $\Or(0,2\ell)$, que  l'on  note  $[\nu_1,\ldots,\nu_\ell]_+$ et  $[\nu_1,\ldots,\nu_\ell]_-$
(voir \cite{KV} pour la façon de les distinguer).
Le résultat suivant est dû à Kashiwara-Vergne \cite{KV}.

\begin{thm}\label{KV}
Reprenons $\mu=(m_1,\ldots,m_n)$, $u$ et $v$  comme dans le théorème \ref{hw}, et $\pi(\mu)$ 
le  module unitaire holomorphe de $\Sp(2n,\bbR)$ correspondant. 

(a) Supposons $m_n>n$. Alors $\pi(\mu)$ est image par la   correspondance de Howe 
pour la paire duale $(\Or(0,2n),\Sp(2n,\bbR))$  de la représentation 
$[m_1-n,\ldots,m_n-n  ]_+$ de $\Or(0,2n)$.

(b) Supposons $m_n=n-a$ avec $0\leq a\leq u$ et posons $\ell=m_n=n-a$. Alors $\pi(\mu)$ est image par la  correspondance de Howe 
pour la paire duale $(\Or(0,2\ell),\Sp(2n,\bbR))$  de la représentation 
$[m_1-\ell,\ldots,m_\ell-\ell  ]_+$ de $\Or(0,2\ell)$.

(c) Supposons $m_n=n-u-b$ avec $2\leq 2b \leq v$ et posons $\ell=m_n=n-u-b$. Alors $\pi(\mu)$ est image par la correspondance de Howe 
pour la paire duale $(\Or(0,2\ell),\Sp(2n,\bbR))$  de la représentation 
 $[m_1-\ell,\ldots , m_{n-u-2b}-\ell, 0,\ldots,0]_-$ de $\Or(0,2\ell)$ (il y a $b$ zéros).

(d) On suppose  $ m_n=n-a $ avec $a\geq 0$   et en plus on suppose
  $ m_n=n-a \geq n+1-u/2$, c'est-à-dire $a\leq u/2-1$ (c'est donc un sous-cas de (b)),   et l'on pose $\ell=m_n-1$.
Alors $\pi(\mu)$ est image de  correspondance de Howe 
pour la paire duale $(\Or(0,2\ell),\Sp(2n,\bbR))$  de la représentation 
 la représentation notée $[m_1-\ell ,\ldots , m_{2\ell-n}-\ell, 0,\ldots,0]_-$ de $\Or(0,2\ell)$ (il y a $n-\ell$ zéros).
\end{thm}

\begin{rmq} \label{rmqKV} Le cas  $(d)$ est un sous-cas de $(b)$. 
Considérons les sous-cas suivants de $(b)$ :

$(b')$  $m_n=n-a=n-u$, c'est-à-dire $a=u$,

$(b'')$ $m_n=n-a$ avec $\frac{u-1}{2}\leq a <  u$. On a donc  $m_{n-a}=n-a$.   
Le cas $(b)$ est  alors la réunion disjointe des sous-cas $(b')$, $(b'')$ et $(d)$.
\end{rmq}

\section{Paramètres et paquets d'Arthur pour $\Sp(2n,\bbR)$}\label{PAQART}
 
 \subsection{Paramètres et paquets  d'Arthur} \label{ArtparSp}
 Soit $G$ le groupe des points réels d'un groupe algébrique réductif connexe $\mathbf G$ défini sur $\bbR$. 
 Rappelons qu'un paramètre d'Arthur (ou $A$-paramètre)  pour le  groupe $G$ est un morphisme
 $\psi:\, W_\bbR\times \SL_2(\bbC)\to {}^LG$,  
où $W_\bbR$ est le groupe de Weil de $\bbR$.

Soit donc $\psi$ un paramètre d'Arthur et notons $A(\psi)$ le groupe des composantes connexes du centralisateur de $\psi$
dans $\widehat G$ (en fait, il faut en général passer à un revêtement {\sl cf.}  \cite{Art13}, chapitre 9, mais ceci est inutile dans le cas 
du groupe $G=\Sp(2n,\bbR)$ qui va nous occuper ici).
Arthur suggère qu'au paramètre $\psi$ est attaché une  combinaison linéaire de représentations  irréductibles   
de $G$ à coefficients dans l'espace des fonctions sur le groupe  $A(\psi)$ à valeurs complexes et  invariantes par conjugaison. On note 
 ${\pi}^A(\psi)$ cette combinaison linéaire. Ces objets  ${\pi}^A(\psi)$ doivent être  compatibles à l'endoscopie. 
 Cela ne suffit pas à les définir dans le cas quasi-déployé. Pour compléter la définition dans ce cas et pour les groupes classiques, 
 Arthur ajoute la compatibilité à l'endoscopie tordue et cela suffit alors.  Limitons nous dans ce qui suit aux cas des groupes classiques.
  Le groupe $A(\psi)$ est alors abélien (c'est même un $2$-groupe), les fonctions invariantes par conjugaison sur ce groupe
   sont donc  des    combinaisons linéaires à coefficients complexes de caractères de ce groupe  et $\pi^A(\psi)$ est  donc une combinaison linéaire à 
   coefficients complexes de représentations irréductibles de $G\times A(\psi)$.

Quand le groupe classique est  quasi-déployé, Arthur montre dans \cite{Art13} que 
 cette combinaison linéaire  est en fait à coefficients dans  les  entiers positifs, et  ${\pi}^A(\psi)$ est 
  donc une représentation semi-simple (par construction). Il montre même mieux :  c'est une représentation unitaire de $G\times A(\psi)$.
  On peut décomposer cette représentation unitaire de $G\times A(\psi)$ en l'écrivant comme somme de produits tensoriels extérieurs
de   représentations irréductibles unitaires de $G$ avec des représentations de dimension finies du groupe $A(\psi)$. Si l'on note 
$\Pi(\psi)$ l'ensemble des  représentations irréductibles unitaires de $G$ qui interviennent dans cette décomposition, on a 
\[\pi^A(\psi)=\bigoplus_{\pi\in\Pi(\psi)} \pi\boxtimes \rho_\pi.\]
  L'ensemble $\Pi(\psi)$ est le {\sl paquet d'Arthur} attaché au paramètre $\psi$ et la dimension de la représentation $\rho_\pi$ est la {\sl multiplicité}
  de $\pi$ dans le paquet. Pour les groupes classiques, on conjecture que cette multiplicité est $1$, et ceci est démontré
  pour certaines familles de paramètres (voir \cite{MR3}  pour une discussion détaillée).

 Prenons maintenant $G=\Sp(2n,\bbR)$. 
Les problèmes qui nous occupent dans cet article sont 

-  déterminer les couples $(\pi(\mu),\psi)$, où $\pi(\mu)$ est un module holomorphe unitaire
 comme dans le théorème \ref{hw} et $\psi$ un $A$-paramètre, 
 tels que $\pi(\mu)\in \Pi(\psi)$.

- Montrer la propriété de multiplicité un, à savoir que la représentation $\rho_{\pi(\mu)}$ de $A(\psi)$ 
est de dimension 1.

- Calculer ce caractère $\rho_{\pi(\mu)}$ de $A(\psi)$.

 Remarquons
que les  $\pi^A(\psi)$ dépendent du choix des facteurs de transfert géométrique  et  en suivant Kottwitz et Shelstad,
 les facteurs de transfert sont normalisés par des choix de données de Whittaker. Ceci est l'objet de la section suivante.
 
\subsection{Donnée de Whittaker}

 Rappelons qu'une donnée de Whittaker d'un groupe quasi-deployé $\mathbf G$ est un couple
$(N,\chi_N)$ où $N=\mathbf{N}(\bbR)$ est le groupe des points réels du  radical unipotent d'un sous-groupe de Borel $\mathbf B$
de $\mathbf G$ défini sur  $\bbR$ et d'un caractère non-dégénéré $\chi_N$ de $N$  (voir \cite{ABV}, chapter 3).

 Pour cela, on considère la réalisation standard de $\Sp(2n,\bbR)$ et 
l'épinglage standard    $\mathbf{spl}_{\Sp(2n,\bbR)}=(\mathbf B_d,\mathbf T_d, \{X_\alpha\}_{\alpha\in \Delta})$.
 En particulier, on a un sous-groupe de Borel $B_d=\mathbf{B}_d(\bbR)$  de $\Sp(2n,\bbR)$, son radical unipotent $N$, son algèbre de Lie $\frn$,
 et une base  $(X_\alpha)_{\alpha\in \Delta} $    de $\frn / [\frn,\frn]$ donnée
par les vecteurs radiciels  pour les racines simples venant de l'épinglage.
Pour avoir une donnée de Whittaker, il suffit alors de choisir un caractère unitaire additif non trivial   $\psi_\bbR$ de $\bbR$
et l'on définit un caractère unitaire  non dégénéré $\chi_N$ de $N$ par 
\begin{equation}\label{chiN} 
 \chi_N \left(\exp \left(\sum_{\alpha\in R(B,T)}x_\alpha   \, X_\alpha \right) \right)=   \psi_\bbR\left(\sum_{\alpha\in \Delta} x_\alpha\right)   \end{equation}

On pose 
\begin{equation}\label{caraddR}
\psi_{\bbR,1} : x \mapsto \exp (2i\pi x) . \end{equation}
C'est un caractère unitaire additif  de $\bbR$ et $\psi_{\bbR,-1}=\psi_{\bbR,1}^{-1}$.
On note  $\mathrm{Wh}_{\pm 1}$ la classe de conjugaison de la donnée de Whittaker $(N,\chi_N)$
définie comme ci-dessus avec $\psi_\bbR=\psi_{\bbR,\pm 1}$.

Il y a deux  classes de conjugaison de données de Whittaker pour le groupe $\Sp(2n,\bbR)$, qui sont 
$\mathrm{Wh}_{\pm 1}$.

Pour un groupe $\mathbf G$ défini sur $\bbR$ et quasi-déployé admettant des séries discrètes, les classes de conjugaison
de données de Whittaker sont en bijection avec les classes de conjugaison de paires de Borel fondamentales de 
type Whittaker $(\mathbf B,\mathbf T)$. Rappelons que cela signifie que $\mathbf T$ est un tore maximal de $\mathbf G$
défini sur $\bbR$ tel que $T=\mathbf T(\bbR)$ soit compact, et $\mathbf B$ est un sous-groupe de Borel contenant 
$\mathbf T$. La paire de Borel $(\mathbf B,\mathbf T)$ est alors fondamentale et une telle paire détermine
par la paramétrisation d'Harish-Chandra une unique série discrète  $\pi_{(\mathbf B,\mathbf T)}$ de $G$ de
 même caractère infinitésimal que la représentation triviale.
 De plus, la condition  d'être \og de type Whittaker\fg \,   signifie que 
 les racines simples de $\mathbf T$ dans $\mathbf B$ sont toutes imaginaires non compacte et la série discrète 
 $\pi_{(\mathbf B,\mathbf T)}$ est alors générique. Elle admet donc un modèle de Whittaker pour une 
certaine donnée de Whittaker, et ceci réalise la bijection. Dans \cite{MR3}, nous avons utilisé ceci pour fixer les 
données de Whittaker sur les groupes classiques, sans rendre cette bijection explicite. Nous le faisons  maintenant ici pour les groupes symplectiques.

Dans l'article, nous avons fixé un sous-groupe de Cartan compact  $T$ de  $\Sp(2n,\bbR)$
et le système de racine $\Delta(\frt,\frg)=\{ \pm e_i\pm e_j, \, 1\leq i<j\leq n, \, \pm 2e_i,  \, 1\leq i \leq n \}$.

Il y a deux classes de conjugaison de  paires de Borel fondamentales de 
type Whittaker dans $\Sp(2n,\bbR)$, admettant comme représentants  $(\mathbf B_{\pm},\mathbf T)$, où
le sous-groupe de Borel $\mathbf B_{+}$ est celui dont les racines simples sont
$ \{e_1+e_2, -(e_2+e_3), (-1)^{n} (e_{n-1}+e_n), (-1)^{n-1}2e_n  \} $
et  $\mathbf B_{-}$ est celui dont les racines simples sont les opposées de celles-ci.

\begin{prop}
La série discrète générique  $\pi_{(\mathbf B_{+},\mathbf T)}$ de $\Sp(2n,\bbR)$ admet 
un modèle de Whittaker pour la donnée $Wh_{1}$ de ce groupe. 
La bijection entre classes de conjugaison de donnée de Whittaker et  classes de conjugaison de paires fondamentales de type Whittaker
est donc 
\[ \mathrm{Wh}_{ 1} \leftrightarrow (\mathbf B_{+},\mathbf T), \quad \mathrm{Wh}_{ -1} \leftrightarrow (\mathbf B_{-},\mathbf T), .  \]
\end{prop}

\dem
Pour $n=1$, c'est un calcul de Wallach (\cite{Wallach}) et pour $n=2$, d'Oda (\cite{Oda}).
Esquissons un argument pour se ramener à l'un de ces deux cas, suivant la parité de $n$. Nous laissons les détails au lecteur.
Supposons $n=2r+1$ impair. Un paramètre  d'Harish-Chandra 
$\underline \lambda$ pour la série discrète générique est  un élément entier régulier de $\frt^*$, positif pour les racines
simples de $\mathbf B_+$, c'est-à-dire  $\underline{ \lambda}= (\lambda_1, -\lambda_2, \lambda_3, \ldots , -\lambda_{2r}, \lambda_{2r+1})$ avec 
$\lambda_1>\lambda_2>\ldots >\lambda_{2r+1}>0$. Faisons dégénérer ce paramètre en un paramètre de limite de séries discrètes
en prenant $\lambda_1=\lambda_2>\lambda_3=\lambda_4>\ldots >\lambda_{2r-1}=\lambda_{2r}>\lambda_{2r+1})$.
Les racines simples  s'annulant sur ce paramètre étant imaginaires non compactes, il y a  bien une limite de séries discrètes
génériques $\pi(\underline{\lambda})$ associée à ce paramètre.  Considérons un sous-groupe de Cartan de 
$\Sp(2n,\bbR)$ isomorphe à $(\bbC^\times)^r\times \U(1)$. Sur le $i$-ème facteur $\bbC^\times$, considérons
le caractère de différentielle $(\lambda_{2i-1}, \lambda_{2i})$. 
Sur le facteur $\U(1)$, on considère le caractère de différentielle $\lambda_{2r+1}$.
Induisons cohomologiquement le produit de ces caractères vers le groupe $\GL_2(\bbR)^r\times \SL_2(\bbR)$, avec la version des foncteurs
d'induction cohomologique préservant le caractère infinitésimal.
On obtient une représentation tempérée  de $\GL_2(\bbR)^r\times \SL_2(\bbR)$ que l'on induit paraboliquement vers $\Sp(2n,\bbR)$.
On obtient donc une représentation standard de ce groupe admettant une fonctionnelle de Whittaker unique à un scalaire près  pour la donnée 
compatible avec celle de $\SL_2(\bbR)$, c'est-à-dire, comme $\lambda_{2r+1}>0$, la donnée $\mathrm{Wh}_1$ d'après le calcul de Wallach.
Or, d'après les résultats du chapitre 11 de \cite{KV} (voir \cite{Mat04}, Thm. 2.2.3  pour une référence commode), on peut obtenir cette représentation standard
en  partant du même caractère de  $(\bbC^\times)^r\times \U(1)$ et en appliquant d'abord une induction parabolique vers
$\U(1,1)^r\times \U(1)$, puis une induction cohomologique (dans le good range) vers $\Sp(2n,\bbR)$.  Sur chaque facteur $\U(1,1)$, on obtient une  série principale
qui est somme de deux limites de séries discrètes et qui s'écrivent respectivement  comme induite
cohomologique de $\U(1,0)\times \U(0,1)$ ou $\U(0,1)\times \U(1,0)$  vers $\U(1,1)$. Par transitivité de l'induction cohomologique, on voit que tous les facteurs
de composition de la représentation standard  obtenue sont des limites de séries discrètes, dont une seule est générique,
celle qui correspond au choix de $\U(1,0)\times \U(0,1)$ pour chaque facteur $\U(1,1)$,  et c'est  $\pi(\underline{\lambda})$.
Ceci montre  que la donnée de Whittaker attachée à $\pi(\underline{\lambda})$ et donc à $\mathbf B_+$ est bien  $\mathrm{Wh}_1$.
Dans le cas où $n$ est pair, on fait le même raisonnement pour se ramener à $\Sp(4,\bbR)$ et au résultat d'Oda.\qed

\begin{rmq} \label{rmswh} La donnée de Whittaker utilisée dans \cite{MR3}, section 9.2 est  $\mathrm{Wh}_{ 1}$ pour tous les groupes symplectiques.
 Ceci sera utilisé dans les propositions \ref{Apsipinm} et \ref{Apsisnm}.
\end{rmq}

\subsection{Décomposition des $A$-paramètres pour  $G=\Sp(2n,\bbR)$}
 On compose un $A$-paramètre $\psi$ comme ci-dessus pour $G=\Sp(2n,\bbR)$  avec la représentation standard de 
${}^LG=\SO(2n+1,\bbC)$ dans $\GL(2n+1,\bbC)$ et l'on note encore $\psi$ le morphisme obtenu, que l'on voit
comme une représentation de $W_\bbR\times \SL_2(\bbC)$. Cette représentation est complètement réductible. 
Pour tout $a\in \bbN^\times$, notons
$R[a]$ la représentation algébrique de $\SL_2(\bbC)$ de dimension $a$, et pour tout $t\in \bbN^\times$, notons
$\delta_t$ la représentation irréductible de $W_\bbR$ de dimension $2$ qui est le paramètre de Langlands de la série
discrète de $\GL_2(\bbR)$ de caractère infinitésimal $(t/2,-t/2)$ (elle est notée $V(0,t)$ dans \cite{MR3}).
La forme générale de la  décomposition de $\psi$ en irréductibles est écrite dans  \cite{MR3}, \S 4.1.
Ici, on ne considère que des paquets ayant des caractères infinitésimaux entiers, et  ceux-ci
se décomposent de la manière suivante :
\begin{equation} \label{decApar}  \psi= \psi_u\oplus \psi_d=\left( \bigoplus_{i=1}^r   \eta_i \boxtimes R[a'_i] \right)
 \oplus \left( \bigoplus_{j=1}^s   \delta_{t_j} \boxtimes R[a_j]\right).
 \end{equation}
Dans la première   somme qui constitue la partie unipotente 
$\psi_u$ du paramètre,  $\eta_i$ désigne un caractère quadratique de $W_\bbR$ (le caractère trivial
que nous notons $1_{W_\bbR}$  ou  le caractère signe que nous notons $\sgn_{W_\bbR}$) et les $a'_i$ sont impairs.
Dans la deuxième somme, qui constitue la partie discrète $\psi_d$ du paramètre,
  les $t_j$ sont dans $\bbN\setminus\{0\}$ et  $t_j+a_j$ est impair.
On a de plus 
\begin{equation*} \label{sum} \sum_{i=1}^r  a'_i+2\sum_{j=1}^s a_j=2n+1, 
\; \text{ et } \;   
  \prod_{i=1}^r  \eta_i=\sgn_{W_\bbR}^{\vert \{ j; a_j \text{ impair } \}  \vert }. \end{equation*}

Le caractère infinitésimal  ({\sl cf. Convention \ref{conv}})   des éléments de $\Pi(\psi)$ est obtenu en rangeant dans l'ordre décroissant les
éléments des ensembles (avec multiplicités) d'entiers
\[  \left\{  \frac{a'_i-1}{2},  \frac{a'_i-3}{2}, \ldots, -\frac{a'_i-1}{2}  \right\}_i , \quad
 \left\{  \frac{t_j+a_j-1}{2},  \frac{t_j+a_j-3}{2}, \ldots, \frac{t_j-a_j+1}{2}  \right\}_j , \]   
\[ \left\{  \frac{-t_j+a_j-1}{2},  \frac{-t_j+a_j-3}{2}, \ldots, \frac{-t_j-a_j+1}{2}  \right\}_j . \]

\begin{lemme}\label{longpsiu} Soit $\psi$ un $A$-paramètre pour $G=\Sp(2n,\bbR)$ qui se décompose comme en (\ref{decApar}).
Supposons que le caractère infinitésimal du $A$-paramètre $\psi$ soit celui d'un module unitaire de plus haut poids
$\pi(\mu)$ du théorème \ref{hw}. Alors : 

$(i)$ la longueur de $\psi_u$ comme représentation de $W_\bbR\times \SL_2(\bbC)$
est $1$ ou $3$, et si elle est de longueur $3$, l'une des composante est de dimension $1$; 

$(ii)$ s'il existe un indice $j\in \{1,\ldots, s\}$ tel que $t_j-a_j+1\leq 0$, alors il est unique et   $\psi_u$ est irréductible, et même de dimension 
$1$ si $t_j-a_j+1< 0$.
\end{lemme}
\dem La formule (\ref{sum}) montre que $\psi_u$ est non nul pour une question de parité, et donc de longueur au moins égale à $1$.
Nous avons vu dans la remarque \ref{carinfhw} que le caractère infinitésimal  est donné
par les entiers $m_1-1,m_2-2,\ldots,m_n-n$ de leur opposés et de $0$. En particulier, $0$ apparaît avec multiplicité
$1$ ou $3$, et la longueur de $\psi_u$ est donc au plus $3$. Mais elle ne peut être $2$ car la contribution de la 
partie discrète à la multiplicité de $0$ dans le caractère infinitésimal est paire (cette contribution provient des indices $j$
dans la partie discrète tels que $t_j-a_j+1\leq 0$, ce qui au passage prouve la première assertion de $(ii)$).
D'autre part, la multiplicité des entiers non nuls dans le caractère infinitésimal est paire. 
Si $\psi_u$ est de longueur $3$ et si les trois composantes ont une dimension strictement supérieure à $1$, alors 
$1$ doit apparaître avec un multiplicité au moins égale à 3, ce qui n'est pas vrai. Ainsi si $\psi_u$ est de longueur $3$ l'une des composantes 
est de  dimension $1$. Pour la seconde assertion de $(ii)$, c'est le même argument,  la multiplicité de $1$  est au moins égale à $3$  si $\psi_u>1$
 et s'il existe  un indice $j$ tel que $t_j-a_j+1<0$, ce qui n'est pas possible.
\qed

\begin{defi}  \label{apsi} Soit $\psi$ un A-paramètre pour le groupe $G=\Sp(2n,\bbR)$, avec 
\begin{equation*}   \psi=\psi_u\oplus \psi_d= \left( \bigoplus_{i=1}^r   \eta_i \boxtimes R[a'_i] \right)
 \oplus  \left( \bigoplus_{j=1}^s   \delta_{t_j} \boxtimes R[a_j] \right).
 \end{equation*}
On définit alors :
\[ a(\psi)=\max(a'_i,\,  i=1,\ldots, r, \, a_j, \, j=1,\ldots,s) \quad \text{  et  } \quad a(\psi_u)=\max( a'_i, \, i=1,\ldots,r).\]
Ainsi $a(\psi)$  (resp. $a(\psi_u)$) est la plus grande dimension d'une représentation de $\SL_2(\bbC)$ intervenant dans le paramètre $\psi$, et 
(resp. dans la partie unipotente $\psi_u$).
\end{defi}
Ces entiers associés au paramètre $\psi$ vont jouer un grand rôle dans la suite.

\

\begin{convention} \label{ordret} (pour l'ordre des indices dans la partie discrète).  Soit $\psi$ un $A$-paramètre pour $\Sp(2n,\bbR)$ comme en (\ref{decApar}).
On range les indices de la partie discrète de telle sorte que la suite $(t_j)_{j=1,\ldots,s}$ soit décroissante, et si $t_{j}=t_{j+1}$, 
alors $a_j\geq a_{j+1}$.\end{convention}

\subsection{Description des paquets d'Arthur pour $\Sp(2n,\bbR)$ d'après \cite{MR3}}\label{ArtpacMR}

Soit $\psi$ un $A$-paramètre  pour $G=\Sp(2n,\bbR)$, que l'on écrit comme dans (\ref{decApar}).
Dans \cite{MR3}, les représentations dans le paquet $\Pi(\psi)$ sont décrites comme les composantes
irréductibles de représentations obtenues par induction cohomologique à partir de représentations unipotentes
de $c$-Levi $L$ de $G$. Rappelons brièvement ceci.
Il n'est pas vrai en général que $\psi_u$ soit un $A$-paramètre pour le groupe $G_u=\Sp(2n_u,\bbR)$ où 
$n_u=\frac{\dim \psi_u-1}{2}$. En effet, ce paramètre est à valeur dans $\Or(2n_u+1)$ et non nécessairement dans 
 $\SO(2n_u+1)$. Il faut donc le corriger par le caractère $\sgn_{W_\bbR}^{\dim(\psi_d)/2}$.
Alors 
\[  \psi'_u=\sgn_{W_\bbR}^{\dim(\psi_d)/2}\otimes \psi_u\]
est un $A$-paramètre pour $G_u=\Sp(2n_u,\bbR)$.
On a donc un paquet $\Pi(\psi'_u)$ associé, et l'on sait qu'il est constitué de représentations faiblement unipotentes de $G_u$
au sens de \cite{KnVo} et que leur multiplicité dans le paquet est $1$ (voir \cite{pourhowe} et \cite{MR5}).

D'autre part, pour tout $j$ indexant un terme de la partie discrète $\psi_d$ de $\psi$, donnons nous un entier 
$c_j$ entre $0$ et $a_j$ et notons 
\begin{equation}\label{csou} \scrC=\{\underline c=(c_1,\ldots, c_s)\}   \end{equation}
l'ensemble des familles d'entiers ainsi obtenues. Chaque $\underline c\in \scrC$ détermine à conjugaison près
une paire    $(L_{\underline c},\frqqq_{\underline c})$ où $L_{\underline c}$ est un 
 $c$-Levi de  $G$, isomorphe à 
\begin{equation}\label{Lc} G_u \times \left( \times _{j=1}^s \U(c_j,a_j-c_j) \right) \end{equation}
et $\frqqq_{\underline c}$ une sous-algèbre parabolique $\theta$-stable de $\frg$.
Ceci est clair par une récurrence immédiate  vu  la description des paires paraboliques maximales en (\ref{qpq}).

Les représentations du paquet $\Pi(\psi)$ sont alors les composantes irréductibles des induites cohomologiques
\[  \caR_{L_{\underline c},\frqqq_{\underline c}}^S(\pi_u\otimes \Lambda)  \]
où $\pi_u$ est une représentation unipotente dans le paquet $\Pi(\psi'_u)$ et $\Lambda$ un caractère du groupe 
$\times _{j=1}^s \U(c_j,a_j-c_j)$ uniquement déterminé par son caractère infinitésimal, et que l'on calcule aisément 
à partir des formules (\ref{lambda1}) et (\ref{y}) par récurrence. La condition pour que l'induction 
parabolique soit dans le weakly fair range est $t_1\geq t_2\geq \cdots \geq t_s$, ce que l'on a supposé.

\section{Paramètres unipotents}\label{parunip}\label{ParUnip}

On considère un $A$-paramètre $\psi$ unipotent pour $G=\Sp(2n,\bbR)$, c'est-à-dire $\psi=\psi_u$, 
ayant le caractère infinitésimal d'un module holomorphe unitaire de $G$. D'après le lemme \ref{longpsiu}, la longueur de 
$\psi$ comme représentation de $W_\bbR\times \SL_2(\bbC)$ est 1 ou 3. 
 Si cette longueur est $1$, c'est-à-dire si $\psi$ est irréductible, on a 
 \[  \psi=1_{W_\bbR} \boxtimes R[2n+1] \]
 et il est bien connu que le paquet $\Pi(\psi)$ est réduit  à la représentation triviale de $G$. 
 La représentation triviale est bien sûr un module unitaire holomorphe, avec les notations du théorème \ref{hw}, 
 c'est le module $\pi(0,\ldots 0)$. 
 
 Supposons donc maintenant que la longueur de $\psi$ est 3. On a alors 
 \[  \psi= (\eta_1 \boxtimes  R[a])\oplus (\eta_2\boxtimes R[b])\oplus (\eta_3\boxtimes R[1]) \]
 avec $a$ et $b$ impairs, $a+b=2n$ et $\eta_1\eta_2\eta_3=1_{W_\bbR}$.  On suppose que $a\geq b\geq 1$ et donc 
 en particulier  $b\leq n$. 
 D'autre part, on suppose que $\eta_1=\eta_2\eta_3=\sgn_{W_\bbR}^{\frac{b+1}{2}}$, c'est-à-dire :
 \begin{equation}\label{psiar}
\psi= (\sgn_{W_\bbR}^{\frac{b+1}{2}} \boxtimes  R[a])\oplus (\eta_3\,  \sgn_{W_\bbR}^{\frac{b+1}{2}}\boxtimes R[b])\oplus (\eta_3\boxtimes R[1]) .
  \end{equation}

Soit  $\pi( \mu)$ un module holomorphe unitaire 
 comme dans le théorème \ref{hw}  et
 supposons que  $\pi(\mu)$ soit  l'image par  
 la correspondance de Howe d'une représentation de dimension finie du groupe compact
 $\Or(0,b+1)$ dont les composantes de la restriction à  $\SO(0,b+1)$ sont  contenues dans le paquet  
 $\Pi (\psi')$ où $\psi'=    \eta_2\boxtimes R[b]\oplus \eta_3\boxtimes R[1]$
 ($\psi'$ est bien un paramètre pour ce groupe, grâce à la condition $\eta_2\eta_3=\sgn_{W_\bbR}^{\frac{b+1}{2}}$, car  $\SO(0,b+1)$
 est forme intérieure du groupe déployé $\SO(\frac{b+1}{2},\frac{b+1}{2})$ si $\frac{b+1}{2}$ est pair, et
 du groupe  quasi-déployé non déployé $\SO(\frac{b+1}{2}-1,\frac{b+1}{2}+1)$ si $\frac{b+1}{2}$ est impair).

 Le caractère infinitésimal de cette représentation de dimension finie est  $\left(  \frac{b-1}{2},  \frac{b-3}{2},\ldots,  0\right)$. 
 C'est donc soit la représentation triviale   $\Triv_{\Or(0,b+1)}$, soit le caractère  quadratique de $\Or(0,b+1)$  donné par le déterminant
 (et que nous notons  simplement $\det_{\Or(0,b+1)}$).
 D'après le théorème \ref{KV}, par inspection, on est dans l'un des cas suivants. 
 Soit $\mu=\left(\frac{b+1}{2},\ldots,\frac{b+1}{2}\right)$ et $\pi(\mu)$ est l'image par la correspondence de Howe de 
 $\Triv_{\Or(0,b+1)}$, soit $\mu=\left( \underbrace{ \frac{b+3}{2},\ldots,\frac{b+3}{2}}_{b+1},\underbrace{\frac{b+1}{2},\ldots,\frac{b+1}{2}}_{n-b-1} \right)$
 et $\pi(\mu)$ est l'image par la correspondence de Howe de $\det_{\Or(0,b+1)}$  et ceci nécessite $n\geq b+1$.
 Dans ce cas, notons cette représentation
 \begin{equation} \label{defsigmana}
 \sigma_{n,\frac{b+1}{2}}=\pi(\mu), \quad \text { où } \mu=
 \left( \underbrace{ \frac{b+3}{2},\ldots,\frac{b+3}{2}}_{b+1},\underbrace{\frac{b+1}{2},\ldots,\frac{b+1}{2}}_{n-b-1} \right).
 \end{equation}

\begin{thm}\label{thmunip}
Soit $\psi$ un $A$-paramètre  unipotent pour $G=\Sp(2n,\bbR)$. 
 
 --- Si $\psi$ est irréductible, c'est-à-dire $\psi=1_{W_\bbR}\boxtimes R[n]$, alors 
  le paquet $\Pi(\psi)$ est constitué de la représentation triviale de $G$.

 --- Si $\psi$ est de longueur $3$, on suppose que $\psi$ s'écrit :
\begin{equation}\label{psiar2}
\psi= (\sgn_{W_\bbR}^{\frac{b+1}{2}} \boxtimes  R[a])\oplus ( \eta_3\,  \sgn_{W_\bbR}^{\frac{b+1}{2}}\boxtimes R[b])\oplus (\eta_3\boxtimes R[1] )
= (\sgn_{W_\bbR}^{\frac{b+1}{2}} \boxtimes  R[a]) \oplus \psi',
  \end{equation}
avec $a\geq b\geq 1$.
Soit $\pi(\mu)$ un module unitaire holomorphe comme dans le théorème \ref{KV}. 
On suppose d'autre part que  $\psi'=(\eta_3\,  \sgn_{W_\bbR}^{\frac{b+1}{2}}\boxtimes R[b])\oplus (\eta_3\boxtimes R[1] )$ 
est un $A$-paramètre pour $\Or(0,b+1)$ tel que le paquet $\Pi(\psi')$ contienne une représentation de dimension finie $E_{\psi'}$
et que $\pi(\mu)$ est l'image par la correspondance de Howe  pour la paire $(\Or(0,b+1),\Sp(2n,\bbR))$ de cette représentation $E_{\psi'}$. 
Alors soit $E_{\psi'}$ est la representation triviale de    $\Or(0,b+1)$, et alors $\pi(\mu)=\pi_n(\frac{b+1}{2})$, 
 soit   $E_{\psi'}$ est le déterminant  de    $\Or(0,b+1)$, et  alors $b+1\leq n$ et $\pi(\mu)=\sigma_{n,\frac{b+1}{2}}$
 ({\sl cf.} Notation \ref{pimn} et (\ref{defsigmana})). 
De plus  $\Pi(\psi)$ contient $\pi(\mu)$.
\end{thm}

Nous avons déjà  démontré les premières assertions avant l'énoncé du théorème. Il reste à voir  dans le deuxième cas que  $\Pi(\psi)$ contient $\pi(\mu)$.
Nous démontrerons un résultat plus général plus loin, le théorème \ref{lemp} et nous renvoyons le lecteur
à sa démonstration.

\section{Cas du caractère infinitésimal régulier}\label{carinfreg}

Soit $\underline \chi=(\chi_1,\ldots,\chi_n)\in \frt^*$ définissant un caractère infinitésimal régulier 
de $\Sp(2n,\bbR)$, c'est-à-dire que 
 la suite  $\nu_1\geq \nu_2\geq \ldots \geq \nu_{2n+1}$
 formée avec les entiers $\chi_1, \ldots, \chi_n$, leurs opposés, et $0$ est sans multiplicité.
 On peut supposer que $\chi_1>\chi_2>\ldots>\chi_n>0$.

On suppose qu'il existe des modules unitaires holomorphes dont le caractère infinitésimal est 
$\underline \chi \in \frt^*$. Comme dans la section \ref{ss:modphp}, notons  
$\mu^*= (-m_n\geq -m_{n-1}\geq \cdots \geq -m_1)$ 
le plus haut poids d'un tel module. Alors $ (m_1-1, m_2-2, \ldots, m_n-n) $
définit le même caractère infinitésimal que $\underline \chi$.

On est alors nécessairement dans le cas $(a)$ ou $(b)'$ du théorème \ref{KV} et de la remarque \ref{rmqKV}. De plus dans le cas (b'),
en posant $\ell=m_n$, on a $m_\ell-\ell>a=n-\ell$   {i.e.} $m_\ell>n$;  en  particulier, $m_\ell>m_n=\ell$ et  $a=u$.
Dans le cas $(a)$, on pose $a=0$.  Le caractère infinitésimal étant fixé, 
le module unitaire holomorphe en question, est entièrement déterminé par l'entier $a$ 
(où de manière équivalente l'entier $\ell$ avec $\ell=n-a)$ et on le note $\pi_a$, où
l'entier $a$ varie de $0$ à $a_{max}$,  qui est obtenu comme suit
si $\chi_1\neq 1$, on pose $a_{max}=0$, et si $\chi_1=1$, soit $a_{max}$ l'entier strictement  positif 
tel que $a_{max}\leq n$, $ (\chi_{n-a_{max}+1},\ldots,\chi_n)=(a_{max},a_{max-1},\ldots, 1)$ et 
si $a_{max}\neq n$, $\chi_{n-a_{max}}>a_{max}+1$.

Le module unitaire $\pi_a$  est le plus haut poids $\mu_a = (-m_n,\ldots ,-m_1)$, avec 
$ (m_1, m_2, \ldots, m_n)=(m_1,\ldots, m_\ell,\ell,\ldots,\ell)$,  
 et 
$  (\chi_1,\ldots,\chi_n)  = (m_1-1,\ldots m_\ell-\ell, a ,a-1,\ldots,1)$. 
Reprenons les notations de la section \ref{Sec2}. Notons $L$ le $c$-Levi de $\Sp(2n,\bbR)$ avec 
\[L\simeq \Sp(2a,\bbR)\times  \underbrace{\U(1,0)\times \cdots \times \U(1,0)}_\ell.\]
On a alors 
\[ \rho=(-1,-2,\ldots, -n),  \quad \rho(\frk)= \left(\frac{n-1}{2},\ldots ,-\frac{n-1}{2}\right), \quad \rho(\fru)=( 0,\ldots,0,-n+\ell-1,\ldots, -n),\]
 \[ \rho(\fru\cap \frp)
  =\left(-\frac{\ell}{2}, \ldots ,-\frac{\ell}{2},-\frac{n+1}{2},\ldots ,-\frac{n+1}{2} \right).\]
On définit un caractère unitaire $\Lambda$ de $L$ en donnant sa différentielle, qui le détermine car $L$ est connexe :
\begin{align*}
  \lambda=& \mu_a-2\rho(\fru\cap \frp) =(0,\ldots ,0, -m_\ell+n+1, \ldots ,-m_1+n+1).
\end{align*}
Remarquons que l'on a bien $\lambda+\rho=(-1,-2,\ldots ,-a,  -m_\ell+n+1-(n-l+1), \ldots ,-m_1+n+1-n)=(-1,-2,\ldots ,-a,  -m_\ell+\ell, \ldots ,-m_1+1)$.

\begin{prop}
On a $\pi_a=A_\frqqq(\lambda)$.
\end{prop}
\dem C'est une conséquence directe de \cite{VZ}, Prop. 6.1.
\qed

\bigskip 

Considérons maintenant un paramètre d'Arthur
\[\psi : W_\bbR\times \SL_2(\bbC) \longrightarrow \SO(2n+1,\bbC) \]
de caractère infinitésimal $\underline \chi$, que l'on décompose
en somme directe de représentations irréductibles de $ W_\bbR\times \SL_2(\bbC)$ :
\[\psi=\oplus \psi_i .\]
Soit $\psi_u$ la partie unipotente de $\psi$. Alors $\psi_u$ est irréductible, c'est-à-dire que 
$\psi_u$ est l'un des $\psi_i$ (ceci à cause de l'hypothèse de régularité de $\underline \chi$).
De plus, $\psi_u$ est de la forme
\[\psi_u=\epsilon\otimes R[2a(\psi_u)+1],\]
 où $\epsilon$ est un caractère quadratique de $W_\bbR$
et $a(\psi_u)$ un entier.

\begin{thm}\label{thmcasreg}
La représentation $\pi_a$ est dans le paquet $\Pi(\psi)$ si et seulement si $a=a(\psi_u)$.
\end{thm}
\dem Les paquets d'Arthur de caractère infinitésimal entier régulier $\underline \chi$ sont des paquets d'Adams-Johnson
(\cite{AdJo}, \cite{AMR}), que l'on connait exactement.

\section{ Le cas où le plus haut poids est scalaire}\label{scal}
Soit $\pi(\mu)$ un module unitaire holomorphe comme dans le  théorème \ref{hw} et  
dans cette section, on suppose que le poids  $\mu=(m_1\geq m_2\geq .... \geq m_n)$ est de la forme 
$\mu=(m,\ldots,m)$, avec $0\leq  m \leq n$ (dans le cas $m>n$, $\pi(\mu)$
est une série discrète holomorphe, et les paquets d'Arthur la contenant sont des paquets d'Adams-Johnson décrit dans \cite{AMR}).
   Rappelons que nous avons noté  cette représentation $\pi_n(m)$ ou simplement 
$\pi(m)$.  Avec les notations du théorème \ref{hw} et du théorème \ref{KV}, 
 on a donc $u=n$ et l'on  est  dans un des cas (b'') ou (d) de la remarque \ref{rmqKV}. 

\medskip 

Si $m \leq \frac{n+1}{2}$, on est dans le cas $(b'')$, et $\pi(m)$ est image par la correspondance de Howe de la représentation 
triviale de $\Or(0,2m)$.
Si $m  \geq \frac{n}{2}+1$,  c'est-à-dire $2(m-1)\geq n$, on est dans le   cas  (d). 
C'est un sous-cas de $(b)$ et donc comme ci-dessus  $\pi(m)$ est image par la correspondance de Howe 
de la représentation 
triviale de $\Or(0,2m)$. Mais c'est aussi  l'image par la correspondence de Howe  de la représentation 
$[1,\ldots,1,0,\ldots,0]_-$ de  $\Or(0,2(m-1))$, où il y a $2(m-1)-n$ fois $0$ et $n-(m-1)$ fois $1$.

Enonçons notre  résultat pour les modules holomorphes de plus haut poids scalaire $\pi_n(m)$.

\begin{thm} \label{mainthm}  Soit $\psi$ un $A$-paramètre pour le groupe $G=\Sp(2n,\bbR)$ se décomposant en 
\[\psi=\psi_u\oplus \psi_d=\psi_u\oplus  \bigoplus_{i=1}^s   (\delta_{t_j}\boxtimes R[a_j] )\]
 où  $\psi_u$ est la partie unipotente de  $\psi$,  et $\psi_d $ sa partie 
discrète. Soit $m$ un entier, avec 
$0\leq m\leq n$. On suppose que  le caractère infinitésimal
associé à ce paramètre est  celui de $\pi_n(m)$, c'est-à-dire constitué des entiers 
$m-1,m-2,\ldots, m-n$, de leur opposés, et de $0$  ({\sl cf.}  Convention \ref{conv}).
 Alors le paquet $\Pi(\psi)$ contient $\pi_n(m)$ si et seulement si l'on est dans un des cas suivants: 

$(i)$ $\dim \psi_u=1$, $2m >n+1$, et quels que soient $i<j$ entre $1$ et $s$, 
\begin{equation}\label{disj}  \left[ \frac{t_i-a_i+1}{2},  \frac{t_i+a_i-1}{2}\right]\cap \left[ \frac{t_j-a_j+1}{2},  \frac{t_j+a_j-1}{2}\right]=\emptyset.\end{equation}

$(ii)$  $\psi$ s'écrit 
 \[ \psi=(\sgn_{W_\bbR}^{\frac{2n+1-a(\psi_u)}{2}})  \boxtimes R[a(\psi_u)]\oplus \psi' \] 
où $\psi'$ est un $A$-paramètre pour le groupe $\Or(2n+1-a(\psi_u))$ tel que le paquet 
$\Pi(\psi')$ contienne une representation de dimension  finie $E_{\psi'}$
et $\pi_n(m)$ est image de $E_{\psi'}$ par la correspondance de Howe.

Dans ces deux cas, la multiplicité de $\pi_n(m)$ dans le paquet $\Pi(\psi)$ est $1$.
\end{thm}

Nous allons compléter ce résultat par le suivant.
\begin{thm}\label{complement}
Supposons que  le couple $(\psi,m)$ vérifie les hypothèses de $(i)$, alors  avec les notations de la section \ref{ArtpacMR},   
$\pi_n(m)$ est égale à  l'induite cohomologique  $\caR^S_{L_{\underline c}, \frqqq_{\underline c}} (\Lambda)$
 obtenue en partant d'une paire $(\frqqq_{\underline c}, L_{\underline c})$, où   $\underline c=(0,\ldots, 0)$, 
de sorte que $\frqqq_{\underline c}$ est  
holomorphe ({\sl cf.} Définition \ref{qholo}),  $L_{\underline c}$ isomorphe à $ \times_{j=1}^k \U(0,a_j) $
({\sl cf.} (\ref{csou}) et (\ref{Lc})) et $\Lambda$ est le produit tensoriel de caractères $\Lambda_j$
du facteur isomorphe à $\U(0,a_j)$ de différentielle 
\begin{equation}\label{jlamb}
\lambda_j = \left(\underbrace{ n-\sum_{k<j} a_k   -\frac{t_j-a_j+1}{2}, \ldots  , n-\sum_{k<j} a_k   -\frac{t_j-a_j+1}{2}} _{a_j} \right).
\end{equation}
(L'induction cohomologique a lieu dans le weakly fair range. Ici, en plus de l'unitarité de l'induite, on a aussi son irréductibilité)

Supposons que  le couple $(\psi,m)$ vérifie les hypothèses de  $(ii)$. Alors  
$a(\psi_u)=2(n-m)+1$ ou bien $a(\psi_u)=2(n-m)+3$, cette dernière possibilité
ne pouvant avoir lieu que si $2m\geq n+2$. Si $a(\psi_u)=2(n-m)+1$  (resp.  $2(n-m)+3$),   
alors $\pi_n(m)$ est  l'image par la  correspondance de Howe de la représentation triviale (resp. de la représentation de dimension finie 
 $[\underbrace{1,\ldots,1}_{2(m-1)-n},\underbrace{0,\ldots,0}_{n-(m-1)}]_-$)    de 
$\Or(0,2n+1-a(\psi_u))=\Or(0,2m)$    (resp. de $\Or(0,2n+1-a(\psi_u))=\Or(0, 2(m-1))$). 
 De plus  $\pi_n(m)$ est égale à  l'induite cohomologique  $\caR^S_{L_{\underline c}, \frqqq_{\underline c}} (\rho\boxtimes \Lambda)$
 obtenue en partant d'une paire $(\frqqq_{\underline c}, L_{\underline c})$, où   $\underline c=(0,\ldots, 0)$, 
de sorte que $\frqqq_{\underline c}$ est  
holomorphe ({\sl cf.} Définition \ref{qholo}),  $L_{\underline c}$ isomorphe à $\Sp(\dim(\psi_u)-1,\bbR)\times\left( \times_{j=1}^k \U(0,a_j) \right)$
({\sl cf.} (\ref{csou}) et (\ref{Lc})) et $\Lambda$ est le produit tensoriel de caractères $\Lambda_j$
de  $\U(0,a_j)$ de différentielle comme en (\ref{jlamb}).
 D'autre part, la   représentation    $\rho$ du facteur $\Sp(\dim(\psi_u)-1,\bbR)$ est la représentation unipotente
 $\pi_{n_u}(m_u)$   avec $n_u=\frac{\dim(\psi_u)-1}{2}$ et $m_u=\frac{\dim(\psi_u)-a(\psi_u)}{2}$,   
  image par la correspondance de Howe de 
la représentation triviale de $\Or(0,\dim(\psi_u)-a(\psi_u))$
  (resp. $\sigma_{n_u,k_u}$, avec $n_u=\frac{\dim(\psi_u)-1}{2}$ et $k_u=\frac{\dim(\psi_u)-a(\psi_u)}{2}$,   
  image par la correspondance de Howe du 
déterminant de $\Or(0,\dim(\psi_u)-a(\psi_u))$).

\end{thm}

\begin{rmq}
La condition $2m>n+1$ du $(i)$ du théorème est la même que celle du cas où $a(\psi_u)=2(n-m)+3$ du $(ii)$. 
Lorsque $2m\leq n+1$, il n'y a donc que les paramètres du $(ii)$ avec $a(\psi_u)=2(n-m)+1$.
\end{rmq}

Démontrons tout de suite une   partie de ce théorème, que l'on met à part sous la forme du lemme suivant. 
\begin{lemme}  \label{rmqapsiu}  Si le couple $(\psi,m)$ vérifie la condition $(ii)$ du théorème \ref{mainthm}, on a alors 
$a(\psi_u)=2(n-m)+1$ ou bien $a(\psi_u)=2(n-m)+3$, cette dernière possibilité
ne pouvant avoir lieu que si $2m\geq n+2$. Si $a(\psi_u)=2(n-m)+1$,    alors $\pi_n(m)$ est l'image de la représentation triviale 
 de $\Or(0,2n+1-a(\psi_u))$ par  la correspondance de Howe, et si $a(\psi_u)=2(n-m)+3$, c'est l'image de la représentation de dimension finie 
 $[\underbrace{1,\ldots,1}_{2(m-1)-n},\underbrace{0,\ldots,0}_{n-(m-1)}]_-$.
\end{lemme}

\dem Supposons que 
$\pi_n(m)$ soit  l'image par la correspondance de Howe d'une représentation de dimension finie $[\alpha_1,\ldots,\alpha_x,\underbrace{0,\ldots,0}_{\ell-x}]_\epsilon$
du groupe  $\Or(0,2\ell)$ avec $\alpha_1\geq \alpha_2\geq \cdots\geq \alpha_x>0$, et $\epsilon=\pm$. On obtient si $\epsilon=+$ :
\[  (\alpha_1+\ell, \ldots,\alpha_x+\ell,\underbrace{ \ell,\ldots,\ell}_{n-x})=(m,\ldots,m), \]
d'où $x=0$,  $m=\ell$ et $a(\psi_u)=2(n-m)+1$.  Maintenant si $\epsilon=-$, on obtient  :
\[  (\alpha_1+\ell, \ldots,\alpha_x+\ell,\underbrace{\ell+1,\ldots,\ell+1}_{2\ell-2x},\underbrace{\ell,\ldots \ldots,\ell}_{n-2(\ell-x)-x})
=(m,\ldots,m). \]
Ceci entraine $m=\ell+1$, puis $\alpha_1=\ldots,\alpha_x=1$ et enfin  $n-2(\ell-x)-x=0$, soit $x=2(m-1)-n$.
On a alors $a(\psi_u)=2(n-m)+3$ et $2m=n+2+x\geq n+2$. Les considérations qui précèdent le théorème montrent que réciproquement, $\pi_n(m)$ est bien 
l'image par la correspondance de Howe de la représentation triviale ou déterminant (selon que $a(\psi_u)=2(n-m)+1$ ou $2(n-m)+3$)   du groupe $\Or(0,2n+1-a(\psi_u))$.
\qed

\bigskip 
La démonstration des théorèmes \ref{mainthm} et \ref{complement} va occuper les  sections \ref{reduc}  à \ref{Multun}.

\section{Les représentations unipotentes $\sigma_{n,k}$}\label{sigmana}

Dans cette section, nous considérons certains modules unitaires holomorphes particuliers, qui ont la propriété d'être de plus unipotents, 
et nous montrons leur appartenance à certains paquets d'Arthur. 

On fixe $n$ et un entier $k$ tel que $2\leq 2k \leq n$. 
Considérons la représentation 
\begin{equation} \label{sigma}
 \sigma_{n,k} = \pi (\underbrace{k+1, \cdots, k+1}_{2k}, \underbrace{k, \cdots, k}_{n-2k}).
\end{equation} 
Nous avons déjà rencontré cette représentation en (\ref{defsigmana}). C'est l'image par la correspondance de Howe de la représentation $\det_{\Or(0,2k)}$.
L'image par la correspondance de Howe de la représentation triviale de $\Or(0,2k)$ est-elle le module holomorphe unitaire 
 $\pi_n(k)$. 
Ces représentations ont été particulièrement étudiées par J. Adams dans \cite{Ad87}, qui les note $\pi^\pm_{2k,n}$, avec $\pi^-_{2k,n}=\sigma_{n,k} $
et  $\pi^+_{2k,n}=\pi_n(k) $.

  Introduisons le paramètre 
 \begin{equation} \label{psia}
\psi=(\sgn_{W_\bbR}^{k}\boxtimes R[2(n-k)+1])\oplus  (\delta_{k-1}\boxtimes R[k]).
\end{equation}
On constate que son caractère infinitésimal est celui de $\sigma_{n,k}$ et de $\pi_n(k)$.

\begin{prop}\label{propsna0}
Le paquet d'Arthur $\Pi(\psi)$ où $\psi$ est le paramètre (\ref{psia}) contient $\sigma_{n,k}$.
\end{prop}
 
 On sait d'après \cite{MR3} que  $\Pi(\psi)$ contient les composantes irréductibles d'une certaine  représentation construite par induction
  cohomologique à partir de la paire   $(\frqqq,L)$ où ici $\frqqq=\frqqq_{p,q}$ avec $(p,q)=(0,k)$
 ({\sl cf.}  section \ref{ppm}), de sorte que 
 $L$ est isomorphe à   $ \Sp(2(n-k),\bbR)\times \U(0,k)$. 
 La représentation que l'on induit est un caractère de $L$ que l'on note $\Lambda$ : sur le facteur 
 isomorphe à $\Sp(2(n-k),\bbR)$, c'est donc le caractère trivial, et sur le facteur isomorphe à 
 $\U(0,k)$, il est donné par sa différentielle
 $(\underbrace{n-k+1,\ldots,n-k+1}_k)$.
 Notons cette  représentation induite  $\tilde{\sigma}_{n,k}$. C'est un cas particulier de la représentation considérée en (\ref{zalum}), avec 
   $y=-n+k-1$ dans (\ref{lambda1}) et $t=k-1$ dans (\ref{y}).
 C'est  donc un $A_\frqqq(\lambda)$ dans le weakly fair range : elle 
  n'est pas irréductible en général mais elle est unitaire et donc semi-simple. La proposition est conséquence immédiate du lemme suivant
  du à J. Adams (\cite{Ad87}, Prop. 5.1). 
  
  \begin{lemme}\label{sna}
 L'induite cohomologique   $\tilde{\sigma}_{n,k}$ se décompose en 
 $\tilde \sigma_{n,k}\simeq \sigma_{n,k}\oplus \pi_n(k) $.
\end{lemme}
    
    \begin{rmq}
    En fait Adams considère une induite à partir d'un $c$-Levi $L$ isomorphe à 
    \[ \Sp(2(n-k),\bbR)\times \underbrace{\U(0,1)\times \cdots \times \U(0,1)}_{k}.\]
 On utilise les résultats sur l'induction cohomologique  par étape en induisant d'abord de $\underbrace{\U(0,1)\times \cdots \times \U(0,1)}_{k} $ à $\U(0,k)$
 pour voir que les énoncés sont équivalents (voir le lemme suivant et sa démonstration). 
     \end{rmq}
     
  \begin{cor} \label{corsna} Considérons le paramètre 
 \begin{equation} \label{psia2}
\psi=(\sgn_{W_\bbR}^{k }\boxtimes R[2(n-k)+1] ) \oplus\bigoplus_{j=1}^s(\delta_{t_j}\boxtimes R[a_j])  .
\end{equation}
où $k=\sum_{j=1}^s a_j$, et supposons que le caractère infinitésimal de $\psi$ soit celui de $\sigma_{n,k}$. 
Alors la représentation $\sigma_{n,k}$ est dans le paquet $\Pi(\psi)$.
\end{cor}

\dem Remarquons tout d'abord que $\sgn_{W_\bbR}^k= \sgn_{W_\bbR}^{\vert \{j; \, 1\leq j\leq s, \,  a_j \text{impair } \}\vert  }$, 
et que le paramètre est bien à valeurs dans
 $\SO(2n+1,\bbC)$.
D'autre part, le caractère infinitésimal est aussi celui du paramètre (\ref{psia}), ce qui force 
 \begin{align*} 
 &( t_1+a_1-1, t_1+a_1-3 ,\dots , t_1-(a_1-1),  \\
& \qquad t_2+a_2-1, t_2+a_2-3 ,\dots,  t_2-(a_2-1), \ldots ,\\
&\qquad \qquad, \ldots,    t_s+a_s-1, t_s+a_s-3 ,\dots,  t_s-(a_s-1) )
= (k-1,\ldots , 1,0) \end{align*} 
Par induction cohomologique par étape, l'induite cohomologique $\tilde \sigma_{n,k}$ considérée dans le lemme précédent
est égale à l'une des induites cohomologiques associées au paramètre (\ref{psia2})
où l'on induit d'une paire parabolique $(\frqqq,L)$ avec $\frqqq$ holomorphe et 
$L$ isomorphe à 
\[\Sp(2(n-k),\bbR)\times\left(  \times_{j=1}^s \U(0,a_j) \right),\]  avec comme représentation induisante
le produit tensoriel de 
la représentation triviale sur le facteur $\Sp(2(n-k),\bbR)$ et des caractères $\Lambda_j$ sur les 
facteurs  $\U(0,a_j)$, avec pour différentielle 
\begin{equation}\label{lambdaj}\lambda_j=\big(\underbrace{n-\sum_{k<j}a_k-\frac{t_j+a_j-1}{2}, \ldots , n-\sum_{k<j}a_k-\frac{t_j+a_j-1}{2}}_{a_j}\big),
\end{equation}
de sorte que quand on induit cohomologiquement ce caractère $\boxtimes_{j=1}^s \Lambda_j$
de $\times_{j=1}^s \U(0,a_j)$ à $\U(0,k)$, on obtienne sur ce groupe le caractère 
de poids $(\underbrace{n-k+1,\ldots,n-k+1}_k)$ qui a servi dans la définition de  $\tilde \sigma_{n,k}$.
\qed 

\medskip 
Montrons maintenant que $\sigma_{n,k}$ est dans d'autres induites cohomologiques, et donc dans d'autres paquets.
Le résultat suivant se trouve aussi dans \cite{Ad87}, plus précisément dans la démonstration de la proposition 5.1, p.135 (Case 2). 

\begin{lemme}\label{lemmesna2}  Considérons la paire parabolique $(\frqqq,L)$ avec $\frqqq=\frqqq_{0,1}$ 
({\sl cf} section \ref{ppm});  en particulier,  $L$ est isomorphe à   $ \Sp(2n-2,\mathbb{R})\times \U(0,1)$. Supposons $2<2k\leq n$. 
Alors la représentation $\sigma_{n,k}$ est égale à   l'induite cohomologique  
 $\caR^S_{\frqqq,L,G}( \sigma_{n-1,k-1}\boxtimes \det^{n-k+1}) $ (cette induite est dans le weakly fair range, donc unitaire).
\end{lemme}

Là encore, on passe de l'énoncé du lemme à celui d'Adams par un argument simple d'induction par étape.
Ce que démontre Adams, c'est que  $\sigma_{n,k}$ est égale à   une induite cohomologique  construite avec 
$\sigma_{n-k+1,1}$.

\begin{prop}\label{propsna}
Soit $\psi=\psi_u\oplus \psi_d= \psi_u\oplus\bigoplus_{j=1}^s (\delta_{t_j}\boxtimes R[a_j])$  
 un $A$-paramètre pour $\Sp(2n,\bbR)$. Supposons que le caractère infinitésimal de $\psi$ soit celui
 de $\sigma_{n,k}$, et supposons de plus que $\psi_u$ contienne le facteur $\sgn_{W_\bbR}^{k}\boxtimes R[2(n-k)+1]$.
 Alors $\sigma_{n,k}$ est contenu dans le paquet $\Pi(\psi)$.
\end{prop}

\dem  Comme $2k\leq n$, on voit facilement que $a(\psi_u)=2(n-k)+1$.
 Ecrivons $\psi=\psi_u\oplus \psi_d= \psi_u\oplus\bigoplus_{j=1}^s(\delta_{t_j}\boxtimes R[a_j])$. 
Si $\psi_u$ est irréductible, on a $\psi_u=\sgn_{W_\bbR}^{k}\boxtimes R[2(n-k)+1] $.
Le lemme \ref{sna} et le corollaire \ref{corsna} donnent alors la conclusion voulue.
On suppose donc que $\psi_u$ est de longueur $3$.
 Posons $\dim \psi_u=2 n_u+1$ de sorte que  $\psi_u$, éventuellement tensorisé par le caractère $\sgn_{W_\bbR}$,  est un $A$-paramètre pour $\Sp(2n_u;\bbR)$.
On sait d'après \cite{MR3} que  $\Pi(\psi)$ contient les composantes irréductibles de la représentation construite par induction
  cohomologique à partir de la paire 
 $(\frqqq,L)$
 ({\sl cf.} section \ref{ArtpacMR}) où $\frqqq$ est holomorphe et 
 $L$ est isomorphe à   $\Sp(2n_u,\bbR)\times \left( \times_{j=1} ^s \U(0,a_j)\right)$ comme suit :
sur le facteur $\Sp(2n_u,\bbR)$,  la représentation que l'on induit  est 
  la représentation unipotente $\sigma_{n_u,a_u}$ où $a_u=k-\sum_ja_j$, de sorte que $n-k=n_u-a_u$ et 
  sur  les facteurs $\U(a_j,0)$, la représentation que l'on induit est 
un caractère de différentielle comme en (\ref{lambdaj}).
Appelons $\rho$ cette induite cohomologique. L'induction étant dans le weakly fair range, $\rho$ est unitaire, mais pas nécessairement
 irréductible. On veut montrer que $\rho$ contient $\sigma_{n,k}$.

L'égalité du  caractère infinitésimal du  paramètre $\psi$ avec 
celui de $\sigma_{n,k} $ implique  que pour tout $j$ avec  $1<j\leq s$, $\frac{t_{j-1}-(a_{j-1}-1)} {2}=1+\frac{ t_j+(a_j-1)}{2}$.
En utilisant des argument d'induction par étape comme dans la démonstration du corollaire \ref{corsna}, on voit que l'on peut remplacer
l'induite cohomologique par une induite cohomologique à partir d'une paire 
$(\frqqq',L')$ où $\frqqq'$ est holomorphe et $L'$ est isomorphe à $\Sp(2n_u, \bbR)\times \left(\times_{j' }\U(0,b_{j'})\right)$
du moment que  $\sum_j a_j=\sum_{j' }b_{j'}$, en adaptant le caractère du produit de groupe  unitaires que l'on induit, bien évidemment.
  En particulier, $\rho$ est aussi obtenue comme induite cohomologique avec la paire parabolique 
  $(\frqqq',L')$ où $\frqqq'$ est holomorphe et $L'$ est isomorphe à $\Sp(2n_u, \bbR)\times \U(0,\sum _{j}a_j)$.
  
On raisonne par récurrence sur $\sum_j a_j$  pour montrer que $\pi$ contient $\sigma_{n,k}$. 
 Par hypothèse de récurrence,  on peut supposer que l'induite cohomologique  construite  avec la paire parabolique 
 $\Sp(2n_u,\mathbb{R})\times \U(0,(\sum_j a_j)-1)$ contient $\sigma_{n-1,k-1}$,  car on amorce la récurrence avec le cas 
  $\sum_i a_i=1$ et le paramètre  unipotent  $\psi_u$, cas traité dans la section 
 \ref{parunip}.
Grâce au  lemme \ref{lemmesna2}, on en déduit  que l'induite cohomologique $\rho$ contient $\sigma_{n,k}$. \qed

Nous allons maintenant énoncer une réciproque en montrant que la liste des paquets contenant 
$\sigma_{n,k}$ donnée dans la proposition est complète si $2k\leq n-1$. 
Lorsque $n=2k$, on a $\sigma_{2k,k}=\pi(k+1)$ et la liste des paquets d'Arthur contenant ce module holomorphe unitaire
est donnée dans le théorème \ref{mainthm}. De plus, on a la propriété de multiplicité un.

\begin{thm} \label{Thmsnk} On suppose que $2k\geq n-1$. Soit $\psi$ un $A$-paramètre pour $G=\Sp(2n,\bbR)$ 
dont le caractère infinitésimal est celui de $\sigma_{n,k}$. Alors 
 $\sigma_{n,k}\in \Pi(\psi)$ si et seulement si $a(\psi_u)=2(n-k)+1$ et $\psi$ contient $\sgn_{W_\bbR}^k\boxtimes R[2(n-k)+1]$.
 De plus, la multiplicité de $\sigma_{n,k}$ dans $ \Pi(\psi)$ est un.
\end{thm}

La démonstration sera donnée dans la section \ref{demThmsnk}.

\section{Paramètres de Langlands des représentations $\pi_n(m)$  et $\sigma_{n,k}$.    \label{langlands}}
On fixe  $m$ un entier naturel. Si $m > n$, la représentation $\pi_n(m)$ est une série discrète et si $m=n$ c'est une limite de séries discrètes, 
donc en particulier c'est une représentation tempérée. Si $m=0$, la représentation $\pi_n(m)$ est la représentation triviale qui n'apparaît dans un paquet $\Pi(\psi)$ 
que si $\psi=\psi_u$ est irréductible. 

On suppose donc que $m\in [1,n[$ et on va décrire les paramètres de Langlands de $\pi_n(m)$.

\begin{prop}\label{parLang} Soit $m$ un entier dans $[1,n[$. Soit $P=MN$ le sous-groupe parabolique standard de $G=\Sp(2n,\bbR)$ dont le facteur de Levi $M$
est isomorphe à $\GL(1,\bbR)^{n-m}\times \Sp(2m,\bbR)$.
Alors  la représentation $\pi_n(m)$ est le quotient de Langlands de la représentation standard
$ \Ind_P^G\left( \sgn^m  |\, |^{n-m}\boxtimes \cdots\boxtimes \sgn^m |\, |\boxtimes \pi_{m}(m)\right)$.
 (L'induction parabolique est ici normalisée de la manière usuelle qui préserve l'unitarité. 
 Sur le $j$-ième facteur $\GL(1,\bbR)$,  $\sgn^m  |\, |^{n-m-j+1}$ est  le caractère $x\mapsto  \sgn^m(x)  | x |^{n-m-j+1}$, et l'on étend trivialement
 la représentation  $\sgn^m  |\, |^{n-m}\boxtimes \cdots\boxtimes \sgn^m |\, |\boxtimes \pi_{0}(m)$ ainsi définie en une représentation de $P$ sans changer
 la notation). 
\end{prop}

Le résultat est sans doute déjà bien connu puisque $\pi_n(m)$ est l'image de la représentation triviale du groupe orthogonal compact $\Or(0,2m)$. 
Nous donnons dans la section  \ref{demparlang}
 une preuve globale. 

\

Traitons maintenant le  cas des représentations $\sigma_{n,k}$ de (\ref{sigma})  où $k$ où $2k\leq n$. 
\begin{prop} \label{laglands2}
La représentation $\sigma_{n,k}$ est le quotient de Langlands de la représentation standard
\[
\Ind_P^G \left( \sgn^k |\, |^{n-k}\boxtimes \cdots\boxtimes \sgn^k |\, |^{k+1} \boxtimes \sgn^k |\, |^{k-1}\boxtimes \cdots\boxtimes \sgn^k |\, |\times \pi_{k+1}(k+1)
\right).\]
Ici $P=MN$ est le parabolique standard de $G=\Sp(2n,\bbR)$ comme dans la proposition précédente.
\end{prop}
De même, la démonstration est dans la section \ref{demparlang}.

\

Soit $\psi$ un A-paramètre.
Rappelons les entiers $a(\psi)$ et $a(\psi_u)$ de la définition \ref{apsi}.

\begin{cor}\label{borne}
Soit $\psi$ un A-paramètre tel que $\pi_n(m)\in \Pi(\psi)$ et on suppose 
que $m$ est un entier dans  $[1,n[$.
Alors $a(\psi)\geq 2(n-m)+1$ avec inégalité stricte si 
$a(\psi)>a(\psi_u)$ ou si $\psi$ ne contient pas 
$\sgn_{W_R}^{m}\boxtimes R[2(n-m)+1]$.
\end{cor}

\dem Le plus grand exposant d'une représentation dans $\Pi(\psi)$ est certainement inférieur ou égal à $(a(\psi)-1)/2$ d'après 
la discussion à la fin de la page 155 et le premier paragraphe de la page 156 de \cite{Art13}.
 Comme cet exposant
 pour $\pi_n(m)$ est $n-m$ on obtient l'inégalité large. On a l'inégalité stricte en reprenant la démonstration de \cite{Art13}. 
On peut remarquer que si $m \geq n$ l'inégalité  $a(\psi)\geq 2(n-m)+1$ est  trivialement vérifiée. 
\qed

\section{Un résultat de réduction}\label{reduc}
On suppose dans cette section que  $\psi\neq \psi_u$. 
Le paramètre $\psi$ a donc une partie discrète $\psi_d=\oplus_{j=1}^s (\delta_{t_j}\boxtimes R[a_j])$ non triviale. Rappelons que 
l'on a ordonné les indices de sorte que la suite $(t_j)_{j=1,\ldots,s}$ soit décroissante, et si $t_j=t_j+1$, on a $a_j\geq a_{j+1}$.

On suppose que le caractère infinitésimal de $\psi$ est celui d'un module holomorphe unitaire $\pi_n(m)$ , avec $m\in \{0,\ldots,n\}$ 
 et  on suppose de plus qu'il existe un indice $j_0\in \{1,\ldots,s\}$  tel que :
\begin{equation}\label{condfac} \frac{t_{j_0}+a_{j_0}-1}{2}\geq m-1 \quad \text{  et  }  \quad \frac{t_{j_0}-(a_{j_0}-1)}{2}\geq 0.\end{equation} 
On prend alors  $j_0$ minimal avec cette propriété.

Le  point $(ii)$ du lemme   \ref{longpsiu} montre que  l'inégalité $\frac{t_{j_0}-(a_{j_0}-1)}{2}\geq 0$ est certainement vérifiée si $\dim(\psi_u)>1$,
 mais elle est plus générale. 

\begin{lemme}  \label{rmqred} 
Avec les hypothèses (\ref{condfac}) sur $t_{j_0}$ et $a_{j_0}$, on a $\frac{t_{j_0}+a_{j_0}-1}{2} \in \{m-1, n-m\}$.
\end{lemme}

\dem
Ceci découle de l'égalité du caractère infinitésimaux  de  $\psi$ et  de  $\pi_n(m)$ force. En effet, si  $\frac{t_{j_0}+a_{j_0}-1}{2} >m-1$,  alors $n-m>m-1$ et 
 la multiplicité de $\frac{t_{j_0}+a_{j_0}-1}{2} $ dans le caractère infinitésimal est $1$. Il ne peut donc pas y avoir de terme unipotent $\eta_i \boxtimes R[a_i]$
 avec $\frac{a_i-1}{2}\geq \frac{t_{j_0}+a_{j_0}-1}{2}$, car cela contredirait cette multiplicité $1$.
 Si  $\frac{t_{j_0}+a_{j_0}-1}{2} \neq n-m$, il y a donc nécessairement un autre indice $j>1$
tel que  $\frac{t_j+a_j-1}{2} =n-m$. Si $j>j_0$,  comme $t_{j_0}\geq t_j$, ceci force $a_j>a_{j_0}$ et $\frac{t_j-(a_j-1)}{2} < \frac{t_{j_0}-(a_{j_0}-1)}{2}$ 
ce qui contredit encore la multiplicité 
$1$. Si $j<j_0$ on a $\frac{t_{j}+a_{j}-1}{2}>\frac{t_{j_0}+a_{j_0}-1}{2}>m-1$ et donc on ne peut pas avoir $t_j-a_j+1\geq 0$, car un tel $j$ contredit la minimalité de $j_0$.
On a donc nécessairement $t_j-a_j+1< 0$, mais là encore, la contribution du terme $\delta_{t_j}\boxtimes R[a_j]$ au caractère infinitésimal
comporte le terme $\frac{t_{j_0}+a_{j_0}-1}{2} $ dont la multiplicité $1$ est à nouveau contredite. \qed

\

Fixons une représentation unitaire irréductible $\sigma$ de $G'=\Sp(2(n-a_{j_0}),\bbR)$. 
On va supposer de plus que $\sigma$ est dans un paquet d'Arthur pour ce groupe.
Considérons une induite cohomologique à partir d'une paire parabolique $(\frqqq,L)$ comme en (\ref{qpq})
avec $L$ isomorphe à  $\Sp(2(n-a_{j_0}),\bbR) \times \U(p,q)$, et  $p+q=a_{j_0}$.
Notons  $\Lambda$ le caractère de $\U(p,q)$, où $p+q=a_{j_0}$ de poids $ \lambda= \lambda_{t_{j_0}}$ comme en 
(\ref{lambda1}) et  (\ref{y}). Soit  $\pi_{p,q}(\sigma,t_{j_0})=\caR^S_{\frqqq,L}(\sigma\boxtimes \Lambda)$.
Cette induction   cohomologique, se fait dans le weakly fair range puisque $t_{j_0}\geq 0$.  
On suppose que  $\pi_{p,q}(\sigma,t_{j_0})$ a  le même caractère infinitésimal que $\pi_n(m)$.

\begin{prop}\label{propred} Avec les hypothèses et les notations précédentes, la représentation $\pi_{p,q}(\sigma,t_{j_0})$ contient $\pi_n(m)$ 
si et seulement si $p=0$, $\frac{t_{j_0}+(a_{j_0}-1)}{2}=m-1$ et $\sigma$ est la représentation $\pi_{n-a_{j_0}}(m-a_{j_0})$.
\end{prop}

\begin{rmq}\label{pointdelicat}
L'hypothèse $t_{j_0}\geq 0$ nous dit que l'induction cohomologique qui définit $\pi_{p,q}(\sigma,t_{j_0})$ à lieu dans le weakly fair range, 
mais la représentation $\sigma\boxtimes \Lambda$ n'est pas faiblement unipotente en général, on ne dispose donc pas  a priori des résultats 
d'unitarité de l'induite et d'annulation 
des $\caR^j_{\frqqq,L}(\sigma\boxtimes \Lambda)$ lorsque $j\neq S$. Or nous avons besoin de ces résultats. Ils vont  découler des 
hypothèses sur  $\sigma$, et d'un résultat de P. Trapa pour les groupes unitaires.
Nous donnons l'argument à la fin de cette section.
\end{rmq}

 \dem  On remarque que l'on a supposé que $\frac{t_{j_0}+a_{j_0}-1}{2}\geq a_{j_0}-1$ et quand $\frac{t_{j_0}+a_{j_0}-1}{2}=m-1$ on a donc $m\geq a_{j_0}$, ce qui donne bien
 un sens à la représentation $\pi_{n-a_{j_0}}(m-a_{j_0})$. Supposons donc que $\pi_{p,q}(\sigma,t_{j_0})$ contient $\pi_n(m)$. 
 Montrons d'abord que $p=0$ et que $\frac{t_{j_0}+(a_{j_0}+1)}{2}=m$.
 Nous avons vu dans la lemme ci-dessus  que $\frac{t_{j_0}+a_{j_0}-1}{2}=m-1$ ou $n-m$, c'est-à-dire $t_{j_0}+a_{j_0}+1=2m$ ou $2(n-m+1)$.
 Supposons tout d'abord $t_{j_0}+a_{j_0}+1=2(n-m+1)$. La première  inégalité dans (\ref{condfac})
 s'écrit alors $n-m+1 \geq m$, ou encore $n+1\geq 2m$. 
  L'inégalité fondamentale (\ref{ineqfond})   donne :
  \begin{equation}\label{in2} m(q-p)\geq (p+q)\left(\frac{t_{j_0}+a_{j_0}+1}{2}\right)-2pq=  (p+q)(n-m+1)-2pq,  \end{equation}
  En combinant avec $n+1\geq 2m$, on  obtient $m(q-p)\geq (p+q)m-2pq$,
 ce qui entraîne $2p(q-m)\geq 0 $. 
 Supposons $p\geq1$. On a alors $q\geq m$. La seconde inégalité dans (\ref{condfac})
 s'écrit  $t_{j_0}+1\geq p+q$, ou encore $  n-m+1   \geq p+q$, soit   $n+1   \geq p+q+m$. 
 On réinjecte dans (\ref{in2}) :
 $q(q-p)\geq m(q-p)\geq (p+q)^2-2pq=p^2+q^2$ d'où $-pq\geq p^2\geq 1$, et l'on aboutit à une contradiction.
 On a donc $p=0$, que l'on réinjecte dans (\ref{in2}).
 Alors $mq\geq q(n-m+1)$,  ce qui implique  $2m\geq n+1$ et ainsi $2m=n+1$.  On aboutit ainsi à la conclusion voulue dans ce cas : $p=0$ et 
 $\frac{t_{j_0}+(a_{j_0}+1)}{2}=n-m+1=m$. 

  L'autre cas est $t_{j_0}+a_{j_0}+1=2m$. La seconde inégalité dans (\ref{condfac})
 s'écrit alors $t_{j_0}+1\geq p+q$, ou encore $m\geq p+q$. 
  L'inégalité fondamentale (\ref{ineqfond})   donne :
  \[m(q-p)\geq (p+q)\left(\frac{t_{j_0}+a_{j_0}+1}{2}\right)-2pq=  (p+q)m-2pq,  \]
  ce que l'on réécrit $2 p(q-m)\geq 0$. 
   Si $p\geq 1$,  alors $2p(q-m)\geq 0$ implique  $q\geq m$, ce qui contredit $m\geq p+q\geq q+1$.
    Donc nécessairement, $p=0$. On a donc $q=a_{j_0}$.
    Ceci montre la nécessité des deux premières conditions de l'énoncé.

    Calculons la multiplicité du $K$-type $E=\bbC_{-m}$ de dimension $1$ et de plus haut poids scalaire $(-m,\ldots,-m)$ dans  
      $\pi_{0,a_{j_0}}(\sigma,t_{j_0})= \caR^S_{\frqqq,L,G}(  \sigma \boxtimes \Lambda)$ avec $\Lambda= \det {}^{n-\frac{t_{j_0}+a_{j_0}-1}{2}}$.
      On utilise la formule de multiplicité du théorème 5.64 de \cite{KnVo}, et le fait mentionné dans la remarque \ref{pointdelicat} qui donne l'annulation
   des induites cohomologique   en degré $j\neq S$ :  
   \[ (-1)^S \dim(\Hom_K(E;  \pi_{0,a_{j_0}}(\sigma,t_{j_0}) ))= (-1)^S \dim(\Hom_K(E;   \caR^S_{\frqqq,L}(\sigma\boxtimes \Lambda)))\]
  \[ = \sum_{r=0}^S(-1)^r   \dim(\Hom_{K\cap L}( H_r(\fru\cap \frk; E); 
   S(\fru\cap \frp) \otimes (\sigma \boxtimes \Lambda) \otimes  \bbC_{2\delta(\fru)})).\]
  Combinée avec la dualité de Poincaré du corollaire 3.8 de \cite{KnVo}, le changement d'indice $r\to S-r$ et la formule de la section 
  \ref{ppm} pour $2\delta(\fru)$, on obtient :
  \begin{align*}&\dim(\Hom_K(E;   \pi_{0,a_{j_0}}(\sigma,t_{j_0})    )) \\
&= \sum_{r=0}^S  (-1)^r   \dim(\Hom_{K\cap L}( H^r(\fru\cap \frk; E);  S(\fru\cap \frp) \otimes    (  \sigma \boxtimes 
\det{}^{n-\frac{t_{j_0}+a_{j_0}-1}{2}} ) \otimes \bbC_{2\delta(\fru\cap \frp)})\\
&= \sum_{r=0}^S  (-1)^r   \dim(\Hom_{K\cap L}( \wedge^r(\fru\cap \frk)\otimes  \bbC_{-m} ;  S(\fru\cap \frp) \otimes     \sigma   \otimes 
\bbC_{-m-a_{j_0}}).\end{align*}

  D'autre part, rappelons que 
  \[ \fru\cap \frp=\bigoplus_{n-a+1\leq i<j\leq n} \left( \frg_{-e_i-e_j}\oplus \frg_{-2e_i} \right)
  \oplus \bigoplus_{ 1\leq i\leq   n-a<j\leq n} \frg_{-e_i-e_j}  =U_1\oplus U_2 \]
   On a donc $S(\fru\cap \frp)=\bigoplus_{\alpha,\beta} S^\alpha(U_1)\otimes S^{\beta}(U_2)$ et l'on regarde 
\[ \Hom_{K\cap L}( \wedge^r(\fru\cap \frk)\otimes  \bbC_{-m} ;  S^\alpha(U_1) \otimes S^\beta(U_2) \otimes     \sigma   \otimes 
\bbC_{-m-a_{j_0}}) \]
En regardant l'action du centre de l'algèbre enveloppante du facteur $\U(a)$, on obtient une condition nécessaire pour la non nullité de
     cet espace, à savoir   $2\alpha+\beta= -r$.
  Ceci force $\alpha=\beta=r=0$ et l'on obtient finalement
  \begin{equation} \label{HomK} \dim(\Hom_K(E;  \pi_{0,a_{j_0}}(\sigma,t_{j_0})  )) =
\dim(\Hom_{\U(n-a_{j_0}) }( \mathbb{C}_{-m+a_1},\sigma)).
\end{equation}
On vérifie que l'hypothèse sur les caractères infinitésimaux force $\sigma$ à avoir comme caractère infinitésimal: $m-a_{j_0}-1, \cdots, n-m$, 
c'est-à-dire le caractère infinitésimal de $\pi_{n-a_{j_0}}(m-a_{j_0})$.

Pour conclure, nous allons utiliser le résultat suivant de Chen-Bo Zhu  \cite{Zhu} :
\begin{lemme}\label{lemmemanquant}
Soit  $\sigma$ une représentation irréductible unitaire de 
$G=\Sp(2n,\bbR)$ ayant même caractère infinitésimal que $\pi_n(m)$. Supposons de plus que $\sigma$ contienne le $K$-type 
scalaire $(\underbrace{-m,\ldots,-m}_{n})$. Alors $\sigma=\pi_n(m)$.
\end{lemme}

Grâce à ce lemme, appliqué à $G'=\Sp(2(n-a_{j_0},\bbR)$, $\pi_{n-a_{j_0}}(m-a_{j_0})$ et $\sigma$, 
  la non nullité du terme de droite dans (\ref{HomK}) est exactement équivalent à ce que $\sigma=\pi_{n-a_{j_0}}(m-a_{j_0})$.
Dans ce cas, la dimension de cet espace est $1$.   Comme $\pi_{0,a_{j_0}}(\sigma,t_{j_0})$ est une induite cohomologique dans le weakly
fair range, elle est unitaire, et certainement non nulle si la dimension de  (\ref{HomK}) est $1$. Le lemme ci-dessus, à nouveau appliqué
aux composantes irréductibles de $\pi_{0,a_{j_0}}(\sigma,t_{j_0})$ montre le résultat voulu, c'est-à-dire que l'une de ces composantes est $\pi_n(m)$.
Ceci termine la démonstration de nécessité des conditions de la proposition. Il est facile de montrer qu'elle sont suffisantes
en utilisant le lemme ci-dessus et (\ref{HomK}).
\qed 
\medskip 

\begin{rmq} \label{rmq2red}Soit $\psi'$ le $A$-paramètre pour $G'=\Sp(2(n-a_{j_0}),\bbR)$ défini par 
 $\psi=(\delta_{t_{j_0}}\boxtimes R[a_{j_0}])\oplus (\sgn_{W_\bbR}^{a_{j_0}}\otimes \psi')$. 
Dans les hypothèses de la proposition, si l'on suppose de plus  $\frac{t_{j_0}+a_{j_0}-1}{2}=m-1$, alors  
 $(\psi, \pi_n(m))$  satisfait les hypothèses du théorème \ref{mainthm} si et seulement si $(\psi', \pi_{n-a_{j_0}}(m-a_{j_0}))$ les satisfait pour $\Sp(2(n-a_{j_0}),\bbR)$.
\end{rmq}
\dem
Evidemment $\dim (\psi_u)=\dim (\psi'_u)$ et $a(\psi_u)=a(\psi'_u)$. Dans le cas  $(ii)$ du théorème, avec
 $a(\psi_u)=a(\psi'_u)=2(n-m)+1$, 
il n'y a aucune difficulté, mais si  
$a(\psi_u)=2(n-m)+3$, il faut vérifier que les hypothèses forcent 
$2(m-a_{j_0})\geq (n-a_{j_0})+2$. Mais si $\frac{a(\psi_u)-1}{2}=n-m+1$ apparaît  dans le caractère infinitésimal de $\psi$, il apparaît avec multiplicité $1$
et ne peut pas être dans l'intervalle $[\frac{t_{j_0}-a_{j_0}+1}{2},\frac{ t_{j_0}+a_{j_0}-1}{2}]=[m-a_{j_0},m-1]$ 
puisque les éléments de cet intervalle sont  dans le caractère infinitésimal grâce à la contribution de  $\delta_{t_{j_0}}\boxtimes R[a_{j_0}]$.
On a donc $m-a_{j_0}\geq n-m+2$ et ainsi $2(m-a_{j_0})\geq (n-a_{j_0})+2$.
On remarque aussi que les inégalités impliquent aussi que $\frac{t_{j_0}-a_{j_0}+1}{2}=m-a_{j_0}>m-a_{j_0}-1\geq (n-a_{j_0})-(m-a_{j_0})$, et donc l'induction
 cohomologique se fait
dans le good range dans ce cas.\qed

\

Il nous reste à démontrer les résultats d'unitarité et d'annulation d'induites cohomologiques de la remarque \ref{pointdelicat}.

\dem Si $\sigma$ est dans un paquet d'Arthur unipotent pour $G'$, elle est en particulier faiblement unipotente (\cite{MR5}), et les résultats
déjà mentionnés de \cite{KnVo} s'appliquent dans ce cas.
Si $\sigma$ n'est pas unipotente  (ici, \og unipotente \fg\,  signifie \og appartenant à un paquet d'Arthur unipotent\fg),
 on sait d'après la description des paquets d'Arthur pour les groupes symplectiques
rappelés dans la section  \ref{ArtpacMR} que $\sigma$ est obtenue comme sous-représentation d'une induite cohomologique $\tilde \sigma$
à partir d'un $c$-Levi  $L'$ de $G'$ isomorphe à produit de groupes unitaires et d'un groupe symplectique $G''$ de rang plus petit.
La représentation que l'on induit est unipotente sur le facteur $G''$, et ce sont des caractères unitaires sur les groupes unitaires.
On est dans le weakly fair range, avec les résultats d'annulation en tout degré sauf un, et d'unitarité qui en découlent.
On va démontrer  l'unitarité en degré $j=S$ et l'annulation en degré $j\neq S$ des induites cohomogiques obtenues en remplaçant
$\sigma$ par $\tilde \sigma$, ce qui suffit.
Grâce au lemme d'induction par étape (\cite{Vgreen}, Cor. 6.3.10) on écrit 
$\caR^j_{\frqqq,L}(\tilde \sigma\boxtimes \Lambda)$ comme une induite cohomologique
à partir d'un $c$-Levi de $G$ isomorphe à  $L'\times \U(p,q)$ (rappelons que ce facteur $\U(p,q)$ vient de l'indice
$j_0$, il aurait été plus clair mais trop lourd  d'écrire $\U(p_{j_0}, q_{j_0})$), et la représentation que l'on induit est 
alors faiblement unipotente, le problème est maintenant que l'on a sorti l'indice $j_0$ du lemme pour le mettre
en premier, et que l'induction n'est plus dans le weakly fair range. C'est ici que l'on va utiliser les résultats de P. Trapa \cite{TrapA}
en décomposant l'induction cohomologique en deux étapes, l'étape intermédiaire
étant constitué du $c$-Levi de $G$ isomorphe au produit d'un seul groupe unitaire et de $G''$.
L'étape intermédiaire peut donc être vue comme ayant lieu uniquement dans ce groupe unitaire.
On induit un produit de caractères, et l'on est dans ce que Trapa appelle le \og mediocre range \fg.
On peut utiliser alors les égalités entre $A_\frqqq(\lambda)$  démontrées par Trapa pour voir que le résultat de cette induction
cohomologique dans le mediocre range est égale à celle obtenue en remettant l'indice $j_0$ à sa place, et qui elle est dans 
le weakly fair range, où l'on a les résultats d'annulation voulus.  On peut donc appliquer à nouveau \cite{Vgreen}, Cor. 6.3.10
sur l'induction par étape pour montrer l'annulation de $\caR^j_{\frqqq,L}(\tilde \sigma\boxtimes \Lambda)$
lorsque $j\neq S$, et de même pour l'unitarité en degré $S$. \qed

\section{Démonstration du théorème \ref{mainthm} (les conditions sont suffisantes) et du théorème \ref{complement}}\label{sensdirect}

Dans cette section, on suppose que le couple $(\psi,m)$ vérifie une des conditions du théorème \ref{mainthm}, et l'on va montrer que $\pi(m)=\pi_n(m)$ est dans 
$\Pi(\psi)$. Rappelons  que l'on a posé
\[\psi=\psi_u\oplus \psi_d=\psi_u\oplus  \bigoplus_{i=1}^s   (\delta_{t_j}\boxtimes R[a_j]).\]

\subsection{Cas $(i)$}\label{71}

On suppose donc que $\dim(\psi_u)=1$,  $2m>n+1$ et que l'on a (\ref{disj}). On va raisonner par récurrence sur   la longueur de  $\psi_d$. Supposons tout d'abord que 
$\psi_d$ soit de longueur $1$, de sorte que 
$\psi=( \delta_t\boxtimes R[n])\oplus (\sgn_{W_\bbR}^n\boxtimes R[1])$. 
Considérons la paire  $(\frqqq, L)$, où  $L=K$ et $\frqqq=\frk\oplus \frp^-$. La condition sur le caractère infinitésimal est 
que $\frac{t+n-1}{2}=\max(m-1,n-m)=m-1$.

On considère l'induite cohomologique 
$\caR^0_{\frqqq, L }(\Lambda)$
où $\Lambda$ est le caractère de $L=K$ de différentielle
\[\lambda=(-m+n+1,\ldots,-m+n+1)=\left(\frac{-t+n+1}{2}, \ldots,\frac{-t+n+1}{2} \right).\]
Ici $\lambda+2\delta(\fru\cap \frp)=\left(-\frac{t+n+1}{2}, \ldots,-\frac{t+n+1}{2} \right)=(-m,\ldots,-m)$ qui est dominant 
pour $\Delta(\frg,\frk)^+$, et c'est donc le 
bottom layer de $\caR^0_{\frqqq, L }(\Lambda)$.
D'autre part,  $\caR^0_{\frqqq, L }(\Lambda)$ est irréductible (c'est un module de Verma généralisé, donc monogène, et il est de plus unitaire, donc semi-simple) et 
 $\caR^0_{\frqqq, L}(\Lambda)=\pi(m)$.  
 D'après \cite{MR3}, $\pi(m)$ est  alors dans 
$\Pi(\psi)$.

On suppose maintenant  
\[ \psi=(\sgn_{W_\bbR}^{\vert \{j;\, 1\leq j\leq s, \, a_j \text{ impair } \}  \vert} \boxtimes R[1])\oplus \bigoplus_{j=1}^s( \delta_{t_j}\boxtimes R[a_j]), \]
avec $s>1$ 
et  le résultat démontré pour 
\[\psi^\flat=(\sgn_{W_\bbR}^{\vert \{j; \, 2\leq j\leq s, \,a_j \text{ impair }\}  \vert} \boxtimes R[1])\oplus \bigoplus_{j=2}^{s} (\delta_{t_j}\boxtimes R[a_j]).\]

Dans \cite{MR3}, les éléments de $\Pi(\psi)$ sont décrits comme composantes irréductibles d'induites cohomologiques
dans le weakly fair range à partir de paires paraboliques $(\frqqq,L)$ comme dans la section \ref{ppm} avec  $\frqqq=\frqqq_{p,q}$ et $p+q=a_1$.
On prend ici $(p,q)=(0,a_1)$, le $c$-Levi $L$ est alors isomorphe à $\Sp(2(n-a_1),\bbR)\times \U(0,a_1)$. 
 Comme représentation induisante, sur le facteur $\Sp(2(n-a_1),\bbR)$, on met le module holomorphe unitaire  de plus haut poids 
 $(\underbrace {m-a_1,\ldots,m-a_1}_{n-a_1})$,
 qui par hypothèse de recurrence est dans $\Pi(\psi^\flat)$, 
 et sur le facteur  $\U(0,a_1)$, on met  le caractère unitaire $\Lambda$ de différentielle $(\underbrace {n-\frac{t_1+a_1-1}{2},\ldots,n-\frac{t_1+a_1-1}{2}}_{a_1})$.
Par induction cohomologique par étapes, on obtient $\pi(m)$.
Remarquons que l'hypothèse (\ref{disj}) entraîne qu'après la première induction cohomologique (celle pour le paramètre $\delta_{t_s}\boxtimes R[a_s]$, 
les autre inductions cohomologiques successives se font dans le good range, et en particulier, elle préservent l'irréductibilité.
 Ceci montre que $\pi(m)$ est bien dans $\Pi(\psi)$ et établit aussi le théorème \ref{complement} dans le cas $(i)$.

\subsection{Cas $(ii)$ avec  $a(\psi_u)=2(m-1)+1$ et $m-1>n-m$} \label{subsectplus3}
On suppose que le couple $(\psi,m)$ vérifie l'hypothèse $(ii)$ du théorème, avec de plus $\psi\neq \psi_u$
 (le cas $\psi=\psi_u$ a été traité dans la section
\ref{parunip}), $\dim\psi_u>1$,   $a(\psi_u)=2(m-1)+1$ et  $m-1>n-m$.
D'après le lemme \ref{rmqapsiu}, on a  $a(\psi_u)=2(n-m)+3$, et donc $n=2(m-1)$. De plus $\pi(m)$ est l'image 
par la correspondence de Howe du déterminant de   $\Or(0,2n+1-a(\psi_u))=\Or(0,2(m-1))$.
On a donc  $\pi(m)=\sigma_{n=2k,k}$ avec les notations de la section \ref{sigmana} et  $k=m-1$.
On conclut alors grâce à la proposition \ref{propsna}.

\subsection{Les conditions du théorème sont suffisantes : fin de la démonstration} \label{findemsuf}
On a déjà traité    les cas   $\psi=\psi_u$ et  $\dim(\psi_u)=1$.
On se place donc dans le cas $(ii)$ du théorème \ref{mainthm} avec $\psi\neq \psi_u$, et l'on raisonne par récurrence sur la longueur 
de $\psi_d$. On a donc  $a(\psi_u)=2(n-m)-1$ ou $2(n-m)+3$.

Si $m-1>n-m$, et si $a(\psi_u)=2(m-1)+1$, on est dans la situation de la section \ref{subsectplus3} et ce cas  a donc déjà été traité.
Si $a(\psi_u)<2(m-1)+1$, les considérations sur le caractère infinitésimal montrent qu'il existe un unique indice $j_0\in \{1,\ldots,s\}$ tel que  $\frac{t_{j_0}+a_{j_0}-1}{2}=m-1$
et  d'après le lemme \ref{longpsiu}, on a       $t_{j_0}-a_{j_0}+1\geq 0$.
 Ainsi $j_0$ vérifie l'hypothèse (\ref{condfac})  et il est minimal avec cette propriété, car unique.  On peut donc appliquer  la proposition \ref{propred}
  et la remarque  \ref{rmq2red}. On écrit $\psi= (\delta_{t_{j_0}} \boxtimes R[a_{j_0}])\oplus (\sgn_{W_\bbR}^{a_{j_0}}\otimes \psi')$. 
 Alors $(\psi', \pi_{n-a_{j_0}}(m-a_{j_0}))$ vérifie les conditions du théorème pour $G'=\Sp(2(n-a_{j_0}),\bbR)$.
  Donc $\Pi(\psi')$ contient  $\pi_{n-a_{j_0}} (m-a_{j_0})$ et $\pi_n(m)$ est une composante irréductible de 
  $\pi_{p,q}(\sigma,t_{j_0})=\caR^S_{\frqqq,L}(\sigma\boxtimes \Lambda)$ d'après la proposition \ref{propred}.
 On conclut que  $\pi_n(m)\in \Pi(\psi)$  d'après la description de $\Pi(\psi)$ de \cite{MR3} rappelée en \ref{ArtpacMR}. 

Si $n-m\geq m-1$, alors on a nécessairement $a(\psi_u)=2(n-m)+1$. De plus, par des considérations sur le caractère infinitésimal
comme dans le lemme \ref{longpsiu} et la remarque \ref{rmqred}, on voit que les segments $[ \frac{t_j-a_j+1}{2},  \frac{t_j+a_j-1}{2} ]$ doivent pas s'intersecter  lorsque $j$ parcourt 
$\{1,\ldots s\}$, et de plus, on a nécessairement 
  $\frac{t_1+a_1-1}{2}=m-1$ et $t_1-a_1+1\geq 0$.
On peut donc appliquer la proposition \ref{propred} et la remarque \ref{rmq2red} et conclure par hypothèse de récurrence comme ci-dessus.
On obtient aussi le   théorème \ref{complement} dans le cas $(ii)$.

\section{Démonstration du théorème \ref{mainthm}: les conditions sont nécessaires}\label{condnec}
\subsection{Cas $\psi=\psi_u$}
On écrit \[\psi=(\eta_1\boxtimes R[a_1])\oplus (\eta_2\boxtimes R[a_2])\oplus (\eta_1\eta_2\boxtimes R[1]),\] avec $a_1\geq a_2$.
 On doit donc démontrer:
\begin{lemme} On suppose que $\pi(m)$ est dans $\Pi(\psi)$. On a alors $a_1=2(n-m)+1$ ou $2(n-m)+3$ et  le caractère $\eta_1$ vaut
 $\sgn_{W_\bbR}^{(a_2+1)/2}$.
\end{lemme}

Supposons d'abord que $a_1=a_2=n$. L'égalité des caractères infinitésimaux de $\pi(m)$ et de $\psi$
force les égalités $n=2m-1=2(n-m)+1=a_1$. D'après  le corollaire \ref{borne} , le paramètre 
doit nécessairement contenir  $\sgn_{W_\bbR}^m\boxtimes R[n]$. On a donc la conclusion voulue dans ce cas.

Si l'on n'est pas dans dans le cas  $a_1=a_2=n$, on a nécessairement  $a_1\geq n+1$, et l'on peut 
appliquer le corollaire \ref{SW} et   les représentations dans $\Pi(\psi)$ sont alors
 des images par la correspondance de Howe de caractères des groupes $\Or(p,q)$ tels que $p+q=a_2+1$ 
 et le discriminant normalisé  (voir section \ref{discethasse}) de la forme quadratique est  $\eta_1=\sgn_{W_\bbR}^{\frac{p-q}{2}}$ 
  Comme, par hypothèse, $\pi(m)$ est dans $\Pi(\psi)$, il existe un couple  $(p,q)$ avec $p+q=a_2+1$ et  et un  caractère $\tau$ 
  de $\Or(p,q)$  tels que $\pi(m)$ soit obtenu comme image de $\tau$ par la correspondance de Howe. 
   Il y a deux caractères du groupe orthogonal si $pq=0$ et quatre sinon, ils sont facilement indexés par le choix de deux signes, 
   $\epsilon$ et $\eta$, définis de la façon suivante : le caractère $\tau(\epsilon,\eta)$ de $\Or(p,q)$
   restreint au sous-groupe compact maximal $\Or(p,0)\times \Or(0,q)$, est égal au produit tensoriel de $\det_{\Or(p,0)}^{\frac{1-\epsilon}{2}}$
et  $\det_{\Or(0,q)}^{\frac{1-\eta}{2}}$
On fixe donc $\epsilon$ et $\eta$ tel que  $\tau=\tau(\epsilon,\eta)$.
On sait alors que $\pi(m)$ contient la représentation de $\U(n)$ de plus haut poids:
\begin{equation}\label{phpUn}
\frac{p-q}{2}+(\underbrace{1,\quad  \cdots \quad ,  1}_{p\frac{1-\epsilon}{2}}, \underbrace{0,\qquad \ldots \qquad, 0}_{n-p\frac{1-\epsilon}{2}-q\frac{1-\eta}{2}},
\underbrace{-1, \quad \ldots \quad, -1}_{q\frac{1-\eta}{2}})
\end{equation}
et ce poids est nécessairement $(-m,\ldots,-m)$.

Si   $p=0$ et  si $\eta=1$, on obtient $q=2m$, d'où $a_2=2m-1$.
En comparant les caractères infinitésimaux, on a nécessairement:
\begin{equation}\label{infsup}
\frac{a_1-1}{2}=\sup(m-1,n-m), \quad \frac{a_2-1}{2}=\inf(m-1,n-m).
\end{equation}
On en déduit que $\inf(m-1,n-m)=m-1$ et donc $\frac{a_1-1}{2}=n-m$, soit $a_1=2(n-m)+1$.
Le caractère $\eta_1$ est donné par le discriminant de la forme quadratique de signature $(0,a_2+1)$, c'est-à-dire 
$\eta_1= \sgn_{W_\bbR}^{(a_2+1)/2}$. On a donc la conclusion voulue dans ce cas. 

Supposons maintenant toujours $p=0$,  et maintenant  $\eta=-1$.
Ceci implique $-q/2-1=-m$, c'est-à-dire $q=2(m-1)$. Ensuite, on obtient $q=n$, puis $a_2=q-1=2(m-1)-1=2m-3$, et avec 
(\ref{infsup}), $\frac{a_2-1}{2}=\inf(m-1,n-m)= m-2=n-m  $ et $a_1=2m-1=2(n-m)+3$. On a encore $\eta_1= \sgn_{W_\bbR}^{(a_2+1)/2}$
comme dans le cas précédent, et là encore, on a obtenu la conclusion voulue.

Si $p\neq 0$, et si $\epsilon=-1$, on a alors $p-q=-2m$, $n-p-q\frac{1-\eta}{2}=0$ et $q\frac{1-\eta}{2}=0$. On en déduit $n=p$, 
$q=2m+n$. Or $p+q=2(m+n)=a_2+1$ et en utilisant (\ref{infsup})
\[a_2+1\leq \frac{a_1-1}{2}+ \frac{a_2-1}{2}+2\leq m-1+n-m+2=n-1,\]
ce qui n'est pas possible. Ce cas ne se produit donc pas.

Si $p\neq 0$, et si $\epsilon=1$, on distingue encore les cas $q\frac{1-\eta}{2}=0$ et $q\frac{1-\eta}{2}\neq 0$.
Dans le premier cas, on obtient $p-q=-2m$. En utilisant (\ref{infsup}), on obtient 
\[a_2-1=p+q-2=2m+2p-2\leq 2(m-1),\]
ce qui contredit $p\neq 0$. Dans le second cas, on a $p-q=2(m-1)$ et $n=q$.
Avec 
 (\ref{infsup}), on obtient cette fois  
\[a_2-1=p+q-2=2(m-1)+2p-2\leq 2(m-1),\]
d'où $p=1$. Ensuite, on en déduit que $a_2-1=2(m-1)$, d'où $a_1=2(n-m)+1=a_2$, ce qui contredit $a_1\geq n+1$.

\begin{rmq}
Dans le cas $n=2m-1$, le calcul ci-dessus montre que  $\pi(m)$ est à la fois l'image par la correspondance de Howe du caractère $\tau(1,-1)$ de $\Or(1,n)$, et 
 de la représentation triviale du groupe compact $\Or(0,n+1)$ ({\sl cf.}   \cite{Zhu}, Thm. 4.2.2). 
\end{rmq}

\subsection{Nécessité des conditions dans le cas où $\dim(\psi_u)=1$}\label{npu}

On suppose ici  que $\psi=\oplus_{j=1}^s (\delta_{t_j}\boxtimes R[a_j] )\oplus( \sgn_{W_\bbR}^n  \boxtimes R[1]) $ et l'on suppose 
que $\Pi(\psi)$ contient un module holomorphe unitaire $\pi(m)$, avec $0\leq m\leq n$.
On veut montrer que $2m>n+1$ et que la condition (\ref{disj}) est satisfaite. 
On raisonne par récurrence sur $s$. Commençons par le cas $s=1$, c'est à dire   $\psi=( \delta_{t_1}\boxtimes R[a_1])\oplus (\sgn_{W_\bbR}^{n}  \boxtimes R[1])$. 
Dans ce cas $n=a_1$  et donc $a(\psi)=n>a(\psi_u)=1$. Le corollaire \ref{borne}
nous donne $n>2(n-m)+1$, soit $2m>n+1$, et la condition  (\ref{disj}) est trivialement vérifiée.

Supposons maintenant $s\geq 2$ et  le résultat établi pour des paramètres avec $\dim(\psi_u)=1$ et au plus $s-1$ termes dans la partie discrète de $\psi$.
Supposons $n-m\geq m-1$, on doit trouver une contradiction. Soit $j_1$  un indice dans $\{1,\ldots,a\}$ tel que 
$\frac{t_{j_1}+a_{j_1}-1 } {2}=n-m$. On suppose que $t_{j_1}-a_{j_1}+1\geq 0$. 
Alors l'ensemble des indices vérifiant la condition (\ref{condfac}) est non vide, et l'on peut prendre $j_0$ minimal 
dans cet ensemble. La proposition \ref{propred} et la remarque \ref{rmq2red} nous donnent alors, avec les notations afférentes
$\pi_{n-a_{j_0}}(m-a_{j_0})\in \Pi(\psi')$. Mais d'autre part 
$ (n-a_{j_0})-(m-a_{j_0})=m-n\geq m-1> m-a_{j_0}-1 $
et ceci contredit l'hypothèse de récurrence appliquée à $\psi'$, qui a une partie discrète de longueur $s-1$, et $\pi_{n-a_{j_0}}(m-a_{j_0})$.
Supposons maintenant que $t_{j_1}-a_{j_1}+1 <  0$. Il existe  un indice $j_2\neq j_1$ tel que 
 $\frac{t_{j_2}+a_{j_2}-1}{2}= m-1$ et nécessairement 
 $\frac{t_{j_2}-a_{j_2}+1}{2}\geq 0$, d'après le lemme (\ref{longpsiu}).
Il n'est pas difficile de voir que $j_2$ est  l'unique indice vérifiant les propriétés  (\ref{condfac}) et l'on peut donc appliquer 
la proposition \ref{propred} et la remarque \ref{rmq2red} : on a alors 
$\pi_{n-a_{j_2}}(m-a_{j_2}) \in \Pi(\psi')$ où $\psi'$ est défini par 
$ \psi=(\delta_{t_{j_2}}\boxtimes R[a_{j_2}] )\oplus (\sgn_{W_\bbR}^{a_{j_2}}\psi').$
Mais là encore 
$  (n-a_{j_2})-(m-a_{j_2})=m-n\geq m-1> m-a_{j_2}-1 $
et ceci contredit l'hypothèse de récurrence appliquée à $\psi'$ et  $\pi_{n-a_{j_0}}(m-a_{j_2})$.

Regardons maintenant le cas  $m-1>n-m$. Il reste à montrer que la condition (\ref{disj}) est satisfaite.
Soit $j_0\in \{1,\ldots,s\}$ tel que  $\frac{t_{j_0}+a_{j_0}-1}{2}=m-1$. Si $t_{j_0}-a_{j_0}+1\geq 0$,  
alors $j_0$ vérifie (\ref{condfac}), et c'est l'unique indice qui peut vérifier ces conditions. On  applique la proposition \ref{propred},
et la remarque \ref{rmq2red} :
il faut alors que $\pi_{n-a_{j_0}}(m-a_{j_0})$ soit dans $\Pi(\psi')$ où $\psi'$ est défini par 
  $ \psi=(\delta_{t_{j_0}}\boxtimes R[a_{j_0}]) \oplus (\sgn_{W_\bbR}^{a_{j_0}}\psi').$
D'après l'hypothèse de récurrence, les conditions (\ref{disj}) sont satisfaites par le paramètre $\psi'$, et de plus $m-a_{j_0}-1>n-a_{j_0}-(m-a_{j_0})=n-m$
On a alors 
$ \frac{t_{j_0} - a_{j_0}+1}{2} =  \frac{t_{j_0} + a_{j_0}-1}{2}-(a_{j_0}-1)=m-1-(a_{j_0}-1 ) =m-a_{j_0} >n-m+1$,
ce qui montre que les conditions (\ref{disj}) sont satisfaites par $\psi$.

On suppose  maintenant $t_{j_0}-a_{j_0}+1 <0$. Comme $s\geq 2$, il existe un indice $j_1\neq j_0$ dans $\{1,\ldots, s\}$
que l'on prend minimal.
 On va donner un argument utilisant la variété associée. C'est un invariant des modules de Harish-Chandra
défini par Vogan dans \cite{VogAss}. Il consiste en une union de $K_\bbC$-orbites  nilpotentes dans $\frp$, où $K_\bbC$ est la complexification du sous-groupe
compact maximal $K$ de $G$. Or, la variété associée des modules unitaires holomorphes est connue, voir par exemple \cite{NOT}, \S 7.3.
C'est une unique $K_\bbC$-orbite, contenue dans $\frp^-$.  
D'autre part, on dispose d'une paramétrisation combinatoire des $K_\bbC$-orbites  nilpotentes dans $\frp$, voir par exemple \cite{Oht}, et l'ordre naturel 
sur les orbites  donné par l'inclusion dans l'adhérence et bien décrit combinatoirement.
Dans le cas du groupe symplectique $G=\Sp(2n,\bbR)$,  les orbites sont paramétrées par certains tableaux de Young signés. Tout d'abord, on a un tableau de Young, 
c'est-à-dire une partition, qui paramètre une orbite nilpotente dans l'algèbre de Lie $\frg$, et pour le groupe symplectique, ces tableaux ont 
$2n$ cases, et les lignes de longueur impaires doivent apparaitre avec une multiplicité paire. Pour paramétrer les $K_\bbC$-orbites  nilpotentes dans $\frp$
contenues dans une orbite nilpotente de $\frg$, on met en plus des signes dans les cases du tableau de Young, qui doivent alterner le long des lignes du tableau.
De plus, pour les lignes d'une longueur impaire donnée (donc un nombre pair de lignes), la moitié doit commencer par $+$ et l'autre moitié
commencer par $-$. Les orbites contenue dans $\frp^-$ sont celles dont les lignes ont au plus deux éléments, et les lignes à deux éléments
commencent toutes par $+$. L'ordre naturel induit sur les orbites contenues dans $\frp^-$ est linéaire, et les orbites ont des dimensions distinctes. 
L'orbite dense dans $\frp^-$ est donc celle correspondant au tableau ayant $n$ lignes de longueur deux, toutes étant égales à $+-$, et l'orbite $\{0\}$ correspond bien sûr au tableau 
à $2n$  lignes de longueur $1$; $n$ étant $+$ et $n$ étant $-$.
La variété associée à $\pi_n(m)$ est paramétrée par le tableau signé suivant : si $2m\geq n$,  la variété associée est l'orbite dense dans $\frp^-$.
Si $2m<n$, elle est paramétrée par le tableau à $2m$ lignes de longueur $2$ toutes égales à  $+-$; et $n-2m$ lignes de longueur $1$ égale à $+$ 
et $n-2m$ lignes de longueur $1$ égale à $-$.

  Les éléments du paquet $\Pi(\psi)$ sont les composantes des induites cohomologiques à partir de  $c$-Levi de la forme 
  $ \times_{j=1}^s \U(p_j,q_j)  $
  avec $p_j+q_j=a_j$. 
  Or, on sait calculer la variété associée de ces induites cohomologiques en termes de leur paramétrisation par les tableaux de Young signés 
  comme ci-dessus, voir par exemple 
  \cite{Trap} qui donne un algorithme simple.
  Si l'un des $p_j$ est non nul,  la variété associée de tels induites 
  comporte au moins une ligne à au moins 4 cases si l'un des $q_j$ est non nul, 
ou bien seulement  des lignes $-+$ si tous les $q_j$ sont nuls.
  Dans ce dernier cas $\pi_n(m)$ n'est certainement pas composante de l'induite, car $(+-)^{n}$ n'est pas dans l'adhérence de 
  $(-+)^n$. Mais le résultat est vrai aussi dans le premier cas, car l'induite cohomologique est ici irréductible, par un critère dû à Bernstein et Vogan : 
  on induit d'un sous-groupe parabolique $\theta$-stable $Q=LU$,  l'orbite de Richardson  $\caO$
 correspondant à ce parabolique est d'adhérence normale, et l'application moment  $T^*(G/Q)\rightarrow \caO$ est birationnelle. 
 La liste des sous-groupes paraboliques vérifiant ce critère est donnée en  \cite{BaCG}, \S14, et l'on vérifie donc que tel est le cas ici.
La démonstration de l'irréductibilité de l'induction cohomologique (dans le weakly fair range, ce qui est le cas ici) lorsque le critère
est vérifié se trouve dans \cite{VogDS}.

  Si tous les $p_j$ sont nuls, on a un $c$-Levi compact, et on peut faire une induction par étape, en passant par $\U(n)$.
  En utilisant l'argument donné dans la démonstration de la remarque  \ref{pointdelicat} à la fin de la section \ref{reduc}
  faisant appel au résultats \cite{TrapA}, on permute s'il le faut des indices pour que les indices $j_0$ et $j_1$
     deviennent adjacents, tout en restant dans le mediocre range.
  Mais le passage par $\U(a_{j_0})\times \U(a_{j_1})$ avec $(t_{j_0},a_{j_0})$, $(t_{j_1},a_{j_1})$ comme ci-dessus donne $0$ dans $\U(a_{j_0}+a_{j_1})$, 
  et l'on aboutit à une contradiction.  \qed

\subsection{Fin de la démonstration de la nécessité des conditions}

Il nous reste le cas où $\dim(\psi_u)>1$ et $\psi\neq \psi_u$. Remarquons que cela force $n>1$. 
D'après le lemme \ref{longpsiu}, on a  $\frac{t_1-a_1+1}{2}\geq 0$ puisque $\dim(\psi_u)>1$  et si de plus  $\frac{t_1+a_1-1}{2} \geq m-1$, 
 on peut invoquer   la réduction de  la section \ref{reduc}. On conclut rapidement dans ce cas grâce à la remarque \ref{rmq2red}. 
  
  On suppose donc que  $\frac{t_1+a_1-1}{2} < m-1$. Supposons d'autre part $n-m\geq m-1$ et montrons que cela mène à une contradiction.
En effet,  $m-1$ a multiplicité 2 dans le caractère infinitésimal. Si un  indice $j\in [2,s]$ de la partie discrète contribue à cette multiplicité, alors 
$a_j> a_1$, et cet indice contribue aussi à la multiplicité de $\frac{t_1+a_1-1}{2} $. De même, si un terme unipotent 
contribue à la multiplicité de  $m-1$, il doit aussi contribuer à celle de  $\frac{t_1+a_1-1}{2} $. Ainsi la multiplicité de  $\frac{t_1+a_1-1}{2} $
est au moins égale à $3$, ce qui contredit le lemme \ref{longpsiu}.

On voit donc que l'hypothèse  $\frac{t_1+a_1-1}{2} < m-1$ implique $m-1>n-m$, ce que l'on suppose maintenant.
  Le même genre d'arguments sur le caractère infinitésimal montre que nécessairement $\frac{t_1+a_1-1}{2}=n-m$.
 Si $a(\psi_u)<2m-1$, il existe un indice $j_0\in \{1,\ldots,s\}$ tel que $\frac{t_{j_0}+a_{j_0}-1}{2}=m-1$. Comme nécessairement 
 $t_{j_0}-a_{j_0}+1\geq 0$, cet indice $j_0$ vérifie la condition (\ref{condfac}), et il est facile de voir 
 que c'est le seul. On peut donc appliquer les résultats de la proposition \ref{propred} et de la remarque \ref{rmq2red}, et conclure dans ce cas.

 Nous allons conclure grâce au lemme suivant. 
 \begin{lemme} On suppose que $a(\psi_u)=2m-1$, $m-1> n-m$ et  que $\frac{t_1+a_1-1}{2}=n-m$.  Alors $\pi(m)$ n'est pas dans $\Pi(\psi)$ 
sauf éventuellement  si $n=2(m-1)$, 
\[a(\psi_u)=n+1=2(m-1)+1=2(n-m)+3\] et 
$\psi$ contient  $\sgn_{W_\bbR}^{n/2}\boxtimes R[n+1]$.
\end{lemme}

\dem On suppose  que  $\pi(m)\in \Pi(\psi)$.  On écrit $\eta\boxtimes R[a(\psi_u)]$ 
pour le terme de $\psi_u$ qui définit $a(\psi_u)$.
 En particulier, comme $ 2m>n+1$, on a   $a(\psi_u)>n$.
 On peut donc appliquer les résultats de la section 
 \ref{siegelweil} qui assurent que la représentation $\pi(m)$ est obtenue 
 par la correspondance de Howe à partir d'une représentation d'un groupe orthogonal $\Or(p,q)$, où  la  forme quadratique à  
  pour   discriminant normalisés $\eta$ et $p+q=2n+1-a(\psi_u)=2(n-m)+2$.

On connaît donc un $K$-type de $\pi(m)$: il existe  des entiers, $x\in [0,\lfloor p/2\rfloor ],y \in [0, \lfloor q/2 \rfloor ]$ et des entiers
 strictement positifs $\alpha_j, \beta_\ell$ pour $j\in [1,x]$, $\ell\in [1,y]$ tels que $\pi(m)$ contienne le $K$-type de plus haut poids:
\begin{equation}\label{degmin}
\frac{p-q}{2}+(\alpha_1, \cdots, \alpha_x, 1, \cdots, 1, 0, \cdots, 0, -1, \cdots, -1, -\beta_{y}, \cdots,-\beta_1),\end{equation}
où le nombre de $1, 0, -1$ est comme en (\ref{phpUn}),  mais pour $G_u=\Sp(\dim(\psi_u)-1,\bbR)$. 
On sait que ce $K$-type est de \og degré minimal \fg,  
 au sens de Howe, où ce  degré est la somme des valeurs absolues des coefficients (sans le $\frac{p-q}{2}$).
 
  Pour qu'une telle 
 représentation de $K$ soit dans $\pi(m)$, il faut que chaque coefficient 
$$ \alpha_1, \cdots, \alpha_x, 1, \cdots, 1, 0, \cdots, 0, -1, \cdots, -1, -\beta_{y}, \cdots,-\beta_1,$$
  soit inférieur ou égal à $-m-\frac{p-q}{2}$.
  On a    
  \[p+q=2n+1-a(\psi_u)=2n+1-(2m-1)= 2(n+1)-2m<4m-2m=2m.\]
  On en déduit  $-m-\frac{p-q}{2}< 0$, et ainsi $x=0$, et il n'y a pas non plus de $1$ ni de $0$  dans (\ref{degmin}).
  Et le nombre de $-1$ dans (\ref{degmin}) est $n-\ell$.

La représentation $\pi(m)$ contient le $K$-type de plus haut poids $(-m\ldots,-m)$, et le degré de ce $K$-type est 
   est $n(m+\frac{p-q}{2})$, qui  doit  être supérieur 
  ou égal au degré du $K$-type de plus haut poids (\ref{degmin}), qui vaut $\sum_\ell \beta_\ell +n-\ell$.
  On obtient donc 
  \begin{equation} \label{1ineg} n(m+\frac{p-q}{2})\geq \sum_\ell \beta_\ell +n-\ell\geq n(m+\frac{p-q}{2}), \end{equation}
  ce qui force tous les coefficients de (\ref{degmin}) à être égaux à $-m-\frac{p-q}{2}$.
  En particulier tous les coefficients sont égaux.
Supposons  $p>0$. On a alors
$n+1>2(n+1)-2m=p+q>q $,
et donc $y\leq \lfloor q/2\rfloor<n$. Ainsi, il y a au moins un coefficient $-1$ dans  (\ref{degmin}), tous les $\beta_\ell$ valent $1$.
De (\ref{1ineg}), on tire à nouveau $1= m+\frac{p-q}{2}$, et comme $p+q=2(n+1)-2m$, on a $q=n$ et $p=n-2m+2$.
Or $2m>n+1$ et on a supposé $p>0$, on aboutit donc a une contradiction.
On a donc nécessairement $p=0$. Comme ci-dessus, 
\[n+1>2(n+1)-2m=p+q=q\geq \lfloor q/2 \rfloor \geq x,\]
 et  donc a nouveau $n>x$, sauf si $n=q=0$ qui est absurde. Tous les $\beta_\ell$ valent $1$ et 
$1= m+\frac{p-q}{2}=m-\frac{q}{2}$, et finalement $q=n=2(m-1)$.\qed 

\section{Démonstration du théorème \ref{mainthm}:  multiplicité un}\label{Multun}
Pour finir la démonstration du théorème, nous allons  démontrer que si $\pi(m)$ est dans le paquet d'Arthur   $\Pi(\psi)$, sa multiplicité est  un.
 Dans certains cas, on sait que les paquets $\Pi(\psi)$ ont la propriété de  multiplicité un, ce qui  règle  la question. 
 C'est le cas si $\psi=\psi_u$ (\cite{pourhowe}) ou si $a(\psi_u)>n$,  ce que nous démontrons  ci-dessous dans la section  \ref{siegelweil}.
  Comme la réduction de la section  \ref{reduc} contrôle bien les multiplicités, le seul  cas restant  
   est celui où $\dim(\psi_u)=1$, $2m> n+1$ et
   \[
\psi=(\delta_{t}\boxtimes R[n])\oplus (\sgn_{W_\bbR}^n\boxtimes R[1]),
\]
où $\frac{t+n-1}{2}= \sup(m-1,n-m)$.
 Le paquet  $\Pi(\psi)$ est décrit dans \cite{MR3} comme l'ensemble des composantes irréductibles   d'induites cohomologiques 
 à partir d'un caractère $\Lambda$ de $\U(p,q)$,  de différentielle $\lambda$ comme en (\ref{lambda1})
 et (\ref{y}), avec ici  $(p,q)$ décrivant l'ensemble des couples tels que $p+q=n$. 
 Ces induites cohomologiques sont dans le weakly fair  range, et sont  irréductibles ou nulles  grâce au critère de Bernstein-Vogan
déjà mentionné à la fin de la section \ref{npu}. 
 Mais elles sont non nulles car elles ont un bottom layer qui n'est pas vide : c'est la représentation de $K$ de plus haut poids
 $\lambda+2\delta(\fru\cap \frp)$ qui vaut ici : 
\begin{equation*}
\left(\underbrace{ \frac{t+n+1}{2}-q, \cdots, \frac{t+n+1}{2}-q}_{p}, \underbrace{-\frac{t+n+1}{2}+p, \cdots, -\frac{t+n+1}{2}+p}_{q} \right),
\end{equation*}
qui est bien un poids dominant puisque $t+n+1-(p+q)=t+1$. Pour que cette représentation soit un  $K$-type de    $\pi(m)$, 
il faut nécessairement, si $p \neq 0$, que $\frac{t+n+1}{2}-q\leq - m$.  Comme $\frac{t+n+1}{2} \geq m-1$, 
il faudrait $2m \leq q+1 \leq p+q=n$. Or par hypothèse $2m\geq n+1$ ce qui donne une  contradiction.
On a donc $p=0$, et l'on conclut grâce aux résultats de la section  \ref{71}.

\section{Démonstration du théorème \ref{Thmsnk}}\label{demThmsnk}

Supposons  que   $\sigma_{n,k}\in \Pi(\psi)$.  On connaît le paramètre de Langlands de $\sigma_{n,k}$  par la proposition  \ref{laglands2},. 
Le plus grand exposant est $(n-k)$ et il existe donc  d'après le corollaire \ref{borne} une représentation de $\SL_2(\mathbb{C})$ intervenant dans $\psi$ 
de dimension supérieure ou égale à $2n+1-2k$. Si cette représentation n'est pas dans $\psi_u$, on doit avoir
$n\geq 2n+1-2k$ ce qui force $n+1\leq 2k$,  mais cela est impossible. D'où $a(\psi_u)\geq 2(n-k)+1$. 
Le coefficient maximal dans le caractère infinitésimal est
 $n-k$ et on ne peut donc pas avoir $a( \psi_u)>2n+1-2k$, d'où le calcul de $a(\psi_u)$. 
 D'après le corollaire \ref{borne}, l'égalité force le fait que $\psi$ contient $\sgn^k\boxtimes R[2(n-k)+1]$. 
 L'énoncé de multiplicité un se démontre avec les mêmes arguments que  dans le cas des représentations $\pi_n(m)$,
 en particulier grâce aux  résultats de la section \ref{siegelweil}, et nous laissons les détails au lecteur.
 \qed

\section{Discriminant et invariant de Hasse normalisés}\label{discethasse}

Soit $F$ un corps local, et $V$ un espace vectoriel de dimension $N$ sur  $F$, muni d'une forme quadratique $Q$ non dégénérée.
On note $(.,.)_F$ le symbole de Hilbert associé au corps $F$.
Les invariants suivants sont attachés à la forme quadratique $Q$.

\

{\bf Le déterminant} : on écrit $Q(x)=\sum_{i=1}^N a_i x_i^2$ dans une base orthogonale de $V$, et l'on définit
$D(Q)$ comme la classe de $\prod_{i=1}^N a_i$ dans $F^\times/(F^\times) ^2$. Cette classe est indépendante de la base orthogonale choisie.

\ 

{\bf Le discriminant} : on multiplie le déterminant par un facteur de normalisation   en posant  $\eta(Q)=(-1)^{\frac{N(N-1)}{2}} D(Q)$.
Le discriminant est donc un élément de $F^\times/(F^\times) ^2$. Par la théorie du corps de classe, on lui associe un caractère quadratique du
groupe de Weil $W_F$ de $F$, que l'on note de la même façon $\eta(Q)$.

\ 

{\bf L'invariant de Hasse}  : on définit  $E(Q)=\prod_{1\leq i<j\leq n}(a_i,a_j)_F$
et l'on normalise en posant 
\begin{equation} \label{hasse}
\epsilon_\delta(Q) = (-\delta,\eta(Q))_F \;    (-1,D(Q))_F^{ \frac{N(N-1)}{2} }  (-1,-1)_F^{\lfloor  \frac{\lfloor   \frac{N}{2}\rfloor +1 }{2} \rfloor } \; E(Q). \end{equation} 
Ici $\delta$ est un élément de $F^\times  /(F^\times)^2$ dont la signification apparaîtra plus tard.
On a donc $\epsilon_\delta(Q)\in \{\pm 1\}$.

\ 

Le discriminant et l'invariant de Hasse ainsi normalisés sont invariants par ajout d'un plan hyperbolique.

\

\begin{prop}
Dans le cas où $F$ est le corps des réels $\bbR$, ces invariants sont tous déterminés par la signature $(p,q)$ de la forme quadratique
$Q$. On note alors $D(p,q)$, $\eta(p,q)$, $\epsilon(p,q)$. On a :

--- si $n=p+q$ est pair, 
\[D(p,q)=(-1)^q,\quad  \eta(p,q)=(-1)^{\frac{p-q}{2}},\quad  \epsilon_\delta(p,q)=(-1)^{\lfloor \frac{\delta(p-q)}{4} \rfloor}.\]

--- si $n=p+q$ est impair, 
\[ D(p,q)=(-1)^q,\quad  \eta(p,q)=(-1)^{\frac{p-q-1}{2}}, \quad \epsilon_\delta(p,q)=(-1)^{\lfloor \frac{\delta(p-q-1)}{4} \rfloor}.\]

\end{prop}

\dem  Calculons ces invariants  lorsque $p+q$ est pair,
qui est le seul cas dont on se servira, en laissant le cas $p+q$ impair au lecteur.
Comme la formule (\ref{hasse}) est invariante par ajout de plans hyperboliques et qu'il en est de même de la formule proposée,
 il suffit de vérifier l'égalité pour les signatures $(p,q)$ tel que $pq=0$. Remarquons d'abord que
  $\eta(p,q)=(-1)^{(p-q)/2}$, $D(p,q)=(-1)^q$, $E(p,q)=(-1)^{q(q-1)/2}$. Si $q=0$, (\ref{hasse}) se simplifie en 
\[ \epsilon_\delta(p,q)=(-1)^{\frac{(1+\delta)}{2}  \frac{p}{2} + \lfloor  \frac{\frac{ p}{2}  +1}{2}\rfloor}\]
Si $p/2$ est pair, cela vaut $(-1)^{p/4}$ quelque soit la valeur de $\delta$ et si $p/2$ est impair cela vaut $(-1)^{\lfloor \frac{p}{4}\rfloor}$ 
si $\delta=+1$ et $(-1)^{\lfloor \frac{p}{4}\rfloor+1}$ si $\delta=-1$; mais comme $\frac{p}{2}$ est impair, 
\[ \lfloor \frac{p}{4}\rfloor+1=\frac{p}{4}+\frac{1}{2}=
-\frac{p}{4}+\frac{p}{2}+\frac{1}{2} \equiv -\frac{p}{4}-\frac{1}{2} \mod 2,\]
 d'où le résultat annoncé.
Si $p=0$, on utilise l'égalité (pour $n$ pair)\[
\epsilon_{\delta}(p,q)\epsilon_{\delta}(q,p)=E(p,q)E(q,p)=(-1)^{\lfloor \frac{q}{2}\rfloor +\lfloor \frac{p}{2}\rfloor}.
\]
D'où pour $p=0$, $\epsilon_{\delta}(p,q)=(-1)^{q/2}\epsilon_{\delta}(q,p)=(-1)^{q/2+\lfloor \delta q/4\rfloor }$. 
Si $q/2$ est pair, cette formule est indépendante de $\delta$ et vaut bien $(-1)^{-q/4}=(-1)^{q/4}$ 
et si $q/2$ est impair la formule dépend de $\delta$ et est celle annoncée   :  $(-1)^{\lfloor -\delta q/4\rfloor }$.
\qed 

\

Remarquons qu'en dimension impaire, on fait clairement un choix arbitraire, $\eta(2,1)=1$, $\eta(1,2)=-1$.

\section{Caractères de $\Or(V,Q)$ et image par la correspondance de Howe  \label{caracO}}

\subsection{Correspondance de Howe locale } \label{Howeloc}
On se place d'abord sur un corps local $F$ de caractéristique $0$.  On fixe un caractère additif $\psi$ de $F$.
Nous n'allons utiliser la correspondance de Howe que dans le cas d'une paire duale orthogonale paire/symplectique.
Soit $V$ un espace  vectoriel de dimension $N=2m$ sur  $F$, muni d'une forme quadratique $Q$ non dégénérée.
Soit $W$ un espace vectoriel symplectique sur $F$ de dimension $2n$.
Les conjectures de Howe étant maintenant démontrées en toute généralité (\cite{Ho}, \cite{Wald89}, \cite{GanTakeda}), 
pour toute représentation irréductible $\pi$ de $\Or(V,Q)$, on note 
$\theta^\psi_{W,V}(\pi)$ son image par la correspondance de Howe, qui est une représentation irréductible 
de $\Sp(W)$, ou bien $0$. Comme tous les espaces symplectiques sur $F$ de dimension $2n$ sont isomorphes,
on note en fait plutôt $\theta^\psi_{n,V}(\pi)$ pour $\theta^\psi_{W,V}(\pi)$ et $\Sp(2n,F)$ pour  $\Sp(W)$.

Rappelons quelques  propriétés  de la correspondance de Howe. Pour toute représentation irréductible $\pi$ de $\Or(V,Q)$,

(1)  si $\theta^\psi_{n,V}(\pi)\neq 0$, alors pour tout $n'\geq n$, $\theta^\psi_{n',V}(\pi)\neq 0$,

(2) si $n$ est suffisament grand, alors  $\theta^\psi_{n,V}(\pi)\neq 0$.

On en déduit que $n_0(\pi)=\min_{n\in \bbN} \{n\vert \, \theta^\psi_{n,V}(\pi)\neq 0\} $ existe, et on appelle cet entier la première occurence
de $\pi$ dans la correspondance.

(3) (Loi de conservation)  $n_0(\pi)+n_0(\pi\otimes \det)=2m$, 

(4) $n_0(\Triv_{\Or(V,Q)})=0$, $n_0(\det_{\Or(V,Q)})=2m$.

Les lois de conservation ont été conjecturées dans \cite{KR94} et démontrées dans \cite{SunZhu}.

\

\subsection{Correspondance de Howe globale}\label{Howeglob}
On se place maintenant sur un corps de nombre  $k$. Si  $v$ est une  place de $k$, on note $k_v$ le corps local correspondant, 
et l'on ajoute un indice $v$  aux notations pour les
 objets sur  $k_v$  lorsque ceux-ci ont  été définis précédemment sur un corps local   quelconque.  Les objets définis sur le corps
de nombre $k$ ou sur l'anneau des adèles $\bbA_k$ seront eux notés en caractère gras. 
 On fixe  un espace vectoriel $\caV$ sur $k$ de dimension $N=2m$ muni d'une forme quadratique $\caQ$
et l'on pose $(\bm V,\bm Q)= \prod_v(\caV_v,\caQ_v)$. On fixe un caractère additif $\bm \psi=(\psi_v)_v$ de $\bbA_k/k$.
Aux places réelles qui nous intéressent, on supposera implicitement que $\psi_v$ est le caractère $\psi_{\bbR,1}$ de la remarque 
\ref{rmqpqqp}.

Pour toute représentation automorphe cuspidale $\bm \pi$ de $\Or(\bm V,\bm Q)$, on note $\theta^{\bm \psi}_{n,\caV}(\bm\pi)$
son relèvement par séries thêta dans $\Sp(2n,\bbA_k)$. C'est une représentation automorphe de $\Sp(2n,\bbA_k)$, ou bien $0$.

L'analogue évident des propriétés (1) et (2) ci-dessus est valable dans le cas global, notons (1')
 et (2') ces propriétés.  On peut alors définir aussi la première occurence $n_0(\bm \pi)$ comme ci-dessus. On a alors aussi les 
 propriétés suivantes
 
 (3') $\theta^{\bm \psi}_{n_0(\bm \pi),\caV}(\bm\pi)$ est cuspidale.
 
 (4') $n_0(\bm \pi)\leq 2m$.
 
 (5') Si $n> n_0(\bm \pi)$ et si $n_0(\bm \pi)+n>2m-1$, alors  $\theta^{\bm \psi}_{n,\caV}(\bm\pi)$ est de carré intégrable
 et si $n_0(\bm \pi)+n=2m-1$, c'est une induite de carré intégrable.
 
En effet, d'après la théorie de Rallis ({\sl cf.} \cite{Rallis}), pour $n\geq n_0(\bm \pi)$,  $\theta^{\bm \psi}_{n,\caV}(\bm\pi)$ 
n'a qu'un seul terme constant cuspidal non nul, celui-ci étant relatif au parabolique de facteur de Levi 
\[  \underbrace{\GL_1\times\cdots \times \GL_1}_{n-n_0(\bm \pi)}\times \Sp_{2n_0(\bm \pi)}.  \]
Ce terme constant vit dans l'induite 
\[\eta(\caQ)  \vert .\vert^{-\left(n-m \right)} \times \cdots \times \eta(\caQ)  \vert .\vert^{-\left(n_0(\bm \pi)+1-m \right)}\times \theta^{\bm \psi}_{
n_0(\bm \pi),\caV}(\bm\pi) ,\]
et l'assertion découle alors de  \cite{MWDec}.  Ici $\eta(\caQ)$ est le discriminant de la forme quadratique $\caQ$, 
qui est définie comme dans le cas local et est un élément de $k^\times/(k^\times)^2$. Par la théorie du corps de classe, on le voit aussi 
comme un caractère quadratique de $k^\times$.

Supposons que la forme quadratique $\caQ$ soit anisotrope à au moins une place, de sorte que le groupe orthogonal correspondant 
soit compact. Alors la caractère trivial de $\Or(\bm V,\bm Q)$ est automorphe cuspidal et son relèvement par série thêta
existe pour tout $n$, c'est-à-dire : 

(6') $n_0(\Triv_{\Or(\bm V,\bm Q)})=0$.

D'après \cite{Rallis}, on a de plus : 

(7')  pour tout  $n\geq 0$,   $\theta^{\bm \psi}_{n+1,\caV}(\Triv_{\Or(\bm V,\bm Q)})$ se réalise 
 comme quotient de l'induite 
 \[ \eta(\caQ)  \vert .\vert^{-\left(n+1-m \right)}  \times \theta^{\bm \psi}_{n,\caV}( \Triv_{\Or(\bm V,\bm Q)}) ,\]
 relative au  sous-groupe parabolique de facteur de Levi 
$  \GL_1\times  \Sp_{2n} $ de $\Sp_{2n+2}$.

\ 

\subsection{Caractères de $\Or(V,Q)$}\label{sgnd}

On se place sur un corps local $F$ de caractéristique $0$ et  $V$ est un espace vectoriel de dimension $N=2m$ sur  $F$, muni d'une forme quadratique $Q$ 
non dégénérée. En général, si $\pi$ est un représentation de $\SO(V,Q)$ stable par l'automorphisme extérieur de ce groupe, 
elle se relève de deux façons à $\Or(V,Q)$, et il n'y a pas de manière canonique de choisir un de ces deux relèvements.
Pour les caractères de $\SO(V,Q)$, nous allons fixer un relèvement. Evidemment, comme relèvement du caractère trivial
de $\SO(V,Q)$, on choisit le caractère trivial de $\Or(V,Q)$.

La norme spinorielle est un morphisme de groupes $NS_Q: \, \Or(V,Q)\rightarrow  F^\times $. 
Elle dépend de la forme quadratique $Q$, et pas simplement de $\Or(V,Q)$ 
(en effet $\Or(V,Q)=\Or(V,\lambda Q)$ pour tout $\lambda\in F^\times$, mais $NS_{\lambda Q}\neq NS_Q$ en général).
Les caractères de $\Or(V,Q)$ sont de la forme 
$ (\eta \circ NS_Q )\otimes \textstyle  \det^\tau$
où  $\eta$ est  un caractère quadratique de $F^\times$ et $\tau\in \bbZ/2\bbZ$. Un calcul immédiat montre que 
\begin{equation} \label{moinsQ}
 \eta \circ NS_{-Q}  =  (\eta \circ NS_Q )\otimes (\eta\circ \textstyle  \det). \end{equation} 

Pour   $\delta\in F^\times /(F^\times)^2$ et $\eta$ caractère quadratique de $F^\times$, posons  
\begin{equation}\label{defetadelta}\eta_\delta(Q) =  \eta \circ NS_{-\delta \eta(Q) Q} . \end{equation}
C'est un caractère de $\Or(V,Q)$, bien défini car ne dépendant pas du choix d'un représentant de $\delta\eta(Q)$ dans $F^\times$ .
 Tous les caractères de $\Or(V,Q)$, $\delta$ étant fixé, sont de la forme 
\begin{equation}\label{touscar}  \eta_\delta(Q) \otimes \textstyle \det^\tau   \end{equation}
où $\eta$ décrit les caractères quadratiques de $F^\times$ et $\tau$ décrit $\bbZ/2\bbZ$.

\ 

\begin{rmq}
Il convient de bien distinguer des notations qui sont proches : $\eta$ est un caractère quadratique de $F^\times$, $\eta(Q)$
est le discriminant de la forme quadratique $Q$, et donc un élément de $F^\times/(F^\times)^2$, et $\eta_\delta(Q)$
est un caractère du groupe $\Or(V,Q)$. La théorie du corps de classe nous permet d'identifier 
caractères quadratiques de $F^\times$ (ou de $W_F$) et éléments de $F^\times/(F^\times)^2$, et selon le contexte
on voit $\eta$ et $\eta(Q)$ comme l'un ou l'autre de ces types d'objets.
\end{rmq}

 Dans le cas global, sur un corps de nombres $k$,  on adapte les notations comme en \ref{Howeglob}.
  On fixe  un espace vectoriel $\caV$ sur $k$ de dimension $N=2m$ muni d'une forme quadratique $\caQ$
et l'on pose $(\bm V,\bm Q)= \prod_v(\caV_v,\caQ_v)$. On fixe aussi   
  et  un 
 caractère adélique $\bm \eta=(\eta_v)_v: \bbA_k^\times /k^\times \to \bbC^\times$.
On définit alors le caractère 
\begin{equation}\label{defetadeltaglob}  \bm \eta_\delta(\bm Q) = \prod_v   \eta_{\delta}(\caQ_v):
\Or(\bm V,\bm Q)\longrightarrow \{\pm 1\}. \end{equation}
C'est un caractère automorphe.

\

On revient au cas local et  on suppose que $F$ est un corps $p$-adique.

\begin{prop} \label{occurence} Soit $V$ un espace vectoriel de dimension $N=2m$ sur  un  corps local non archimédien 
$F$ de caractéristique $0$,  muni d'une forme quadratique $Q$ non dégénérée. Soient  
$\eta$ un caractère quadratique de $F^\times$, $\delta$ un élément de $F^\times/(F^\times )^2$ et $\eta_\delta(Q)$ 
le caractère de $\Or(V,Q)$ défini en (\ref{defetadelta}).

(i)
On suppose que $\eta \notin\{\eta(Q),1\}$. Alors la première occurence de
 $\eta_\delta(Q) \otimes \det^\tau$ dans la correspondance de Howe est pour $\Sp(2m,F)$, en particulier cela ne dépend pas de $\tau$.

(ii) On suppose que $\eta=1$, la première occurence du caractère trivial de $\Or(V,Q)$  est pour
 $\Sp(0,F)=\{1\}$ et la première occurence du caractère déterminant est pour $\Sp(4m,F)$.

(iii) Quel que soit $\eta$ l'un au moins des deux caractères $\eta_\delta(Q) \otimes \det^\tau$, $\tau=0,1$,  a une image dans $\Sp(2m,F)$.
\end{prop}

%Nous allons utiliser le résultat suivant sur la filtration de Kudla (voir  \cite{pourkudla}). On utilise les notations de \cite{elementaire}
 %pour les modules de Jacquet
%\begin{lemme} \label{KudlaFilt}
%Avec les notations qui précèdent, notons $n_0=n_0(\eta_\delta(Q)\otimes \det^\tau)$ la première occurence de
 %$\eta_\delta(Q) \otimes \det^\tau$ dans la correspondance de Howe. 
% Supposons $n\geq \max(n_0,m)$, de sorte que   
 %$\theta^\psi_{n,V}(\eta_\delta(Q)\otimes \det^\tau)$ soit non nul et supposons que 
 %$\mathrm{Jac}_{\eta(Q)\vert\, \vert ^{n-m}} (\theta^\psi_{n,V}(\eta_\delta(Q)\otimes \det^\tau))$ soit non nul.
 %Alors 
 %\[\textstyle \mathrm{Jac}_{\eta(Q)\vert\, \vert ^{n-m}} (\theta^\psi_{n,V}(\eta_\delta(Q)\otimes \det^\tau))=\theta^\psi_{n-1,V}(\eta_\delta(Q)\otimes \det^\tau),\]
%et en particulier $\theta^\psi_{n,V}(\eta_\delta(Q)\otimes \det^\tau) $ est non nul {\sl i.e.} $n_0\geq n-1$.
%\end{lemme}
%{\color{blue} est-ce que c'est correct maintenant? j'ai ajouté $n\geq m$ comme hypothèse}

\underline {\sl Démonstration de la proposition}. Remarquons que $(ii)$ est la propriété (4) de \ref{Howeloc}, et que $(iii)$ découle directement de la propriété (5).
Montrons (i)  par des arguments globaux.  
On considère un corps de nombre $k$ et un espace quadratique $(\caV,\caQ)$ sur $k$ avec une  place $v_1$ de $k$  telle 
 que $k_{v_1}=F$,  $(\caV_{v_1},\caQ_{v_1})=(V,Q)$.
On suppose que la forme quadratique est anisotrope en au moins une place $v_2$, et donc que le groupe orthogonal correspondant est compact
(une telle place est forcément archimédienne si $N=2m>4$, soit $m>2$).
On fixe aussi un caractère quadratique adélique $\bm \eta=(\eta_v)_v $ de $\bbA_k^\times /k^\times$ tel que et $\eta_{v_1}=\eta$  
et une famille $\bm \tau=(\tau_v)_v$ d'éléments de $\{0,1\}$ avec un nombre fini et pair de $v$ tels que $\tau_v=1$ et $\tau_{v_1}=\tau$. 
On pose $\det^{\bm \tau}=\prod_{v\vert \tau_v=1} \det_v$.

 Le caractère $\bm \eta_\delta(\bm Q)\otimes \det^{\bm \tau}$ de $\Or(\bm V,\bm Q)$  est alors automorphe
 cuspidal  (grâce à l'anisotropie en $v_2$) et
on peut appliquer les résultats énoncés en \ref{Howeglob}. En particulier ce caractère 
  admet un relèvement par séries thêta  $\theta^{\bm \psi}_{n,\caV}(\bm \eta_\delta(\bm Q)\otimes \det^{\bm \tau})$
     non nul  si $n$ est assez grand, 
et  ce relèvement est une représentation automorphe de carré intégrable (propriétés (2') et (5') de \ref{Howeglob}).  

On connait les composantes locales de  $\theta^{\bm \psi}_{n,\caV}(\bm \eta_\delta(\bm Q)\otimes \det^{\bm \tau})$
 aux places finies non ramifiées, et ceci détermine un paramètre d'Arthur global \cite{Art13}.
Pour un tel $n$, on en déduit    ({\sl cf. }\cite{pourkudla})   localement que $\theta^{\psi}_{n,V}(\eta_\delta(Q) \otimes \det^\tau)\neq 0$
et que cette représentation appartient  au paquet  d'Arthur de paramètre  :
\[ \psi:= (\eta \boxtimes R[1] )\oplus (\eta\eta(Q)  \boxtimes R[2m-1])\oplus (\eta(Q)\boxtimes R[2(n-m)+1]).
\]

Mais les représentations associées à ces paramètres  sont bien connues, elles ont été étudiées dans  \cite{elementaire}. En particulier l'hypothèse 
sur $\eta$ assure que  le troisième bloc de $\psi$ est différent des deux premiers, et il est alors montré en {\sl loc. cit.} que si 
$\theta^{\psi}_{n,V}(\eta_\delta(Q) \otimes \det^\tau)$  apparaît dans $\Pi(\psi)$ alors, pour tout $n'\geq m$, il existe une représentation $\rho$ dans $\Pi(\psi')$ avec 
\[ \psi'=( \eta \boxtimes R[1]) \oplus (\eta\eta(Q) \boxtimes R[2m-1])\oplus( \eta(Q)\boxtimes R[2(n'-m)+1]),\]
telle que $\theta^{\psi}_{n,V}(\eta_\delta(Q) \otimes \det^\tau)$ soit l'unique quotient irréductible de l'induite:
\[
 \eta(Q) |.|^{n-m} \times   \eta(Q) |.|^{n-m-1} \times\ldots \times   \eta(Q) \vert .\vert ^{n'-m+1} \times \rho
\]
Un résultat sur la filtration de Kudla de  \cite{pourkudla} assure que $\rho$ est égale à 
$ \theta^{\psi}_{n',V}( \eta_\delta(Q) \otimes  \det^\tau)$. En particulier, on peut 
prendre $n'=m$ et l'on voit  la première occurence $n_0(\eta_\delta(Q) \otimes \det^\tau)$ du caractère est certainement inférieure ou égale à $m$.
 Mais comme ceci  est vrai pour les deux valeurs de
 $\tau$ la loi de conservation (4) de \ref{Howeloc}  assure que  $n_0(\eta_\delta(Q) \otimes \det^\tau)=m$.
\qed 

\

On se place maintenant dans le cas  $F=\bbR$ avec $Q$ est de signature $(p,q)$ avec $p+q=N=2m$ pair.
Les caractères du groupe $\Or(p,q)$ se factorisent par le groupe  de ses composantes connexes
qui est isomorphe à $\bbZ/2\bbZ$ si $pq=0$ et $\bbZ/2\bbZ\times \bbZ/\bbZ$ si $pq\neq 0$.
Dans le premier cas, il y a donc deux caractères, le caractère trivial et le déterminant, et dans le second, il faut en ajouter $2$.
Tous ces caractères sont déterminés par leur restriction au sous-groupe compact maximal $\Or(p,0)\times \Or(0,q)$ de $\Or(p,q)$.
On note $\eta_\delta(p,q)$ plutôt que $\eta_\delta(Q)$ pour le caractère de $\Or(p,q)$ construit en (\ref{defetadelta}).

Il y a deux caractères quadratiques de $\bbR^\times$  le caractère trivial et le caractère signe, notés $1_{\bbR^\times}$ et  $\sgn_{\bbR^\times}$
ou simplement $1$ et $\sgn$.
Si l'on prend $\eta=1$, alors $\eta_\delta(p,q)=\Triv_{\Or(p,q)}$  pour tout $\delta$ et l'on 
 obtient ainsi en (\ref{touscar}) les deux caractères $\Triv_{\Or(p,q)}$ et $\det_{\Or(p,q)}$ de $\Or(p,q)$.
Si l'on prend $\eta=\sgn$  et $\delta=1$, alors $\eta_\delta(p,q)=\sgn_1$ admet comme 
restriction à $\Or(p,0)\times \Or(0,q)$ le caractère 
\[   \textstyle \det_{\Or(p,0)}^{\frac{p-q}{2}+1}\boxtimes  \det_{\Or(0,q)}^{\frac{p-q}{2}} .\]
Le quatrième et dernier caractère est alors $\eta_\delta(p,q)\otimes \det =\sgn_1\otimes \det$ qui 
admet comme 
restriction à $\Or(p,0)\times \Or(0,q)$ le caractère 
\[   \textstyle \det_{\Or(p,0)}^{\frac{p-q}{2}}\boxtimes  \det_{\Or(0,q)}^{\frac{p-q}{2}+1} .\]
Remarquons que lorsque $pq=0$, les caractères $\eta_\delta(p,q)$ et $\eta_\delta(p,q)\otimes \det $ sont bien définis, on retrouve simplement
les caractères $\Triv_{\Or(p,q)}$ et $\det_{\Or(p,q)}$ (pas nécessairement dans cet ordre). Par exemple pour $\Or(2,0)$ et $\Or(0,2)$ :
\[\sgn_{1} = \Triv_{\Or(2,0)} ,\quad  \sgn_{-1} = \textstyle \det_{\Or(2,0)} ,      \]
\[\sgn_{1} = \textstyle \det_{\Or(0,2)} ,\quad  \sgn_{-1} = \Triv_{\Or(0,2)}.       \]

On détermine $\eta_\delta(p,q)$ lorsque $\eta=\sgn$ et $\delta=-1$ grâce à (\ref{moinsQ}).

Déterminons l'image par la correspondance de Howe de ces caractères. La correspondance est ici fixée par le choix du caractère additif $\psi_{\bbR,1}$
en (\ref{caraddR}).

Le caractère trivial $\Triv_{\Or(p,q)}$ a une image dans $\Sp(2n,\bbR)$ pour tout $n$. Celle-ci contient le $\U(n)$-type :
\[ \frac{p-q}{2}+( \underbrace{0,\ldots,0}_{n}  ) . \]

Le caractère $\det_{\Or(p,q)}$ a une image dans $\Sp(2n,\bbR)$ pour  $n\geq p+q=N=2m$. Celle-ci contient le $\U(n)$-type :
\[ \frac{p-q}{2}+( \underbrace{1,\ldots,1}_{p} ,\underbrace{0,\ldots,0}_{n-p-q},  \underbrace{-1,\ldots,-1}_{q} ) . \]

Si $\frac{p-q}{2}$ est impair, la restriction  du caractère $\sgn_{1}$ de $\Or(p,q)$
au sous-groupe compact maximal $\Or(p)\times \Or(q)$ est $\Triv_{\Or(p)}\boxtimes  \det_{\Or(q)}$.
L'image de ce caractère par la correspondance de Howe dans le groupe $\Sp(2n,\bbR)$ est 
non nulle si $n\geq q$ et alors cette image contient le $\U(n)$-type :
\[\frac{p-q}{2}+(\underbrace{0,\ldots,0}_{n-q}, \underbrace{-1,\ldots,-1}_{q} ). \]
Toujours avec  $\frac{p-q}{2}$  impair, la restriction  du caractère $\sgn_{1}\otimes \det $ de $\Or(p,q)$
au sous-groupe compact maximal $\Or(p)\times \Or(q)$ est $\det_{\Or(p)}\boxtimes  \Triv_{\Or(q)}$.
L'image de ce caractère par la correspondance de Howe dans le groupe $\Sp(2n,\bbR)$ est 
non nulle si $n\geq p$ et alors cette image contient le $\U(n)$-type :
\[\frac{p-q}{2}+(\underbrace{1,\ldots,1}_{p}, \underbrace{0,\ldots,0}_{n-p} ). \]

Si $\frac{p-q}{2}$ est pair, la restriction  du caractère $\sgn_{1}$ de $\Or(p,q)$
au sous-groupe compact maximal $\Or(p)\times \Or(q)$ est $\det_{\Or(p)}\boxtimes  \Triv_{\Or(q)}$
L'image de ce caractère par la correspondance de Howe dans le groupe $\Sp(2n,\bbR)$ est 
non nulle si $n\geq p$ et alors cette image contient le $\U(n)$-type :
\[\frac{p-q}{2}+(\underbrace{1,\ldots,1}_{p}, \underbrace{0,\ldots,0}_{n-p} ). \]
Toujours avec  $\frac{p-q}{2}$  pair, la restriction  du caractère $\sgn_{1}\otimes \det $ de $\Or(p,q)$
au sous-groupe compact maximal $\Or(p)\times \Or(q)$ est $\Triv_{\Or(p)}\boxtimes  \det_{\Or(q)}$.
L'image de ce caractère par la correspondance de Howe dans le groupe $\Sp(2n,\bbR)$ est 
non nulle si $n\geq q$ et alors cette image contient le $\U(n)$-type :
\[\frac{p-q}{2}+(\underbrace{0,\ldots,0}_{n-q}, \underbrace{-1,\ldots,-1}_{q} ). \]

Si $p=0$, l'image du caractère trivial de $\Or(0,q)=\Or(0,2m)$ est la représentation unitaire de plus haut poids scalaire $-m$
notée $\pi_{n}(m)$ dans l'article et l'image du déterminant est la représentation notée $\sigma_{n,m}$ 
(la condition d'existence est $2m\leq n$). Si l'on prend $q=0$ et $p=N=2m$, on obtient
des représentations unitaires de plus bas poids qui sont les contragrédientes de celles-ci.

Si $p=q=1$, on a une image dans $\SL(2,\bbR)$
de plus haut  poids $1$ : c'est encore  la limite de série discrète holomorphe (avec les conventions usuelles, pas celles de la première partie de cet article)
 de $\SL(2,\bbR)$.

\begin{rmq}
%On aurait pu changer $\delta$ en $-\delta$ dans la définition $\eta_{\delta}$, cela n'a en fait pas grande importance.
%Le choix d'un élément $\delta\in F^\times /(F^\times )^2$ permet  de paramétrer les caractères additifs de $F$
%lorsqu'on a fixé un caractère de référence. Changer $\delta$ en $-\delta$ revient à changer le caractère de reférence en son inverse. 
%Ce caractère de référence intervient dans la définition de la correspondance de Howe.
%{\color{red}  c'est pas très clair...}

Si on prend $\delta=-1$, on obtient les mêmes énoncés en échangeant les rôles de $p$ et $q$, 
\end{rmq}

\subsection{$A$-paquets}

Nous reprenons les notations de la section  précédente pour un corps local $F$. Lorsque
$\theta^\psi_{n,V}( \eta_\delta(Q)\otimes \det^\tau)$  est non nulle, nous allons maintenant donner un paquet d'Arthur la contenant.

\begin{thm} \label{lemp} On suppose   $n \geq  2m-1$. La représentation $\theta^\psi_{n,V}( \eta_\delta(Q)\otimes \det^\tau)$, 
si elle n'est pas nulle, est dans le $A$-paquet de paramètre 
\[\psi=
(\eta\boxtimes R[1]) \oplus (\eta \eta(Q)\boxtimes R[2m-1]) \oplus (\eta(Q)  \boxtimes R[2(n-m)+1]). \]
\end{thm}

\dem   Remarquons que si   $F$ est  un corps  $p$-adique et que $\eta\notin\{\eta(Q),1\}$, on a 
 déjà montré le résultat  dans la démonstration de la proposition 
 \ref{occurence}  (i) (avec même une meilleure hypothèse sur $n$).
Nous allons donner une démonstration générale avec  le même genre d'argument global, mais plus fin.

On considère un corps de nombre $k$, un espace quadratique $(\caV,\caQ)$ sur $k$ 
de dimension $2m$, un caractère adélique quadratique $\bm \eta=(\eta_v)_v$ de $\bbA_k^\times /k^\times$
 et une famille $\bm \tau=(\tau_v)_v$ d'éléments de $\{0,1\}$ avec un nombre fini et pair de $v$ tels que $\tau_v=1$ avec les propriétés
 suivantes : 
 
 - en une place $v_0$, $k_{v_0}=F$,  $(\caV_{v_0},\caQ_{v_0})=(V,Q)$,  $\eta_{v_0}=\eta$, et $\tau_{v_0}=\tau$.
 
 - en toute place $v\notin\{ v_0,v_1\}$, on choisit $\tau_v$ de sorte que la première occurence de 
 $\eta_\delta(\caQ_v)\otimes \det^{\tau_v}$ dans la correspondance de Howe soit au moins égale à $m$
 (ce qui est possible d'après la proposition \ref{occurence}).
 
 - en une place $v_1$ réelle, $(\caV_{v_1},\caQ_{v_1})$ est de signature $(m,m)$ et $\eta_{v_1}=\sgn$.
 Nous avons vu dans la section précédente que la première occurence du caractère 
 $\eta_\delta(\caQ_{v_1})\otimes \det^\tau$ est pour $n_0=m$, et ceci pour tout $\tau\in \{0,1\}$. On choisit alors $\tau_{v_1}=\{0,1\}$ tel
 que le nombre de $\tau_v$ égaux à $1$ soit pair.
 
-  en   une place $v_2$, la forme quadratique est anisotrope, et donc que le groupe orthogonal correspondant est compact
(une telle place est forcément archimédienne si $N=2m>4$, soit $m>2$).
 
 Le caractère $\bm \eta_\delta(\bm Q)\otimes \det^{\bm \tau}$ de $\Or(\bm V,\bm Q)$  est alors automorphe
 cuspidal  (grâce à l'anisotropie en $v_2$) et
on peut appliquer les résultats énoncés en \ref{Howeglob}.
Localement, en toute place $v$ et pour tout $n\geq 2m-1$, on a $\theta^{\psi_v}_{n,\caV_v}(\eta_\delta(\caQ_v)\otimes \det^{\tau_v})\neq 0$.
Le résultat principal de \cite{yam}
(outre la non annulation locale, il y a une condition sur les fonction $L$ qui est vérifiée ici)
nous dit alors que le relèvement global $\theta^{\bm \psi}_{n,\caV}( \bm \eta_\delta(\bm Q)\otimes \det^{\bm \tau})$
est non nul
et  ce relèvement est une représentation automorphe et  de carré intégrable (propriétés (2') et (5') de \ref{Howeglob}).  
En effet, l'hypothèse faite à la place $v_1$ implique que la première occurence d'un relèvement global est au moins 
égale à $m$.  La théorie d'Arthur \cite{Art13}
s'applique et l'on en déduit un paramètre d'Arthur global déterminé par la correspondance de Howe aux place non ramifiées. 
Localement, ceci nous donne l'appartenance de $\theta^\psi_{n,V}( \eta_\delta(Q)\otimes \det^\tau)$
au paquet voulu.
\qed

\begin{rmq}\label{blocsegaux1}
On peut se demander à quelle condition on a  $\theta^\psi_{n,V}(\eta_\delta(Q)\otimes \det^\tau)\neq 0$ pour $n\geq 2m-1$.
La première occurence par la correspondance de Howe d'un caractère de $\Or(V,Q)$ 
est toujours pour $n_0\leq 2m$, et n'est égale à $2m$ que pour le caractère déterminant.
On a donc $ \theta^\psi_{n,V}(\eta_\delta(Q)\otimes \det^\tau)\neq0$ pour $n\geq 2m-1$ sauf si  $n=2m-1$,
$\eta=1$ et $\tau=1$.
Le paramètre d'Arthur du théorème est alors 
\[\psi= (1 \boxtimes R[1]) \oplus ( \eta(Q)\boxtimes R[2m-1]) \oplus (\eta(Q)  \boxtimes R[2m-1]), \]
et l'on constate que le deuxième et troisième blocs sont égaux.
\end{rmq}

\begin{rmq}\label{blocsegaux2} Une autre possibilité d'avoir deux blocs égaux dans le paramètre d'Arthur du lemme 
précédent est lorsque $m=1$ et $\eta(Q)=1$. On a alors 
\[\psi=
(\eta\boxtimes R[1]) \oplus (\eta \boxtimes R[1]) \oplus (1  \boxtimes R[2n-1]) ,\]
En dimension $2$, il n'existe pas de forme quadratique $Q$ de discriminant trivial et d'invariant de Hasse non trivial.
\end{rmq}

\section{Calcul des  $\rho_\pi$}\label{calcrhopi}

On se place sous les hypothèses du théorème \ref{lemp}, dont l'énoncé fait intervenir un paramètre d'Arthur local unipotent
\[\psi=
(\eta\otimes R[1]) \boxplus (\eta \eta(Q)\boxtimes R[2m-1]) \boxplus (\eta(Q)  \otimes R[2(n-m)+1]). \]

Le groupe $A(\psi)$ s'identifie naturellement à $\mathbf S(\Or(1)\times\Or(1)\times\Or(1))$, sauf dans les quelques cas
suivants   où l'on a des multiplicités dans $\psi$ : 

$\eta=1$, $n=2m-1$, $\eta(Q)\neq 1$ ou $m\geq 2$, $A(\psi)=\mathbf S(\Or(1)\times\Or(1))$,

$\eta=\eta(Q)$, $n=m$, $\eta\neq 1$ ou $m\geq 2$,  $A(\psi)=\mathbf S(\Or(1)\times\Or(1))$.

$\eta=\eta(Q)= 1$, $n=2m-1$,  $m=1$, $A(\psi)=\{1\}$,

\
En dehors de ces cas, on identifie $A(\psi)$ au groupe
\[\{(\epsilon_1,\epsilon_2,\epsilon_3)\in \{\pm 1\}^3\vert \, \epsilon_1\epsilon_2\epsilon_3=1\}.\]
Dans tous les cas, on identifie $\widehat {A(\psi)}$ à l'ensemble des fonctions 
sur les blocs de $\psi$, à valeurs dans $\{\pm 1\}$, égales sur les blocs égaux, modulo changement de signe. 
On peut voir ces fonctions
comme un triplet   $(\epsilon_1,\epsilon_2,\epsilon_3)\in \{\pm 1\}^3$ défini à un signe global près, et l'on peut choisir  
un représentant naturel avec la condition $\epsilon_1\epsilon_2\epsilon_3=1$.

Le théorème \ref{lemp} affirme que si $n\geq 2m-1$ et si  la représentation $\theta^\psi_{n,V}( \eta_\delta(Q)\otimes \det^\tau)$ est  non nulle, alors elle est 
dans le paquet $\Pi(\psi)$. On rappelle ({\sl cf.} section  \ref{ArtparSp}) 
qu'Arthur  définit  alors une représentation $\rho_{\theta^\psi_{n,V}( \eta_\delta(Q)\otimes \det^\tau)}$ de $A(\psi)$, et 
l'on sait d'après la propriété de multiplicité un de ces paquets que cette représentation est de dimension un. On l'identifie comme
 ci-dessus à un  triplet de signes $\pm 1$ que l'on note   
\begin{equation}\label{eppi}  (\epsilon_1(n,\eta_\delta(Q),\tau) ,\epsilon_2(n,\eta_\delta(Q),\tau),\epsilon_3(n,\eta_\delta(Q),\tau)) .\end{equation}
Rappelons que la définition de $\rho_{\theta^\psi_{n,V}( \eta_\delta(Q)\otimes \det^\tau)}$ dépend du choix
d'une donnée de Whittaker. On en fixe une non ramifiée.

\begin{lemme}\label{Epique} On suppose que le corps local $F$ est non archimédien et  que $\eta(Q)=\eta=1$. On a alors 
\[( \epsilon_1(n,\eta_\delta(Q),\tau),\epsilon_2(n,\eta_\delta(Q),\tau), \epsilon_3(n,\eta_\delta(Q),\tau) )= ((-1)^\tau,(-1)^\tau\epsilon_\delta(Q), \epsilon_\delta(Q)), \]
et dans ce cas, ni $\eta_\delta(Q)$ ni $\epsilon_\delta(Q)$ ne dépendent de $\delta$.
\end{lemme}

\begin{rmq}
Nous verrons plus loin que cette formule est aussi valide dans le cas où  $\eta(Q)$ et $\eta$
sont non ramifiés et non égaux (corollaire \ref{corepique}). La dernière assertion dit que la formule est indépendante du choix de la donnée de Whittaker.
\end{rmq}

\dem
Si $\eta=1$, $\eta_\delta(Q)=\Triv_{\Or(V,Q)}$ et si $\eta(Q)=1$, la formule (\ref{hasse}) montre que $\epsilon_\delta(Q)$ ne dépend pas de $\delta$.
Le résultat est un cas particulier du lemme suivant,  qui est conséquence de la filtration de Kudla et des 
relations entre modules de Jacquet et paquets d'Arthur étudiés dans \cite{elementaire}.

\begin{lemme} On suppose que le corps local $F$ est non archimédien.

(i) On suppose que $\eta(Q)=1$.  Alors 
\[\epsilon_3(n,\eta_\delta(Q),\tau)=\epsilon_1(n,\eta_\delta(Q),\tau)\epsilon_2(n,\eta_\delta(Q),\tau)=\epsilon_\delta(Q).\]

(ii) On suppose que $\eta=1$.  Alors 
\[\epsilon_1(n,\eta_\delta(Q),\tau)=\epsilon_2(n,\eta_\delta(Q),\tau)\epsilon_3(n,\eta_\delta(Q),\tau)=(-1)^\tau.\]
\end{lemme}

\

\begin{rmq} Dans le lemme \ref{Epique}, le triplet $ ((-1)^\tau,(-1)^\tau\epsilon_\delta(Q), \epsilon_\delta(Q))$ ne correspond 
pas à un caractère de $A(\psi)$ lorsque les coordonnées correspondant à deux blocs de $\psi$ identiques ne sont pas identiques.
Nous avons vu que $\psi$ a deux blocs identiques lorsque $m=1$ et $\eta(Q)=1$, ou bien lorsque $n=2m-1$ et $\eta=1$.
Dans le premier cas, la remarque \ref{blocsegaux2} nous dit que nécessairement $\epsilon_\delta(Q)=1$, et donc le problème ne se pose pas.
Dans le second cas, il y a un problème éventuel  lorsque $\eta=1$ et $\tau=1$, mais c'est justement le
cas considéré dans la remarque \ref{blocsegaux1} où on a alors $\theta^{\psi}_{2m-1,V}(\eta_\delta(Q) \otimes \det^\tau)=0$, de sorte que le problème disparaît aussi.
\end{rmq}

On énonce maintenant des  résultats analogues lorsque $F=\bbR$. 
Ici $V$ est donc un espace vectoriel réel de dimension $N=2m$ muni d'une forme quadratique
$Q$ de signature $(p,q)$,  $p+q=N=2m$. On fixe aussi un caractère quadratique $\eta$ de $\bbR^\times$.
\begin{thm} \label{propR}
On suppose $n\geq 2m-1$.  La représentation $\theta^\psi_{n,V}( \eta_\delta(p,q)\otimes \det^\tau)$, 
si elle n'est pas nulle, est dans le $A$-paquet de paramètre 
\[\psi=
(\eta\otimes R[1]) \oplus (\eta \eta(p,q)\otimes R[2m-1]) \oplus (\eta(p,q)  \otimes R[2(n-m)+1]) \]
et le caractère de $A(\psi)$ attaché à $\theta^\psi_{n,V}( \eta_\delta(p,q)\otimes \det^\tau)$ est alors 
\[ ((-1)^\tau,(-1)^\tau\epsilon_\delta(p,q), \epsilon_\delta(p,q)). \]
\end{thm}

La première assertion est contenue dans le théorème \ref{lemp}. 
Nous allons démontrer le reste par des arguments globaux. On fixe donc un  un corps de nombre $k$, un caractère 
 adélique $\bm \eta=(\eta_v)_v$, 
 un espace quadratique  $(\caV, \caQ)$ sur $k$ de dimension $N=2m$ que l'on suppose anisotrope en au moins une place 
 (une telle place est forcément réelle  si $N>4$) 
 et une famille $(\tau_v)_v$ d'éléments de $\{0,1\}$
avec un nombre fini et pair de $v$ tels que $\tau_v=1$.

\begin{prop}\label{formprod} On suppose que $n\geq  2m$. Alors pour  toute place $v$ de $k$,
  $\theta^{\psi_v}_{n,V_v}( \eta_{\delta}(\caQ_v)\otimes \det^{\tau_v})$ 
est non nulle et appartient au  paquet $\Pi(\psi_v)$
où $\psi_v$ est le paramètre d'Arthur local
\[\psi_v=
(\eta_v\otimes R[1]) \oplus (\eta_v \eta(\caQ_v)\otimes R[2m-1]) \oplus (\eta(\caQ_v)  \otimes R[2(n-m)+1]). \]

On fixe une donnée de Whittaker globale non ramifiée aux places finies.
Le caractère de $A(\psi_v)$ associé à $\theta^{\psi_v}_{n,V_v}( \eta_{\delta}(\caQ_v)\otimes \det^{\tau_v})$ est donné par un triplet de signes
$ (\epsilon_{1,v},\epsilon_{2,v},\epsilon_{3,v} )  $
On a alors pour tout $i\in \{1,2,3\}$ la formule de produit :
$\prod_v \epsilon_{i,v}=1$.

Le même résultat est vrai si $n=2m-1$ et si l'on suppose que pour toute place $v$,  $\theta^{\psi_v}_{2m-1,V_v}( \eta_{\delta}(\caQ_v)\otimes 
\det^{\tau_v}) $ 
est non nulle et que de plus, 
en une place $v'$,  $\theta^{\psi_v}_{n',V_{v'}}( \eta_{\delta}(\caQ_{v'})\otimes \det^{\tau_{v'}})=0$  si $n'< m$.
\end{prop}
\dem On sait déjà d'après le théorème \ref{lemp} que $\theta^{\psi_v}_{n,V_v}( \eta_{\delta}(\caQ_v)\otimes \det^{\tau_v})$ 
est dans le paquet $\Pi(\psi_v)$ si elle est non nulle  et aussi que tel est le cas si $n\geq 2m$ (propriété (4') de la section \ref{Howeloc}).
Comme on a supposé que $\caQ$  est anisotrope à au moins une place, le  caractère $\bm \eta_\delta(\bm Q)\otimes \det^{\bm \tau}$
 de $\Or(\bm V,\bm Q)$ est automorphe cuspidal, et $\theta^{\bm \psi}_{n,\caV}(\bm \eta_\delta(\bm Q)\otimes \det^{\bm \tau})$, lorsqu'elle 
 est non nulle,  est une représentation automorphe de composantes locales les  $\theta^{\psi_v}_{n,V_v}( \eta_{\delta}(\caQ_v)\otimes 
\det^{\tau_v}) $.

 Supposons $n\geq 2m$. On sait alors que  $\theta^{\bm \psi}_{n,\caV}(\bm \eta_\delta(\bm Q)\otimes \det^{\bm \tau})$ est non nulle
 car on est dans le  stable range, et de carré intégrable  (propriété (5') de \ref{Howeglob}). 
Arthur lui associe alors un paramètre global
dont  les localisés $\psi_v$ soient comme dans l'énoncé de la proposition.  La formule  de produit n'est alors que la formule de multiplicité d'Arthur. 
Sous l'hypothèse $n=2m-1$ on n'est pas dans le rang stable. Toutefois l'hypothèse de non annulation des 
$\theta^{\psi_v}_{2m-1,V_v}( \eta_{\delta}(\caQ_v)\otimes  \det^\tau_v) $ 
permet d'appliquer \cite{yam} (l'hypothèse sur les   fonctions $L$ est facile à vérifier) et l'on a $n_0 =
n_0(\bm \eta_\delta(\bm Q)\otimes \det^{\bm \tau}) \leq 2m-1$.
 L'hypothèse locale en $v'$ implique que $n_0\geq m$.  Ainsi  $\theta^{\bm \psi}_{2m-1,\caV}(\bm \eta_\delta(\bm Q)\otimes \det^{\bm \tau})$   est non nulle et 
  c'est donc une représentation automorphe   de carré intégrable
  d'après la propriété  (5') de \ref{Howeglob} et l'on peut conclure comme ci-dessus par la théorie d'Arthur.
 \qed

\underline{\sl Démonstration du théorème \ref{propR}}. On démontre d'abord le cas où la forme est anisotrope, {\sl i.e.} $pq=0$. 
On globalise la situation sur $\bbQ$. On note $v_0$ la place réelle.
On suppose que la donnée de Whittaker globale se localise en la donnée $\mathrm{Wh}_\delta$ en $v_0$
On fixe  un caractère adélique $\bm \eta=(\eta_v)_v$ tel que $\eta_{v_0}=\eta$ et un espace quadratique  $(\caV,\caQ)$ de dimension 
$N=2m$ sur $\bbQ$ vérifiant les propriétés suivantes :

a)  $\eta_{v_0}=\eta$, $(\caV_{v_0},\caQ_{v_0})$ est de signature $(p,q)$ (avec $pq=0$).

b) Il existe deux places finies $v_1$ et $v_2$ telle que $\eta(\caQ_{v_1})=\eta_{v_1}=1$ et  $\eta(\caQ_{v_2})\neq   \eta_{v_2}$ avec 
 $\eta(\caQ_{v_2})$ et $ \eta_{v_2}$ non triviaux.

c) L'invariant de Hasse  $\epsilon_\delta(\caQ_v)$   est $1$ en toute place $v$ différente de $v_0$ et $v_1$, et en $v_0$ et $v_1$, les invariants de Hasse coïncident
(c'est nécessaire par la formule du produit), c'est-à-dire $\epsilon_{\delta}(\caQ_{v_0})= \epsilon_{\delta}(\caQ_{v_1})$.

On fixe une famille $\tau_v$ avec $\tau_{v_0}=\tau_{v_1}=\tau$ et telle que pour pour tout $v\notin\{ v_0,v_1,v_2\}$, la première occurence 
de $n_{0,v}=n_0(\eta_{\delta}(\caQ_v)\otimes \det^{\tau_v})$
vérifie  $n_{0,v}\leq m$. Si le nombre de $v\notin\{ v_0,v_1,v_2\}$ pour lequel $\tau_v=1$
est impair, on pose $\tau_{v_2}=1$, et  $\tau_{v_2}=0$ sinon (de sorte que l'on ait bien un nombre pair de places $v$ telle que $\tau_v=1$).

On suppose tout d'abord que $n\geq 2m$. Dans ce cas, $\theta_{n,\caV_v}^{\psi_v}(\eta_\delta(\caQ_v)\otimes \det^{\tau_v})$ est 
non nul pour toute place $v$ et les triplets  $(\epsilon_1(n,\eta_{\delta}(\caQ_v),\tau_{v}) ,  \epsilon_2(n,\eta_{\delta}(\caQ_v),\tau_{v}),   \epsilon_3(n,\eta_{\delta}(\caQ_v),\tau_{v}) )$ sont définis.
 Pour $i=1,2,3$, on pose $E_i=\prod_{v\neq v_0,v_1} \epsilon_i(n,\eta_{\delta}(\caQ_v),\tau_{v}) $.
Les hypothèses  sur $v_2$  assurent d'après la proposition \ref{occurence} (i) que $\theta^{\psi_{v_2}}_{n',\caV_{v_2}}(\eta_{\delta}(\caQ_{v_2}),\tau_{v_2})= 0$
 si $n'<m$, et les hypothèses de  la proposition \ref{formprod}  (avec $v'=v_2$) sont satisfaites et l'on a donc    pour $i=1,2,3$: 
\[E_i=\epsilon_i(n,\eta_{\delta}(\caQ_{v_0}),\tau_{v_0}) \epsilon_i(n,\eta_{\delta}(\caQ_{v_1}),\tau_{v_1}) .\]
La formule à démontrer est, d'après les hypothèses sur $v_1$ et le lemme \ref{lemp} :   
\[ \epsilon_i(n,\eta_{\delta}(\caQ_{v_0}),\tau_{v_0})=  \epsilon_i(n,\eta_{\delta}(\caQ_{v_1}),\tau_{v_1}),\quad  i=1,2,3,\] et ceci est donc équivalent à 
$E_i=1$, $i=1,2,3$.

Pour montrer que $i=1,2,3$, $E_i=1$, on va se ramener au cas $m=1$ et $n=1$. Supposons d'abord $m>1$. 
On introduit un espace quadratique $(\caV',\caQ')$ sur $\bbQ$ de dimension $2m-2$,  un caractère 
adélique quadratique $\bm \eta'$  et une famille $\bm \tau'=(\tau'_v)_v$ d'élément de $\bbZ/2\bbZ$
  avec les propriétés suivantes :
$\eta(\caQ'_v)=\eta(\caQ_v)$, $\epsilon_\delta(\caQ'_v)=\epsilon_\delta(\caQ'_v)$, 
$\eta'_v=  \eta'_v$ et $\tau'_v = \tau_v$ pour toute place $v\notin\{v_0,v_1\}$.
En $v_0$, si  $m-1>1$, on garde une forme anisotrope, et si $m=1$, on remarque qu'en la place $v_2$, 
comme $\eta'(\caQ_{v_2})=\eta(\caQ_{v_2})\neq 1$, la forme est anisotrope, de sorte que dans tout les cas, la forme 
 $(\caV',\caQ')$ est anisotrope à au moins une place. Aux places $v_0$ et $v_1$, on impose juste 
 $\epsilon_{\delta}(\caQ_{v_0})=\epsilon_{\delta}(\caQ_{v_1})$ et $\tau_{v_0}=\tau_{v_1}$,
mais l'on peut choisir ces valeurs respectivement dans $\{\pm 1\}$ et $\{0,1\}$.

On sait que pour toute place $v\notin\{v_0,v_1\}$,  $\theta_{n,\caV_v}^{\psi_v}(\eta_\delta(\caQ_v)\otimes \det^{\tau_v})$ 
est dans le paquet de paramètre
$\psi=(\eta_v\boxtimes R[1])\oplus (\eta_v \eta(\caQ_v)\boxtimes R[2m-1])\oplus (\eta(\caQ_v)\boxtimes R[2(n-m)+1])$
et de plus, les résultats  de \cite{elementaire} nous disent que $\theta_{n,\caV_v}^{\psi_v}(\eta_\delta(\caQ_v)\otimes \det^{\tau_v})$ 
est dans une induite de la forme $ \eta_v  \eta(\caQ_v)\vert\, \vert ^{m-1} \times \rho $, 
où $\rho$ est une représentation irréductible de $\Sp(2(n-1),\bbQ_v)$  dans le paquet de paramètre
$\psi'=(\eta_v\boxtimes R[1])\oplus (\eta_v \eta(\caQ_v)\boxtimes R[2m-3])\oplus (\eta(\caQ_v)\boxtimes R[2(n-m)+1])$.
En utilisant la filtration de Kudla, on voit de plus que $\rho=\theta_{n-1,\caV'_v}^{\psi_v}(\eta_\delta(Q'_v)\otimes \det^{\tau'_v})$.
D'autre part, le caractère de $A(\psi')$ associé par Arthur à $\theta_{n-1,\caV'_v}^{\psi_v}(\eta_\delta(Q'_v)\otimes \det^{\tau'_v})$
est donné par le triplet $(\epsilon_1(n,\eta_{\delta}(\caQ_v),\tau_{v}) ,  \epsilon_2(n,\eta_{\delta}(\caQ_v),\tau_{v}),   \epsilon_3(n,\eta_{\delta}(\caQ_v),\tau_{v}) )$.
Ceci nous ramène par récurrence au cas où $m=1$ (et $n$ a été remplacé par $n-m+1$).

On suppose donc $m=1$, et si $n>1$, on fait un raisonnement analogue pour remplacer $n$ par $n-1$ : 
 pour toute place $v\notin\{v_0,v_1\}$,  $\theta_{n,\caV_v}^{\psi_v}(\eta_\delta(\caQ_v)\otimes \det^{\tau_v})$ 
est dans le paquet de paramètre
\[\psi=(\eta_v\boxtimes R[1])\oplus (\eta_v \eta(\caQ_v)\boxtimes R[1])\oplus (\eta(\caQ_v)\boxtimes R[2n-1])\]
de plus $\theta_{n,\caV_v}^{\psi_v}(\eta_\delta(\caQ_v)\otimes \det^{\tau_v})$  est dans une induite de la forme
$   \eta(\caQ_v)\vert\, \vert ^{n-1} \times \rho $, 
où $\rho$ est une représentation irréductible de $\Sp(2(n-1),\bbQ_v)$  dans le paquet de paramètre
\[\psi'=(\eta_v\boxtimes R[1])\oplus( \eta_v \eta(\caQ_v)\boxtimes R[1])\oplus (\eta(\caQ_v)\boxtimes R[2n-3])\]
De plus  $\rho=\theta_{n-1,\caV_v}^{\psi_v}(\eta_\delta(\caQ_v)\otimes \det^{\tau_v})$ et le 
caractère de $A(\psi')$ associé par Arthur à $\theta_{n-1,\caV_v}^{\psi_v}(\eta_\delta(\caQ_v)\otimes \det^{\tau_v})$
est encore $(\epsilon_1(n,\eta_{\delta}(\caQ_v),\tau_{v}) ,  \epsilon_2(n,\eta_{\delta}(\caQ_v),\tau_{v}),   \epsilon_3(n,\eta_{\delta}(\caQ_v),\tau_{v}) )$.
On se ramène ainsi au cas $m=n=1$.

 On va alors de  nouveau utiliser la formule de produit de la proposition \ref{formprod}
mais dans le sens inverse, c'est-à-dire que l'on va montrer que $E_i=1$ pour $i=1,2,3$
en montrant que $ \epsilon_i(1,\eta_{\delta}(\caQ_{v_0}),\tau_{v_0})=  \epsilon_i(1,\eta_{\delta}(\caQ_{v_1},\tau_{v_1}),\quad  i=1,2,3$, mais avec ici $n=m=1$.
Pour pouvoir utiliser cette   formule de produit, il faut de plus 
qu'aux places $v_0$ et $v_1$, on ait choisi $\tau_{v_0}$, $\tau_{v_1}$ tels que $\theta^{\psi_{v_0}}_{1,\caV_{v_0}} (\eta_{\delta}(\caQ_{v_0})\otimes \det^{\tau_{v_0}}  ) $ et 
 $\theta^{\psi_{v_1}}_{1,\caV_{v_1}} (\eta_{\delta}(\caQ_{v_1})\otimes \det^{\tau_{v_1}}  ) $   soient non nulles, ce qui impose $\tau_{v_1}=0$.
On a donc $\tau_{v_0}=0$ et il faut que  $\theta^{\psi_{v_0}}_{1,\caV_{v_0}} (\eta_{\delta}(\caQ_{v_0}) ) $ soit non nulle. Si $\eta(\caQ_{v_0})=1$, c'est-à-dire si 
la signature de $\caQ_{v_0}$ est $(1,1)$, c'est le cas, car tous les caractères de $\Or(1,1)$ sauf le déterminant qui est ici exclu ont une image dans $\SL(2,\bbR)$ 
par la correspondance de Howe. Si $\eta_{v_0}$ est trivial, 
$\theta^{\psi_{v_0}}_{1,\caV_{v_0}} (\eta_{\delta}(\caQ_{v_0})  ) $  est la représentation sphérique et le caractère d'Arthur associé est trivial, c'est-à-dire 
$ \epsilon_i(1,\eta_{\delta}(\caQ_{v_0}),0)=  1,\quad  i=1,2,3$, et d'après le lemme \ref{Epique}, on a aussi 
$ \epsilon_i(1,\eta_{\delta}(\caQ_{v_1}),0)=  1,\quad  i=1,2,3$ ce qui conclut ce cas. On regarde ensuite le cas où  
 $\eta_{v_0}$ est le caractère signe.
On a vu dans la section \ref{caracO} que si $\delta=+1$, $\theta^{\psi_{v_0}}_{1,\caV_{v_0}} (\eta_{\delta}(\caQ_{v_0}) ) $  contient le $\U(1)$-type $\vert .\vert $
et si $\delta=-1$ il contient le $\U(1)$-type $\vert .\vert^{-1} $. Dans le premier cas, on a donc la limite de séries discrètes
holomorphe, et dans le second, une limite de séries discrète anti-holomorphe.
Le paramètre d'Arthur est ici tempéré, c'est un paramètre de Langlands, et on connait la paramétrisation de ces paquets
par les caractères de $A(\psi)$ (qui rappelons-le, dépend du choix de $\delta$ via le modèle de Whittaker). 
On a fait en sorte qu'avec le choix du modèle de Whittaker
correspondant à $\delta$, on ait  le caractère trivial, c'est-à-dire à nouveau $ \epsilon_i(1,\eta_{\delta}(\caQ_{v_0}),0)=  1,\quad  i=1,2,3$.
D'après le lemme \ref{Epique}, on a aussi 
$ \epsilon_i(1,\eta_{\delta}(\caQ_{v_1}),0)=  1,\quad  i=1,2,3$  car   $\epsilon_{\delta}(\caQ_{v_1})=\epsilon_{\delta}(\caQ_{v_0})=1$   ce qui conclut ce cas.

Considérons maintenant le cas où  $\eta(\caQ_{v_0})=-1$, c'est-à-dire que la signature de $\caQ_{v_0}$ est $(0,2)$ ou $(2,0)$. 
Si $\eta_{v_0}$ est trivial, $\eta_{\delta}(\caQ_{v_0})$ est le caractère trivial et $\theta^{\psi_{v_0}}_{1,\caV_{v_0}} (\eta_{\delta}(\caQ_{v_0})) $ 
 est non nulle, c'est la limite de série discrète 
holomorphe si la signature de $\caQ_{v_0}$ est $(2,0)$, antiholomorphe si la signature est $(0,2)$.
L'invariant de Hasse $\epsilon_\delta(\caQ_{v_0})$ est égal à $1$ si $\delta=1$ et si la signature  de $\caQ_{v_0}$ est $(2,0)$
ou bien si $\delta=-1$ et si la signature  de $\caQ_{v_0}$ est $(0,2)$,  et cet invariant de Hasse vaut $-1$ si 
$\delta=-1$ et si la signature  de $\caQ_{v_0}$ est $(2,0)$
ou bien si $\delta=1$ et si la signature  de $\caQ_{v_0}$ est $(0,2)$. Le paquet d'Arthur est dans ce cas un paquet tempéré et l'on conclut comme ci-dessus.

Si $\eta_{v_0}=-1$, $\theta^{\psi_{v_0}}_{1,\caV_{v_0}} (\eta_{\delta}(\caQ_{v_0}) ) $  
n'existe que si $\eta_\delta(\caQ_{v_0})$ est le caractère trivial, c'est-à-dire si
$\delta=1$ et la signature  de $\caQ_{v_0}$ est $(2,0)$
ou bien si $\delta=-1$ et si la signature  de $\caQ_{v_0}$ est $(0,2)$. Comme dans le cas précédent, 
$\theta^{\psi_{v_0}}_{1,\caV_{v_0}} (\eta_{\delta}(\caQ_{v_0}) ) $  est donc la limite de série discrète 
holomorphe si la signature de $\caQ_{v_0}$ est $(2,0)$, antiholomorphe si la signature est $(0,2)$ et les invariants de Hasse sont les mêmes. On conclut donc
de la même façon.

Ceci termine la démonstration de la proposition dans le cas où $n\geq 2m$. Nous  en tirons le résultat suivant qui nous permettra de conclure dans le cas 
$n=2m-1$.

\begin{cor}  \label{corepique}
On se replace dans la situation du début de la section avec un corps local  $F$  non archimédien de caractéristique $0$,
 un espace quadratique $(V,Q)$ de dimension $2m$
et un caractère quadratique $\eta$ de $F^\times$. On suppose    que $\eta(Q)$ et $\eta$ sont non ramifiés et non égaux et $n\geq 2m$. On a alors 
$( \epsilon_1(n,\eta_\delta(Q),\tau),\epsilon_2(n,\eta_\delta(Q),\tau), \epsilon_3(n,\eta_\delta(Q),\tau) )= ((-1)^\tau,(-1)^\tau\epsilon_\delta(Q), \epsilon_\delta(Q))$. 
\end{cor}

\dem On commence par le cas où $n\geq 2m$ et on considère une situation globale comme dans la proposition 
\ref{formprod} avec un corps de nombres $k$ et  une place $v_1'$ avec $k_{v_1}=F$
et une place $v_0$ réelle où la forme quadratique est anisotrope.  La formule de la proposition  est démontrée pour $v_0$. On 
utilise la formule de produit ainsi que celle obtenue en changeant simultanément les valeurs de 
$\tau_{v_0}$ et $\tau_{v'_1}$. Comme 
\[ \epsilon_1(n,\eta_\delta(\caQ_{v_0}), \tau_{v_0}) \epsilon_1(n,\eta_\delta(\caQ_{v_0}), \tau_{v_0}+1)
=\epsilon_2(n,\eta_\delta(\caQ_{v_0}), \tau_{v_0}) \epsilon_2(n,\eta_\delta(\caQ_{v_0}), \tau_{v_0}+1)=-1, \] 
\[\epsilon_3(n,\eta_\delta(\caQ_{v_0}), \tau_{v_0}) \epsilon_3(n,\eta_\delta(\caQ_{v_0}), \tau_{v_0}+1)=1 ,\] 
on obtient les mêmes identités avec $v'_1$:
\[ \epsilon_1(n,\eta_\delta (\caQ_{v'_1}), \tau_{v'_1}) \epsilon_1 (n,\eta_\delta(\caQ_{v'_1}), \tau_{v'_1}+1)=
\epsilon_2(n,\eta_\delta(\caQ_{v'_1}), \tau_{v'_1}) \epsilon_2 (n,\eta_\delta(\caQ_{v'_1}), \tau_{v'_1}+1)=-1, \] 
\[\epsilon_3(n,\eta_\delta(\caQ_{v'_1}), \tau_{v'_1}) \epsilon_3 (n,\eta_\delta(\caQ_{v'_1}), \tau_{v'_1}+1)=1 \] 

On fait de même en changeant cette fois simultanément les valeurs de $\eta_\delta (\caQ_{v_0})$ et $\eta_\delta (\caQ_{v'_1})$
et l'on obtient :
\[ \epsilon_1(n,\eta_\delta (\caQ_{v'_1}), \tau_{v'_1}) \epsilon_1 (n,-\eta_\delta(\caQ_{v'_1}), \tau_{v'_1})=1\]
\[\epsilon_2(n,\eta_\delta(\caQ_{v'_1}), \tau_{v'_1}) \epsilon_2 (n,-\eta_\delta(\caQ_{v'_1}), \tau_{v'_1})= 
\epsilon_3(n,\eta_\delta(\caQ_{v'_1}), \tau_{v'_1}) \epsilon_3 (n,-\eta_\delta(\caQ_{v'_1}), \tau_{v'_1})=-1 \] 
De plus, dans le cas où $\tau_{v'_1}=0$ et $\epsilon_\delta(\caQ_{v'_1})=1$, on sait que le 
caractère associé est trivial (correspondance de Howe non ramifiée). Ceci donne la formule voulue.
On utilise alors les arguments de \cite{elementaire} sur les  modules de Jacquet et la filtration de Kudla pour la place finie $v'_1$ 
comme dans la démonstration de la proposition \ref{occurence}  
pour montrer que ces identités sont aussi valides lorsque $n=2m-1$.
\qed

\

Revenons à la démonstration de la proposition dans le cas où $n=2m-1$.  Si l'on cherche à faire la même démonstration
que dans le cas $n\geq 2m-1$, on rencontre un obstacle qui est la possible annulation des 
$\theta_{n,\caV_v}^{\psi_v}(\eta_\delta(\caQ_v)\otimes \det^{\tau_v})$. 
En $v=v_0$, cette image est non nulle par hypothèse, et l'on peut toujours s'arranger pour que ce soit le cas à toutes les  places autres que $v_1$. En effet
les hypothèses en $v_2$ montrent que c'est le cas pour cette place.  En $v_1$ l'image est nulle si $\tau_{v_1}=1$.
On considère alors une place $v'_1$ comme dans le corollaire ci-dessus, et 
l'on change simultanément les valeurs de $\tau_{v_1}$ et $\tau_{v'_1}$
en $\tau_{v_1}+1=0$ et $\tau_{v'_1}+1$, de sorte que maintenant 
$\theta_{n,\caV_{v_1}}^{\psi_{v_1}}(\eta_\delta(\caQ_{v_1}))$  est non nulle, et l'on utilise la formule de produit pour en déduire 
la formule du lemme pour la place $v_0$. Ceci termine la démonstration de cette formule dans le cas anisotrope.

Maintenant, on ne suppose plus que la forme est anisotrope. On globalise alors de la manière suivante : $k$ est une
extension quadratique de $\bbQ$  avec exactement deux places réelles, $v_0$ et $v_0'$. En $v_0$, on met la forme quadratique qui nous intéresse,
et en $v_0'$, on met une forme anisotrope.
On fixe ensuite trois places finie, $v_1$, $v'_1$ et $v_2$. En la place $v_2$, on met les mêmes hypothèses que précédemment.
En $v_1$ et $v'_1$, on suppose $\eta_{v_1}$, $\eta_{v'_1}$, $\eta(\caQ_{v_1})$, $\eta(\caQ_{v'_1})$
triviaux. On suppose que $\epsilon_{\delta}(\caQ_{v_0})=\epsilon_{\delta}(\caQ_{v_1})$,  
$\epsilon_{\delta}(\caQ_{v'_0})=\epsilon_{\delta}(\caQ_{v'_1}) $,  $\tau_{v_0}=\tau_{v_1}$ et  $\tau_{v'_0}=
\tau_{v'_1}$. On fixe $\tau_{v_2}$ pour que le nombre de places  $v$ avec $\tau_v=1$ soit pair. 
Aux autres places finies, on fait les mêmes hypothèses que précédemment.
L'hypothèse de compacité à au moins une place est assurée par $v'_0$. 
Les produits $\epsilon_i(n,\eta_{\delta}(\caQ_{v'_0}),\tau_{v'_0})  \epsilon_i (n,\eta_{\delta}(\caQ_{v'_1}),\tau_{v'_1})$, $ i=1,2,3$
valent $1$ d'après ce que l'on a montré ci-dessus. On a alors pour $i=1,2,3$ d'après les formules de produits
\[  \epsilon_i( n,\eta_{\delta}(\caQ_{v_0}),\tau_{v_0})  \epsilon_i (n,\eta_{\delta}(\caQ_{v_1}),\tau_{v_1})
=E_i:=\prod_{v\neq v_0,v'_0, v_1,v'_1} \epsilon_i(n,\eta_{\delta}(\caQ_{v}),\tau_{v}). \]
Il s'agit alors de montrer que $E_i=1$ pour $i=1,2,3$, et l'on raisonne comme ci-dessus.\qed

\section{Calcul de $\rho_\pi$ pour les modules unitaires de plus haut/bas poids}\label{rhopiR}

On spécialise le résultat de la proposition \ref{propR} au cas où la forme quadratique $Q$ sur l'espace réel $V$  est anisotope.  
Rappelons les notations : la dimension de $V$ est $N=2m$, la signature de $Q$ est $(p,q)$ avec donc ici  $p=0$ ou $q=0$.
 L'image du caractère trivial et  du déterminant de $\Or(0,2m)$ dans $\Sp(2n,\bbR)$ par la correspondance de Howe
sont  respectivement les modules unitaires de plus haut poids noté $\pi_n(m)$ et $\sigma_{n,m}$. 
 L'image du caractère trivial et  du déterminant de $\Or(2m,0)$ dans $\Sp(2n,\bbR)$ par la correspondance de Howe
sont  respectivement les  les modules unitaires de plus bas poids $\pi_n(m)^*$ et $\sigma_{n,m}^*$.

\begin{prop} Fixons sur $\Sp(2n,\bbR)$ la donnée de Whittaker $\mathrm{Wh}_{\delta}$, $\delta\in \{\pm 1\}$ et soit $\tau\in \{0,1\}$.
On suppose $n\geq 2m-1+\tau$.
 L'image 
$\pi$ du caractère  $\det^\tau$ de $\Or(2m,0)$ dans $\Sp(2n,\bbR)$ par la correspondance de Howe 
appartient aux paquets d'Arthur $\Pi(\psi)$ de paramètre $\psi$
\[\psi=  (  \sgn_{W_\bbR}^{\tau'} \boxtimes R[1]) \oplus ( \sgn_{W_\bbR}^{\tau'}  \sgn_{W_\bbR}^{m}
  \boxtimes R[2m-1])\oplus  \sgn_{W_\bbR}^{m} \boxtimes R[2(n-m)+1], \quad (\tau'\in \{0,1\}). \]
et le caractère $\rho_\pi$ de $A(\psi)$ associé est 
\[   ((-1)^{\tau +\tau'(\frac{1+\delta}{2} +m )  },  (-1)^{\tau +\tau'(\frac{1+\delta}{2} +m ) +\lfloor \delta m/2  \rfloor } ,   (-1)^{\lfloor \delta m/2  \rfloor } ). \]

On a un résultat analogue si l'on remplace $\Or(2m,0)$ par $\Or(0,2m)$, en remplaçant $\lfloor \delta m/2  \rfloor $ par $\lfloor- \delta m/2  \rfloor $.
\end{prop}

\dem Si $\tau'=0$, c'est exactement l'énoncé de de la proposition \ref{propR}. Si $\tau'=1$,  il faut déterminer  $\tau''\in \{0,1\}$ tel
 que $\sgn_\delta\otimes \det^{\tau''}=\det^{\tau}$ sur $\Or(2m,0)$. 
Or on a calculé dans la section \ref{sgnd} que $\sgn_\delta$ est égal à $\det_{\Or(2m,0)}^{m+1}$ si $\delta=1$ et 
$\det_{\Or(2m,0)}^{m}$ si $\delta=-1$, d'où le résultat.\qed

\begin{cor}\label{carunip} On fixe  la donnée  Whittaker $\mathrm{Wh}_\delta$ sur $\Sp(2n,\bbR)$. 
Pour le paramètre 
\[\psi=  (  \Triv_{W_\bbR} \boxtimes R[1]) \oplus ( \sgn_{W_\bbR}^{m}
  \boxtimes R[2m-1])\oplus  \sgn_{W_\bbR}^{m} \boxtimes R[2(n-m)+1]   \]
les caractères  $\rho_\pi$ de $A(\psi)$ pour  $\pi=\pi_n(m), \, \sigma_{n,m} , \, \pi_n(m)^*, \, \sigma_{n,m}^*$ sont respectivement 
 \[(1 ,  (-1)^{\lfloor - \delta m/2  \rfloor } ,   (-1)^{ \lfloor - \delta m/2  \rfloor } ), 
   (-1  ,  (-1)^{1 +\lfloor - \delta m/2  \rfloor } ,   (-1)^{\lfloor - \delta m/2  \rfloor } ) , \]
 \[(1 ,  (-1)^{\lfloor  \delta m/2  \rfloor } ,   (-1)^{ \lfloor  \delta m/2  \rfloor } ), 
  (-1  ,  (-1)^{1 +\lfloor  \delta m/2  \rfloor } ,   (-1)^{\lfloor \delta m/2  \rfloor } ). \]

Pour le paramètre 
 \[\psi=  (  \sgn_{W_\bbR} \boxtimes R[1]) \oplus ( \sgn_{W_\bbR}^{m+1}
  \boxtimes R[2m-1])\oplus  \sgn_{W_\bbR}^{m} \boxtimes R[2(n-m)+1]   \]
les caractères  $\rho_\pi$ de $A(\psi)$ pour  $\pi=\pi_n(m), \, \sigma_{n,m} , \, \pi_n(m)^*, \, \sigma_{n,m}^*$ sont respectivement 
\begin{align*}  &((-1)^{(\frac{1+\delta}{2} +m )  },  (-1)^{(\frac{1+\delta}{2} +m ) +\lfloor -\delta m/2  \rfloor } ,   (-1)^{\lfloor -\delta m/2  \rfloor } ),\\
  & ((-1)^{1+\frac{1+\delta}{2} +m   },  (-1)^{1 +\frac{1+\delta}{2} +m  +\lfloor -\delta m/2  \rfloor } ,   (-1)^{\lfloor- \delta m/2  \rfloor } ), \\
&   ((-1)^{(\frac{1+\delta}{2} +m )  },  (-1)^{(\frac{1+\delta}{2} +m ) +\lfloor \delta m/2  \rfloor } ,   (-1)^{\lfloor \delta m/2  \rfloor } ),\\
  & ((-1)^{1+\frac{1+\delta}{2} +m   },  (-1)^{1 +\frac{1+\delta}{2} +m  +\lfloor \delta m/2  \rfloor } ,   (-1)^{\lfloor \delta m/2  \rfloor } ). \end{align*}

\end{cor}

Passons maintenant au cas d'un paramètre d'Arthur non nécessairement unipotent. On utilise les résultats de 
\cite{MR3} qui nous ramènent au cas unipotent du corollaire précédent, et 
ce qui suit est donc juste un exercice de traduction, où la seule difficulté est de faire attention aux choix de données de Whittaker
 ({\sl cf.} remarque \ref{rmswh}).

On se place donc dans les hypothèses des théorèmes
\ref{mainthm} et \ref{complement}. On fixe $m\leq n$, et  
l'on suppose donc que l'on a un paramètre $\psi=\psi_u\oplus \psi_d$
avec $\psi_d=\bigoplus_{i=1}^s (\delta_{t_i}\boxtimes R[a_i])$ et pour tout $i\in [1,s-1]$, $t_i\geq t_{i+1}$, tel que 
$\pi_n(m)\in \Pi(\psi)$. On est dans l'une des situations suivantes : 

(1) $\psi_u$ est irréductible.

(2) $\psi_u$ contient $\sgn_{W_\bbR}^m\boxtimes R[2(n-m)+1]$ et n'est pas irréductible.

(3) $\psi_u$ contient $\sgn_{W_\bbR}^{m-1}\boxtimes R[2(n-m)+3]$ et n'est pas irréductible.

On identifie le caractère $\rho_{\pi_n(m)^*}$  et de  $\rho_{\pi_n(m)}$  de $A(\psi)$ à une application de l'ensemble des blocs de $\psi$
dans $\{\pm 1\}$ telle que les valeurs prises sur des blocs égaux soient égales, et le produit de toutes ces valeurs est $1$.
Pour les applications à l'endoscopie et aux formules de multiplicité d'Arthur, on a besoin de connaître
les valeurs prises sur les blocs discrets  $\delta_{t_i}\boxtimes R[a_i]$, que l'on va noter 
$\epsilon( \delta_{t_i}\boxtimes R[a_i])$  et pour les blocs unipotents, lorsqu'il y en a trois, les 
produits des valeurs sur chaque couple de deux blocs. Notons $\epsilon_3$
la valeur de  $\rho_{\pi_n(m)^*}$ sur le bloc $\sgn_{W_\bbR}^{m-1}\boxtimes R[2(n-m)+3]$ du cas (3)
ou sur le bloc $\sgn_{W_\bbR}^m\boxtimes R[2(n-m)+1]$ du cas (2). Dans ces deux cas, 
notons $\eta_1\boxtimes R[1]$ et $\eta_2\boxtimes R[2a-1]$ (avec nécessairement $2a-1\leq 2(n-m)+1$ les deux autres blocs unipotents de $\psi$
et $\epsilon_1$ et $\epsilon_2$ les  valeurs de  $\rho_{\pi_n(m)}$ sur ces blocs, respectivement).
Pour tout $i\in [1,s]$, posons $a_{<i}= \sum_{j<i} a_j$, et $\delta_i=\delta(-1)^{a_{<i}}$, $\delta'=\delta (-1)^{\sum_{i=1}^s a_i}$.

\begin{prop} \label{Apsipinm}Fixons sur $\Sp(2n,\bbR)$ la donnée de Whittaker $\mathrm{Wh}_{\delta}$, $\delta\in \{\pm 1\}$.
Le caractère $\rho_{\pi_n(m)^*}$ de $A(\psi)$ est donnée par les formules suivantes.

Pour tout $i\in [1,s]$, $\epsilon( \delta_{t_i}\boxtimes R[a_i])=(-1)^{\lfloor  \delta_i  \frac{a_i}{2}  \rfloor } $.
  
Dans les cas (2) et (3), $\epsilon_1\epsilon_2=(-1)^{\lfloor \delta' \frac{a}{2}  \rfloor } $, 

Dans le cas (2),  $\epsilon_2\epsilon_3=1$ si $\eta_2=\sgn_{W_\bbR}^m$ et $\delta'(-1)^{a+1}$ si $\eta_2=\sgn_{W_\bbR}^{m+1}$.

Dans le cas (3),  $\epsilon_2\epsilon_3=-1$ si $\eta_2=\sgn_{W_\bbR}^{m-1}$ et $\delta'(-1)^{a}$ si $\eta_2=\sgn_{W_\bbR}^{m}$.
\end{prop}

\begin{rmq} $(i)$ Si $a=1$, on ne peut pas distinguer les indices $1$ et $2$  comme définis avant l'énoncé. Dans ce cas, $\epsilon_1\epsilon_2=\delta'$ 
et on vérifie que les formules
sont cohérentes par permutation des indices $1$ et $2$.

$(ii)$ On vérifie que $\delta'=(-1)^{m-a}\delta$ dans le cas (2) et $\delta'=(-1)^{m-a-1}\delta$ dans le cas (3). 

\end{rmq}

\dem La démonstration est une traduction immédiate du corollaire \ref{carunip} en tenant compte de la proposition \ref{complement}
 qui ramène le cas général au cas unipotent.\qed

On se place maintenant dans les hypothèses du théorème \ref{Thmsnk}.
On a un paramètre $\psi$ comme ci-dessus dans le cas (2), et on adopte les même notations. On démontre de la même façon le résultat suivant.

\begin{prop}\label{Apsisnm} Fixons sur $\Sp(2n,\bbR)$ la donnée de Whittaker $\mathrm{Wh}_{\delta}$, $\delta\in \{\pm 1\}$.
Le caractère $\rho_{\sigma_{n,m}^*}$ de $A(\psi)$ est donnée par les formules suivantes. 

Pour tout $i\in [1,s]$, $\epsilon( \delta_{t_i}\boxtimes R[a_i])=(-1)^{\lfloor  \delta_i  \frac{a_i}{2}  \rfloor } $, 
  $\epsilon_1\epsilon_2=(-1)^{\lfloor\delta' \frac{a}{2}  \rfloor } $, 
  $\epsilon_2\epsilon_3=-1$ si $\eta_2=\sgn_{W_\bbR}^m$ et $\delta(-1)^{m}$ si $\eta_2=\sgn_{W_\bbR}^{m+1}$.

\end{prop}

\section{Démonstration des théorèmes \ref{parLang} et \ref{laglands2}}\label{demparlang}
 On reprend les notations et hypothèses du théorème  \ref{parLang}
 On globalise la situation sur $k =\mathbb{Q}$ en prenant un espace quadratique 
 $(\caV,\caQ)$ de dimension $2m$ tel que à la place réelle, la forme quadratique soit de signature $(0,2m)$ et pour cette situation globale, 
 on adopte les notations de la section \ref{Howeglob}.
 
Pour tout $n\geq m$, on obtient par une récurrence immédiate de la propriété (7') de   \ref{Howeglob} que 
$\theta^{\psi}_{n,\caV}(\Triv_{\Or(\bm V,\bm Q)})$ se réalise 
 comme quotient de l'induite 
 \[ \eta(\caQ)  \vert .\vert^{-\left(n-m \right)}  \times  \cdots \times   \eta(\caQ)  \vert .\vert \times \theta^{\bm \psi}_{m,\caV}( \Triv_{\Or(\bm V,\bm Q)}) ,\]
 relative au  sous-groupe parabolique de facteur de Levi 
$  \underbrace{\GL_1\times   \cdots \times \GL_1}_{n-m}\times  \Sp_{2m} $ de $\Sp_{2n}$.
 A la place réelle, cela dit que $\pi_n(m)$ est quotient de Langlands de la représentation standard (\ref{parLang}). 
  Or cette représentation induite  a un unique quotient irréductible, le quotient de Langlands.

On se place maintenant dans les hypothèses du théorème \ref{laglands2}.
On remarque que dans le cas $n=2k$, $\sigma_{n,k}=\pi(m)$ avec $m=k+1$ et le résultat vient d'être démontré.
On globalise comme dans la démonstration  (avec ici $(\caV,\caQ)$ de dimension $2k$ et de signature $(0,2k)$ à la place réelle)
 qui précède et l'on  remarque que tout caractère automorphe cuspidal $\bm \sigma$  de 
 $\Or(\bm V,\bm Q)$ dont la restriction à $\Or(0,2k))$ est le déterminant a une image dans la correspondance thêta exactement dans le rang stable, 
 c'est-à-dire pour  $n\geq 2k$. La théorie de Rallis  (propriété (5') de   \ref{Howeglob}), utilisée exactement comme dans la preuve  précédente nous 
ramène alors à $n=2k$ et $\pi_n(k+1)$ ce qui permet de conclure. \qed

\section{Siegel-Weil local\label{siegelweil}}

Soit $k$  un corps de nombres  et l'on adopte les notations de la section \ref{Howeglob}
pour les objets attachés à ce corps. On a le résultat global
 suivant (Siegel-Weil). Soit $\bm \pi$ une représentation cuspidale irréductible de $\Sp(2n,\mathbb{A})$. 
Soit $\bm \eta$ un caractère quadratique de $\bbA^\times/k^\times $  et $s_0$ un demi-entier strictement positif. 
On suppose que la fonction $L$ partielle  $L^S(\bm \eta\times \bm \pi,s)$ a un pôle en $s=s_0+1$ et on suppose $s_0$ maximal avec cette propriété. 
Alors $\bm \pi$ s'obtient par la correspondance thêta à partir d'un groupe orthogonal de dimension $2(n-s_0)$, en particulier il faut que cet entier soit positif 
ou  nul. Cela se produit quand  la représentation  automorphe  $\bm \pi^{GL}$ de $\GL_{2n+1}$ obtenue  par  transfert endoscopique tordu
 (après stabilisation ({\sl cf}. \cite{Art13}) est une induite parabolique du produit tensoriel  du caractère 
 $\bm \eta\circ \det$ de $\GL_{2s_0+1}$ avec une  représentation convenable de $\GL_{2(n-s_0)}$, et que $2s_0+1$ est 
 supérieur ou égal aux dimensions des représentations de $\SL_2(\mathbb{C})$ intervenant dans le paramètre
 d'Arthur de $\bm \pi^{GL}$.

Nous allons énoncer un analogue local de ce résultat, sous
les hypothèses   suivantes :  on considère un paramètre d'Arthur  $\psi$ pour $\Sp(2n,\bbR)$ et 
  on suppose d'une part que  $a(\psi)=a(\psi_u)\geq n+1$ et   d'autre part que $\psi_u$ est de longueur au plus trois 
   et si la longueur est trois alors au moins une des sous-représentations  est de dimension un.
On écrit $\psi=\eta\boxtimes R[a(\psi_u)] \oplus \psi'$, ce qui définit $\psi'$.  On globalise 
 la situation en considérant un corps de nombre  $k$ avec au  moins trois places réelles,  $v_0$, $v_1$, $v_2$, $v_0$ étant celle qui nous intéresse. 
 On fixe aussi une place finie $w$. 
 On fixe $\pi_{v_0}$ une représentation dans $\Pi(\psi)$. En $v_i$, $i=1,2$ on considère un  paramètre local
 $\psi_{v_i}= \sgn_{W_\bbR}^{\frac{2n+1-a(\psi_u)}{2}} \boxtimes R[a(\psi_u)]\oplus \psi'_{v_i}$ tel que $\psi'_{v_i}$ ait le caractère infinitésimal de la représentation triviale
   et en $v_1$   on considère la représentation 
 $\pi_{v_1}$ de $\Sp(2n,\mathbb{R})$ image par la correspondance de Howe de la représentation triviale du groupe compact $\Or(2n+1-a(\psi_u))$ 
 tandis qu'en $v_2$, on définit $\pi_{v_2}$ comme l'image de la représentation $\det$ de ce même groupe.  
Ceci est bien défini car on a supposé que $a(\psi_u)\geq n+1$ d'où $a(\psi_u)> 2n+1-a(\psi_u)$,  ce qui est la condition pour que le caractère 
 non trivial du groupe compact ait une image dans la correspondance de Howe. On a vu que ces représentations, $\pi_{v_i}$, pour $i=1,2$,
  sont bien dans $\Pi(\psi_{v_i})$ ({\sl cf.} section \ref{sensdirect}). 
 
 On écrit le paramètre $\psi$ comme en (\ref{decApar}), c'est-à-dire 
 $$\psi=\bigoplus_{j=1}^s(\delta_{t_j} \boxtimes R[a_j])\oplus \bigoplus_{i=1}^r(\eta_{i} \boxtimes R[a'_i]).$$
 
 On fixe un paramètre d'Arthur global $\Psi=\{ (\rho_j,a_j)_{j=1,\ldots,s}, (\rho'_i,a'_i)_{i=1,\ldots,r} \}$ 
 où les $\rho_j$ (resp. les $\rho'_{i}$) sont des représentations cuspidales autoduales de $\GL_2$ (resp. de $\GL_1$)
 dont les localisés en $v_0$, $v_1$ et $v_2$ sont respectivement $\psi_{v_0}=\psi$, $\psi_{v_1}$ et $\psi_{v_2}$.
 Le paramètre global $\Psi'$ est obtenu à partir de $\Psi$ en enlevant 
 le  bloc $(\rho'_i,a_i)$  qui se localise en $v_0$ en $\eta\boxtimes R[a(\psi_u)] $.
 
A toutes les places finies $v\neq w$, le localisé $\psi_v$ de  $\Psi$ est non ramifié, et l'on prend  une représentation 
$\pi_v\in \Pi(\psi_v)$,  presque partout  non ramifiée.
A la place finie $w$  on impose que les $\rho_{j,w}$ et les $\rho'_{i,w}$ soient des représentations cuspidales toutes distinctes 
(respectivement de $\GL_2(k_w)$ et $\GL_1(k_w)$.

  Dans ces conditions le localisé $\psi_w$ de $\Psi$ en $w$      est un paramètre comme ceux étudiés en \cite{elementaire}
  et même plus précisément un paramètre dual d'un paramètre de séries discrètes, c'est-à-dire un morphisme trivial sur la 
 le facteur  $\SL_2(\mathbb{C})$ du groupe de Weil-Deligne de $k_w$. 
 Il est  montré dans \cite{elementaire} qu'il y a une bijection entre les représentations de 
  $\Pi(\psi_w)$ et les caractères du centralisateur de $\psi_w$. D'autre part, le centralisateur global s'envoie bijectivement dans le
   centralisateur local. 
On fixe    aussi une telle représentation    $\pi_w$ dans $\Pi(\Psi_w)$ mais de sorte que la formule de multiplicité d'Arthur soit satisfaite, ce qui est possible
d'après la remarque qui vient d'être faite. Il existe donc  une représentation automorphe de carré intégrable $\bm \pi$ de $\Sp(2n, \mathbb{A}_k)$ 
attaché au paramètre $\Psi$.

 On globalise le caractère $\eta$ en un caractère $\bm \eta$, de sorte qu'en $v_1$ et en $v_2$ il vaille $\sgn_{W_\bbR}^{\frac{2n+1-a(\psi_u)}{2}}$. 
  
  \begin{lemme} Si $a(\psi_u)=n+1=\dim(\psi')+1$, la représentation $\bm \pi$ est cuspidale
 et si  $a(\psi_u)>n+1$, il existe un entier $b\leq (a(\psi_u)-\dim(\psi')-1)/2$ 
(éventuellement $b=0$) et une représentation automorphe cuspidale $\bm \sigma$ de $\Sp(2\dim(\psi')+2b,\mathbb{A})$ de paramètre d'Arthur
égal à  l'union de $\Psi'$ et de  $(\bm \eta,\dim(\psi')+1+2b)$,  tels que $\bm \pi$ se réalise exactement 
 dans les résidus des séries d'Eisenstein :
$$
\left(\prod_{j\in [1,a_0]} \left(s_j-\frac{a(\psi_u)+1}{2}+j\right)\, 
E(\eta_\mathbb{A}|\, |^{s_{1}}\times \cdots \times \eta_\mathbb{A}|\; |^{s_{a_0}}\times \bm \sigma) \right)_{s_{a_0}
=\frac{a(\psi_u)+1}{2}-a_0, \cdots, s_{1}=\frac{a(\psi_u)-1}{2}},
$$
Ici $a_0=\frac{a(\psi_u)-\dim(\psi')-1}{2}-b$ et les $s_j$ sont des nombres complexes et   les évaluations se font dans l'ordre écrit de la gauche vers la droite. 

\end{lemme}

\dem Il est vraisemblable que $b=0$ mais on ne le démontre pas. 
On sait d'après \cite{manuscripta} que si $\pi$ n'est pas cuspidale il existe $(\rho,a)$ dans le paramètre d'Arthur de $\bm \pi$ tel que 
$a>1$ et une représentation automorphe irréductible $\bm \pi'$, ici nécessairement de carré intégrable, de $\Sp(2n-2d_\rho,\mathbb{A})$ 
dont le paramètre d'Arthur est obtenu à partir de celui de $\bm \pi$ en  remplaçant $(\rho,a)$ par $(\rho,a-2)$, 
et  tel que $\bm \pi$ se réalise comme
 quotient de l'espace des  résidus des séries d'Eisenstein 
 $$\left(  \left(s-\frac{a-1}{2} \right)\,  E(\rho |\, |^s\times \bm \pi')  \right)_{s=\frac{a-1}{2}}.$$
  Ainsi en toute place $v$ la composante locale $\pi_v$ est un quotient de $\rho_v| \; |^{(a-1)/2}\times \pi'_v$. On considère la place $v_1$. 
  Par dualité hermitienne et en utilisant l'unitarité, $\pi_{v_1}$ est un sous-module de l'induite $\rho_{v_1}|\; |^{-(a-1)/2}\times \pi'_{v_1}$.
   On réalise $\pi_{v_1}$ comme quotient de Langlands pour l'induite $\nu\times \pi_0$ où $\nu$ est le caractère de $\GL(\frac {a(\psi_u)-1}{2},\mathbb{R})$
    quotient de la série principale $\eta_{v_1}|\; |^{\frac {a(\psi_u)-1}{2}} \times \cdots \times \eta_{v_1}|\; |$ et 
    $\pi_0$ est une représentation tempérée de $\Sp(\dim (\psi'),\mathbb{R})$. On inclut 
  $\rho_{v_1}|\; |^{-(a-1)/2}$   dans une série principale de $\GL(2,\mathbb{R})$ 
 dont les exposants sont soit $-\frac{a(\psi_u)-1}{2}$ si 
 $(\rho,a)=(\eta,a(\psi_u))$ 
  et  sinon dans $[-\dim(\psi')/2+1,\dim(\psi')/2-1]$ dont soit le premier soit le deuxième (au moins) est strictement négatif. 
On envoie aussi $\pi'_{v_1}$ dans une série principale et on utilise la réciprocité de Frobenius. Ainsi parmi les exposants caractères 
(au sens de Hecht et Schmid) provenant de l'induite de $\nu\times \pi_0$ on doit trouver les exposants de 
 $\rho_{v_1}|\; |^{-(a-1)/2}$.  
   On a la formule de \cite{HS} 8.24 (a) pour calculer ces exposants et on constate que la seule possibilité est que 
$\rho_{v_1} $    soit le caractère $\eta_{v_1}$ avec comme exposant
 $-(a(\psi_u)-1)/2$. Cela force $(\rho,a)$ à coïncider avec le bloc que l'on a retiré de $\Psi$ pour définir $\Psi'$. 
On peut alors remplacer $n$ par $n-1$ et $a(\psi_u)$ par $a(\psi_u)-2$, on vérifie aisément que $\pi'_{v_1}$ est encore l'image
 de la représentation triviale du groupe compact $\Or(\dim(\psi'))$ et en $v_2$ c'est l'image du déterminant. 
 La procédure s'arrête si on arrive sur une représentation cuspidale et le lemme est alors facile ou si on arrive à 
 $a(\psi_u)=n+1=\dim(\psi')+1$ avec $n$ pair. On peut encore faire la procédure ci-dessus, si la représentation n'est pas cuspidale en utilisant 
 $v_1$ on voit que le seul exposant possible est $-\dim(\psi')/2$ mais il est trop négatif pour $\pi_{v_2}$.
  D'où le fait que la représentation est cuspidale si $a(\psi_u)=\dim(\psi')+1$. \qed

\begin{cor} \label{SW} Soit $\psi$ un $A$ paramètre pour $\Sp(2n,\bbR)$. On suppose que  $a(\psi)=a(\psi_u)\geq n+1$, 
que $\psi_u$ est de longueur au plus trois  et si $\psi_u$ est de longueur trois alors une des sous-représentations (au moins) est de dimension un.
Alors les représentations intervenant dans $\Pi(\psi)$ sont obtenues comme image par la correspondance de Howe à partir d'une représentation 
d'un groupe orthogonal $\Or(p,q)$ avec $p+q=\dim (\psi')$  contenue dans le paquet  $\Pi(\psi')$.
 Toute représentation de  $\Pi(\psi)$ a  multiplicité un.
\end{cor}

\dem
Pour la première assertion, on revient à la construction globale en reprenant les notations du lemme précédent. 
On applique la formule de Siegel-Weil de Kudla-Rallis (cf. \cite{KR94})     parce qu'étant donnés les paramètres de 
$\bm \sigma$ il est facile de voir que la fonction $L$ partielle de $\bm \eta\times \bm \sigma$ a un pôle en $s=\frac{a(\psi_u)+1}{2}-b$.
 Donc $\bm \sigma$ est dans l'image par série thêta d'une forme automorphe cuspidale d'un groupe orthogonal sur une forme de discriminant
  $\eta$ et de dimension $\dim (\psi')$. Cela reste alors vrai pour $\bm \pi$ grâce à la théorie des tours de Witt de Rallis (\cite{Rallis}). 
  Cela force l'analogue local et donc la première assertion.

Pour les multiplicités, on note $\bm \sigma'$ la représentation automorphe cuspidale de $\Or(\dim(\psi'),\mathbb{A})$
 (ici la forme du groupe orthogonal dépend de $\bm \pi$) dont l'image par série thêta est $\bm \pi$. La méthode de \cite{pourhowe},
  permet de conclure que la multiplicité de $\bm \pi$ dans l'espace des formes automorphes de carré intégrable de $\Sp(2n,\mathbb{A})$   est un
  si $\bm \sigma'$ intervient avec multiplicité un dans l'espace des formes automorphes de carré intégrable pour le groupe orthogonal. 
  Bien sûr pour  avoir multiplicité un, il suffit que la restriction de $\bm \sigma'$ au groupe spécial orthogonal intervienne avec multiplicité un. 
  On sait que les multiplicités globales se calculent en multipliant les multiplicités locales par un facteur global calculé par \cite{Art13}, 
  du moins cela calcule la multiplicité de la représentation des fonctions invariantes sous l'action d'un automorphisme extérieur
   ({\sl cf.} \cite{Art13} et \cite{taibi}).
  Les multiplicités locales sont un dès que l'on suppose qu'en les places archimédiennes le caractère infinitésimal de $\psi'$ est régulier, ce qui est loisible.
   Le facteur global vaut un si $\psi'$ contient au moins une représentation de dimension un de $W_\mathbb{R}\times \SL_2(\mathbb{C})$
    puisque l'on globalise par un caractère automorphe. Ceci règle le cas où $\psi_u$ a trois sous-représentations irréductibles.
     Si $\psi'$ n'a pas cette propriété on a le facteur deux, on ne peut pas y échapper. Mais les formules de multiplicités d'Arthur montrent que
      si en une place, la représentation du groupe spécial orthogonal n'est pas isomorphe à sa conjuguée sous le groupe orthogonal,
      alors on a quand même multiplicité un globale. Pour se ramener à ce cas là, on ajoute simplement dans la globalisation 
      une quatrième place réelle, $v_3$ où on s'arrange pour que le groupe orthogonal soit le groupe compact
      $\Or(\dim(\psi'_{v_3}/2))$  et où $\psi'_{v_3}$ a un caractère infinitésimal
       très régulier (c'est-à-dire ne contenant pas $0$). On n'a pas le choix de la représentation en cette place, c'est la représentation 
       de dimension finie avec le bon caractère infinitésimal mais il y a deux représentations pour $\SO$ non isomorphe conjuguées sous le 
       groupe orthogonal : celle qui correspond au caractère infinitésimal: $a_1>a_2\cdots > a_{\dim(\psi')/2}>0$ et celle où le dernier coefficient 
        $ a_{\dim(\psi')/2}$ est remplacé par son opposé.
   \qed

\bibliographystyle{smfalpha}

\bibliography{MR6}

\providecommand{\bysame}{\leavevmode ---\ }
\providecommand{\og}{``}
\providecommand{\fg}{''}
\providecommand{\smfandname}{et}
\providecommand{\smfedsname}{\'eds.}
\providecommand{\smfedname}{\'ed.}
\providecommand{\smfmastersthesisname}{M\'emoire}
\providecommand{\smfphdthesisname}{Th\`ese}
\begin{thebibliography}{EHW83}

\bibitem[ABV92]{ABV}
{\scshape J.~Adams, D.~Barbasch {\normalfont \smfandname} D.~A. Vogan, Jr.} --
  \emph{The {L}anglands classification and irreducible characters for real
  reductive groups}, Progress in Mathematics, vol. 104, Birkh\"auser Boston,
  Inc., Boston, MA, 1992.

\bibitem[Ada87]{Ad87}
{\scshape J.~Adams} -- {\og Unitary highest weight modules\fg}, \emph{Adv. in
  Math.} \textbf{63} (1987), no.~2, p.~113--137.

\bibitem[AJ87]{AdJo}
{\scshape J.~Adams {\normalfont \smfandname} J.~F. Johnson} -- {\og Endoscopic
  groups and packets of nontempered representations\fg}, \emph{Compositio
  Math.} \textbf{64} (1987), no.~3, p.~271--309.

\bibitem[AMR]{AMR}
{\scshape N.~Arancibia, C.~M{\oe}glin {\normalfont \smfandname} D.~Renard} --
  {\og Paquets d'{A}rthur des groupes classiques et unitaires, cas
  cohomologique\fg}, prépublication, \url{https://arxiv.org/abs/1507.01432},
  à paraître aux Annales de la Faculté des Sciences de Toulouse.

\bibitem[Art13]{Art13}
{\scshape J.~Arthur} -- \emph{The endoscopic classification of
  representations}, American Mathematical Society Colloquium Publications,
  vol.~61, American Mathematical Society, Providence, RI, 2013, Orthogonal and
  symplectic groups.

\bibitem[Bar89]{BaCG}
{\scshape D.~Barbasch} -- {\og The unitary dual for complex classical {L}ie
  groups\fg}, \emph{Invent. Math.} \textbf{96} (1989), no.~1, p.~103--176.

\bibitem[CL]{CL}
{\scshape G.~Chenevier {\normalfont \smfandname} J.~Lannes} -- {\og Formes
  automorphes et voisins de {K}neser des réseaux de {N}iemeier\fg},
  prépublication, \url{https://arxiv.org/abs/:1409.7616}.

\bibitem[CR15]{CR}
{\scshape G.~Chenevier {\normalfont \smfandname} D.~Renard} -- {\og Level one
  algebraic cusp forms of classical groups of small rank\fg}, \emph{Mem. Amer.
  Math. Soc.} \textbf{237} (2015), no.~1121, p.~v+122.

\bibitem[EHW83]{EHW}
{\scshape T.~Enright, R.~Howe {\normalfont \smfandname} N.~Wallach} -- {\og A
  classification of unitary highest weight modules\fg}, Representation theory
  of reductive groups ({P}ark {C}ity, {U}tah, 1982), Progr. Math., vol.~40,
  Birkh\"auser Boston, Boston, MA, 1983, p.~97--143.

\bibitem[GT16]{GanTakeda}
{\scshape W.~T. Gan {\normalfont \smfandname} S.~Takeda} -- {\og A proof of the
  {H}owe duality conjecture\fg}, \emph{J. Amer. Math. Soc.} \textbf{29} (2016),
  no.~2, p.~473--493.

\bibitem[How89]{Ho}
{\scshape R.~Howe} -- {\og Transcending classical invariant theory\fg},
  \emph{J. Amer. Math. Soc.} \textbf{2} (1989), no.~3, p.~535--552.

\bibitem[HS83]{HS}
{\scshape H.~Hecht {\normalfont \smfandname} W.~a. Schmid} -- {\og Characters,
  asymptotics and {${\mathfrak n}$}-homology of {H}arish-{C}handra modules\fg},
  \emph{Acta Math.} \textbf{151} (1983), no.~1-2, p.~49--151.

\bibitem[Igu62]{Igu}
{\scshape J.-i. Igusa} -- {\og On {S}iegel modular forms of genus two\fg},
  \emph{Amer. J. Math.} \textbf{84} (1962), p.~175--200.

\bibitem[Jak83]{Jak}
{\scshape H.~P. Jakobsen} -- {\og Hermitian symmetric spaces and their unitary
  highest weight modules\fg}, \emph{J. Funct. Anal.} \textbf{52} (1983), no.~3,
  p.~385--412.

\bibitem[KR94]{KR94}
{\scshape S.~S. Kudla {\normalfont \smfandname} S.~Rallis} -- {\og A
  regularized {S}iegel-{W}eil formula: the first term identity\fg}, \emph{Ann.
  of Math. (2)} \textbf{140} (1994), no.~1, p.~1--80.

\bibitem[KV78]{KV}
{\scshape M.~Kashiwara {\normalfont \smfandname} M.~Vergne} -- {\og On the
  {S}egal-{S}hale-{W}eil representations and harmonic polynomials\fg},
  \emph{Invent. Math.} \textbf{44} (1978), no.~1, p.~1--47.

\bibitem[KV95]{KnVo}
{\scshape A.~W. Knapp {\normalfont \smfandname} D.~A. Vogan, Jr.} --
  \emph{Cohomological induction and unitary representations}, Princeton
  Mathematical Series, vol.~45, Princeton University Press, Princeton, NJ,
  1995.

\bibitem[Mat04]{Mat04}
{\scshape H.~Matumoto} -- {\og On the representations of {${\rm Sp}(p,q)$} and
  {${\rm SO}^*(2n)$} unitarily induced from derived functor modules\fg},
  \emph{Compos. Math.} \textbf{140} (2004), no.~4, p.~1059--1096.

\bibitem[M{\oe}g06]{elementaire}
{\scshape C.~M{\oe}glin} -- {\og Sur certains paquets d'{A}rthur et involution
  d'{A}ubert-{S}chneider-{S}tuhler g\'en\'eralis\'ee\fg}, \emph{Represent.
  Theory} \textbf{10} (2006), p.~86--129.

\bibitem[M{\oe}g08]{manuscripta}
\bysame , {\og Formes automorphes de carr\'e int\'egrable non cuspidales\fg},
  \emph{Manuscripta Math.} \textbf{127} (2008), no.~4, p.~411--467.

\bibitem[M{\oe}g11]{pourkudla}
\bysame , {\og Conjecture d'{A}dams pour la correspondance de {H}owe et
  filtration de {K}udla\fg}, Arithmetic geometry and automorphic forms, Adv.
  Lect. Math. (ALM), vol.~19, Int. Press, Somerville, MA, 2011, p.~445--503.

\bibitem[M{\oe}g17]{pourhowe}
\bysame , {\og Paquets d'{A}rthur sp\'eciaux unipotents aux places
  archim\'ediennes et correspondance de {H}owe\fg}, Representation theory,
  number theory, and invariant theory, Progr. Math., vol. 323,
  Birkh\"auser/Springer, Cham, 2017, p.~469--502.

\bibitem[MRa]{MR5}
{\scshape C.~M{\oe}glin {\normalfont \smfandname} D.~Renard} -- {\og Sur les
  paquets d'{A}rthur aux places réelles, translation\fg}, prépublication,
  \url{https://arxiv.org/abs/1704.05096}.

\bibitem[MRb]{MR3}
\bysame , {\og Sur les paquets d'{A}rthur des groupes classiques r{é}els\fg},
  \emph{Journal of the European Mathematical Society, {\em {à} para{î}tre }}.

\bibitem[MW95]{MWDec}
{\scshape C.~M{\oe}glin {\normalfont \smfandname} J.-L. Waldspurger} --
  \emph{Spectral decomposition and {E}isenstein series}, Cambridge Tracts in
  Mathematics, vol. 113, Cambridge University Press, Cambridge, 1995, Une
  paraphrase de l'\'Ecriture [A paraphrase of Scripture].

\bibitem[NOT01]{NOT}
{\scshape K.~Nishiyama, H.~Ochiai {\normalfont \smfandname} K.~Taniguchi} --
  {\og Bernstein degree and associated cycles of {H}arish-{C}handra
  modules---{H}ermitian symmetric case\fg}, \emph{Ast\'erisque} (2001),
  no.~273, p.~13--80, Nilpotent orbits, associated cycles and Whittaker models
  for highest weight representations.

\bibitem[Oda94]{Oda}
{\scshape T.~Oda} -- {\og An explicit integral representation of {W}hittaker
  functions on {${\rm Sp}(2;{\bf R})$} for the large discrete series
  representations\fg}, \emph{Tohoku Math. J. (2)} \textbf{46} (1994), no.~2,
  p.~261--279.

\bibitem[Oht91]{Oht}
{\scshape T.~Ohta} -- {\og The closures of nilpotent orbits in the classical
  symmetric pairs and their singularities\fg}, \emph{Tohoku Math. J. (2)}
  \textbf{43} (1991), no.~2, p.~161--211.

\bibitem[Ral84]{Rallis}
{\scshape S.~Rallis} -- {\og On the {H}owe duality conjecture\fg},
  \emph{Compositio Math.} \textbf{51} (1984), no.~3, p.~333--399.

\bibitem[She15]{Sh15}
{\scshape D.~Shelstad} -- {\og On elliptic factors in real endoscopic transfer
  {I}\fg}, Representations of reductive groups, Progr. Math., vol. 312,
  Birkh\"auser/Springer, Cham, 2015, p.~455--504.

\bibitem[SZ15]{SunZhu}
{\scshape B.~Sun {\normalfont \smfandname} C.-B. Zhu} -- {\og Conservation
  relations for local theta correspondence\fg}, \emph{J. Amer. Math. Soc.}
  \textbf{28} (2015), no.~4, p.~939--983.

\bibitem[Ta{\"\i}]{taibi}
{\scshape O.~Ta{\"\i}bi} -- {\og Arthur's multiplicity formula for certain
  inner forms of special orthogonal and symplectic groups\fg}, {à} para{î}tre
  dans Journal of the European Mathematical Society.

\bibitem[Ta{\"\i}17]{taibi2}
\bysame , {\og Dimensions of spaces of level one automorphic forms for split
  classical groups using the trace formula\fg}, \emph{Ann. Sci. \'Ec. Norm.
  Sup\'er. (4)} \textbf{50} (2017), no.~2, p.~269--344.

\bibitem[Tra01]{TrapA}
{\scshape P.~E. Trapa} -- {\og Annihilators and associated varieties of
  {$A_{\frak q}(\lambda)$} modules for {$\mathrm U(p,q)$}\fg}, \emph{Compositio
  Math.} \textbf{129} (2001), no.~1, p.~1--45.

\bibitem[Tra05]{Trap}
\bysame , {\og Richardson orbits for real classical groups\fg}, \emph{J.
  Algebra} \textbf{286} (2005), no.~2, p.~361--385.

\bibitem[Tsu84]{Tsu}
{\scshape R.~Tsushima} -- {\og An explicit dimension formula for the spaces of
  generalized automorphic forms with respect to {${\rm Sp}(2,\,{\bf Z})$}\fg},
  Automorphic forms of several variables ({K}atata, 1983), Progr. Math.,
  vol.~46, Birkh\"auser Boston, Boston, MA, 1984, p.~378--383.

\bibitem[Tsu87]{Tsuy}
{\scshape S.~Tsuyumine} -- {\og On the {S}iegel modular function field of
  degree three\fg}, \emph{Compositio Math.} \textbf{63} (1987), no.~1,
  p.~83--98.

\bibitem[Vog81]{Vgreen}
{\scshape D.~A. Vogan, Jr.} -- \emph{Representations of real reductive {L}ie
  groups}, Progress in Mathematics, vol.~15, Birkh\"auser, Boston, Mass., 1981.

\bibitem[Vog88]{VogDS}
\bysame , {\og Irreducibility of discrete series representations for semisimple
  symmetric spaces\fg}, Representations of {L}ie groups, {K}yoto, {H}iroshima,
  1986, Adv. Stud. Pure Math., vol.~14, Academic Press, Boston, MA, 1988,
  p.~191--221.

\bibitem[Vog91]{VogAss}
\bysame , {\og Associated varieties and unipotent representations\fg}, Harmonic
  analysis on reductive groups ({B}runswick, {ME}, 1989), Progr. Math., vol.
  101, Birkh\"auser Boston, Boston, MA, 1991, p.~315--388.

\bibitem[VZ84]{VZ}
{\scshape D.~A. Vogan, Jr. {\normalfont \smfandname} G.~J. Zuckerman} -- {\og
  Unitary representations with nonzero cohomology\fg}, \emph{Compositio Math.}
  \textbf{53} (1984), no.~1, p.~51--90.

\bibitem[Wal90]{Wald89}
{\scshape J.-L. Waldspurger} -- {\og D\'emonstration d'une conjecture de
  dualit\'e de {H}owe dans le cas {$p$}-adique, {$p\neq 2$}\fg}, Festschrift in
  honor of {I}. {I}. {P}iatetski-{S}hapiro on the occasion of his sixtieth
  birthday, {P}art {I} ({R}amat {A}viv, 1989), Israel Math. Conf. Proc.,
  vol.~2, Weizmann, Jerusalem, 1990, p.~267--324.

\bibitem[Wal03]{Wallach}
{\scshape N.~R. Wallach} -- {\og Generalized {W}hittaker vectors for
  holomorphic and quaternionic representations\fg}, \emph{Comment. Math. Helv.}
  \textbf{78} (2003), no.~2, p.~266--307.

\bibitem[Yam14]{yam}
{\scshape S.~Yamana} -- {\og L-functions and theta correspondence for classical
  groups\fg}, \emph{Invent. Math.} \textbf{196} (2014), no.~3, p.~651--732.

\bibitem[Zhu03]{Zhu}
{\scshape C.-B. Zhu} -- {\og Representations with scalar {$K$}-types and
  applications\fg}, \emph{Israel J. Math.} \textbf{135} (2003), p.~111--124.

\end{thebibliography}

  \end{document}